\documentclass{article}
\usepackage{amsfonts}
\usepackage[section]{placeins}

\def\Q{{\mathbb Q}}

\def\A{{\mathcal A}}
\def\d{{\bf d}}

\def\mod{\mbox{\ mod\ }}
\def\sqN{{\sqrt{\frac N 2}}}  

\def\vector#1#2#3{\left(
\begin{array}{l} 
#1\\ 
#2\\ 
#3 
\end{array}\right)}

\def\matrix#1#2#3#4#5#6#7#8#9{{
 \left( \begin{array}{ccc}
#1 & #2 & #3 \\
#4 & #5 & #6 \\
#7 & #8 & #9 
\end{array} \right)}}

\def\determinant#1#2#3#4#5#6#7#8#9{{
 \left\vert \begin{array}{ccc}
#1 & #2 & #3 \\
#4 & #5 & #6 \\
#7 & #8 & #9 
\end{array} \right\vert}}

\usepackage{pstricks}
\begin{document}

\psset{unit=3cm}
\def\ThreeTiling{%
\pspicture(1, 1)
\qline(0,0)(1,0)            
\qline(0,0)(0.5,0.866025404)    
\qline(0.5,0.866025404)(1,0)    
\qline(0,0)(0.5,0.288675135)     
\qline(0.50,0.866025404)(0.50,0.288675135)  
\qline(0.50,0.288675135)(1,0)     
\endpspicture}

\def\ThreeTilingA{%
\pspicture(2, 1)
\qline(0,0)(2,0)            
\qline(0,0)(0.5,0.866025404)    
\qline(0.5,0.866025404)(2,0)    
\qline(0,0)(1,0.577350269)     
\qline(1,0.577350269)(1,0)  
\endpspicture}

\def\FourTilingA{%
\psset{unit=2cm}
\pspicture(2,1.2)
\qline(0,0)(1.73205081,0)              
\qline(1.73205081,0)(1.73205081,1)    
\qline(0.866025404,0.5)( 0.866025404,0)   
\qline( 0.866025404,0.5)(1.73205081,0)   
\qline( 0,0)(1.73205081,1) 
\qline( 0.866025404,0.5)(1.73205081,0.5)  
\endpspicture
\psset{unit=3cm}
}

\def\FourTilingB{%
\psset{unit=2cm}
\pspicture(2,1.2)
\qline(0,0)(1.73205081,0)              
\qline(1.73205081,0)(1.73205081,1)    
\qline(0.866025404,0.5)( 0.866025404,0)   
\qline( 0.866025404,0.5)(1.73205081,0)   
\qline( 0,0)(1.73205081,1) 
\qline(1.29903811,0.25)(1.73205081,1)
\endpspicture
\psset{unit=3cm}
}

\def\FourTilingC{%
\psset{unit=2cm}
\pspicture(2,1.2)
\qline(0,0)(1.73205081,0)              
\qline(1.73205081,0)(1.73205081,1)    
\qline(0.866025404,0.5)( 0.866025404,0)   
\qline( 0.866025404,0.5)(1.73205081,0)   
\qline( 0,0)(1.73205081,1) 
\qline(1.73205081,0)(1.29903811,0.75)
\endpspicture
\psset{unit=3cm}
}

\def\FourTiling{%
\pspicture(1,1)
\qline(0,0)(1,0)            
\qline(0,0)(0.5,0.866025404)    
\qline(0.5,0.866025404)(1,0)    
\qline(0.25,0.433012602)(0.75,0.433012602)  
\qline(0.25,0.433012602)(0.5,0)  
\qline(0.5,0)(0.75,0.433012602)   
\endpspicture}

\def\FiveTiling{%
\psset{unit=2cm}
\pspicture(2.5,1.2)
\qline(0,0)(2.5,0)              
\qline(0,0)(2,1)    
\qline(2,1)(2.5,0)   
\qline(2,1)(2,0)     
\qline(1,0.5)(2,0.5)  
\qline(2,0.5)(1,0)     
\qline(1,0)(1,0.5)    
\psset{unit=3cm}
\endpspicture}

\def\FiveTilingB{%
\psset{unit=2cm}
\pspicture(2.5,1.2)
\qline(0,0)(2.5,0)              
\qline(0,0)(2,1)    
\qline(2,1)(2.5,0)   
\qline(2,1)(2,0)     
\qline(1,0.5)(2,0.5)  
\qline(1,0.5)(2,0)     
\qline(1,0)(1,0.5)    
\psset{unit=3cm}
\endpspicture}

\def\SixTiling{%
\pspicture(1, 1)
\qline(0,0)(1,0)            
\qline(0,0)(0.5,0.866025404)    
\qline(0.5,0.866025404)(1,0)    
\qline(0,0)(0.75, 0.433012702)     
\qline(0.50,0.866025404)(0.50,0)  
\qline(0.25, 0.433012702)(1,0)     
\endpspicture}

\def\SixTilingA{%
\psset{unit=2cm}
\pspicture(1, 1)
\qline(0,0)( 3.46410162,0)            
\qline(0,0)(1.73205081,1)        
\qline(1.73205081,1)( 3.46410162,0)    
\qline(1.73205081,1)(1.73205081,0)  
\qline( 0.866025404,0.5)(1.15470054,0)   
\qline(1.15470054,0)(1.73205081,1)
\qline(1.73205081,1)(2.30940108,0)      
\qline(2.30940108,0)(2.59807621,0.5)    
\endpspicture}

\def\EightTiling{%
\pspicture(1,1)
\qline(0,0)(1,0)            
\qline(0,0)(0.5,0.866025404)    
\qline(0.5,0.866025404)(1,0)    
\qline(0.25,0.433012602)(0.75,0.433012602)  
\qline(0.25,0.433012602)(0.25,0)  
\qline(0.75,0.433012602)(0.75,0)  
\qline(0.5,0.866025404)(0.5,0)   
\qline(0.25,0.433012602)(0.5,0)  
\qline(0.5,0)(0.75,0.433012602)   
\endpspicture}

\def\NineTiling{%
\pspicture(1,1)
\qline(0,0)(1,0)            
\qline(0,0)(0.5,0.866025404)    
\qline(0.5,0.866025404)(1,0)    
\qline(0.167777,0.288675135)(0.833333,0.288675135)  
\qline(0.333333,0.577350269)(0.6777777, 0.577350269) 
\qline(0.1677777,0.288675135)(0.333333,0)
\qline(0.333333,0.577350269)(0.6777777,0)
\qline(0.333333,0)(0.677777,0.577350269)
\qline(0.677777,0)(0.833333,0.288675135)
\endpspicture}

\def\NineTilingA{%
\psset{unit=1cm}
\pspicture(6,3.3)
\qline(0,0)(6,0)            
\qline(0,0)(6,3)            
\qline(6,0)(6,3)            
\qline(2,0)(2,1)         
\qline(4,0)(4,2)          
\qline(5,0)(5,2)
\qline(4,2)(6,2)        
\qline(2,1)(4,1)
\qline(2,0)(4,1)        
\qline(4,0)(5,2)
\qline(5,0)(6,2)
\endpspicture
\psset{unit=3cm}
}

\def\TwelveTiling{%
\pspicture(1,1)
\qline(0,0)(1,0)            
\qline(0,0)(0.5,0.866025404)    
\qline(0.5,0.866025404)(1,0)    
\qline(0.25,0.433012602)(0.75,0.433012602)  
\qline(0.25,0.433012602)(0.5,0)  
\qline(0.5,0)(0.75,0.433012602)   
\qline(0,0)(0.25,0.144337567)
\qline(0.25,0.144337567)(0.25,0.433012602)
\qline(0.25,0.144337567)(0.5,0)
\qline(0.5,0)(0.75,0.144337567)
\qline(0.75,0.144337567)(0.75,0.433012602)
\qline(0.75,0.144337567)(1,0)
\qline(0.5,0.288675035)(0.5,0)
\qline(0.5,0.288675035)(0.25,0.433012602)
\qline(0.5,0.288675035)(0.75,0.433012602)
\qline(0.5,0.577350169)(0.5,0.866025404)
\qline(0.5,0.577350169)(0.25,0.433012602)
\qline(0.5,0.577350169)(0.75,0.433012602)
\endpspicture}

\def\ThirteenTiling{%
\pspicture(2.9,1.2)
\qline(0,0)(2.39871747,0)    
\qline(0,0)(1.7325081,1)    
\qline(1.7325081,1)(2.39871747,0)  
\qline(1.7325081,0)(1.7325081,1)  
\qline(2.06538414,0)(2.06538414,0.5)   
\qline(0.577350269,0)(0.577350269,0.333333333)  
\qline(1.15470054,0)(1.15470054, 0.666666667)   
\qline(1.15470054,0.67777777)(1.7325081,0.67777777) 
\qline(0.577350269,0.33333333)(1.7325081,0.33333333) 
\qline(1.7325081,0.5)(2.06538414,0.5)   
\qline(1.7325081,0.666666667)(0.577350269,0)  
\qline(1.7325081,0.333333333)(1.15470054,0)  
\qline(1.7325081,0.5)(2.06538414,0)
\endpspicture}

\def\SixteenTiling{%
\pspicture(1,1)
\qline(0,0)(1,0)            
\qline(0,0)(0.5,0.866025404)    
\qline(0.5,0.866025404)(1,0)    
\qline(0.125,0.216506351)(0.875,0.216506351)  
\qline(0.25,0.433012702)(0.75,  0.433012702) 
\qline(0.375,0.649519052)(0.625,0.649519052 ) 
\qline(0.125,0.216506351)(0.25,0)
\qline(0.25,0.433012702)(0.5,0)
\qline(0.375,0.649519052)(0.75,0)
\qline(0.25,0)(0.625,0.649519052 )
\qline(0.5,0)(0.75,  0.433012702)
\qline(0.75,0)(0.875,0.216506351)
\endpspicture}


\def\TwentySevenTiling{%
\psset{unit=0.5cm}
\pspicture(10.3923048,9.4)
\qline(0,0)(10.3923048,0)             
\qline(0,0)(5.19615242,9)             
\qline(5.19615242,9)(10.3923048,0)    
\qline(0,0)(1.73205081,1)             
\qline(1.73205081,1)(3.46410162,0)
\qline(3.46410162,0)(5.19615242,1)
\qline(5.19615242,1)(6.92820323,0)
\qline(6.92820323,0)(8.66025404,1)
\qline(8.66025404,1)(10.3923048,0)
\qline(1.73205081,1)(8.66025404,1)  
\qline(1.73205081,1)(1.73205081,3)  
\qline(5.19615242,1)(5.19615242,3)
\qline(8.66025404,1)(8.66025404,3)
\qline(1.73205081,3)(3.46410162,4)  
\qline(3.46410162,4)(5.19615242,3)
\qline(5.19615242,3)(6.92820323,4)
\qline(6.92820323,4)(8.66025404,3)  
\qline(3.46410162,4)(3.46410162,6)  
\qline(3.46410162,6)(5.19615242,7)  
\qline(5.19615242,7)(6.92820323,6)  
\qline(6.92820323,6)(6.92820323,4)  
\qline(5.19615242,7)(5.19615242,9)  
\qline(1.73205081,1)(5.19615242,7)  
\qline(5.19615242,7)(8.66025404,1)
\qline(3.46410162,4)(6.92820323,4)
\qline(3.46410162,4)(5.19615242,1)
\qline(5.19615242,1)(6.92820323,4)
\qline(3.46410162,2)(1.73205081,1)  
\qline(3.46410162,2)(3.46410162,4)
\qline(3.46410162,2)(5.19615242,1)
\qline(6.92820323,2)(5.19615242,1)  
\qline(6.92820323,2)(6.92820323,4)
\qline(6.92820323,2)(8.66025404,1)
\qline(5.19615242,5)(3.46410162,4)  
\qline(5.19615242,5)(6.92820323,4)
\qline(5.19615242,5)(5.19615242,7)
\endpspicture}

\def\TwelveTilingA{
\pspicture(1,1)
\psset{unit=0.005cm} 
\qline(409.5,10)(609.25,125.3257)
\qline(609.25,125.3257)(675.8334,10)
\qline(675.8334,10)(409.5,10)
\qline(809,240.6514)(809,10)
\qline(809,10)(675.8334,10)
\qline(675.8334,10)(809,240.6514)
\qline(809,240.6514)(609.25,125.3257)
\qline(609.25,125.3257)(675.8334,10)
\qline(675.8334,10)(809,240.6514)
\qline(809,240.6514)(609.25,125.3257)
\qline(609.25,125.3257)(542.6667,240.6514)
\qline(542.6667,240.6514)(809,240.6514)
\qline(409.5,10)(409.5,240.6514)
\qline(409.5,240.6514)(542.6667,240.6514)
\qline(542.6667,240.6514)(409.5,10)
\qline(409.5,10)(609.25,125.3257)
\qline(609.25,125.3257)(542.6667,240.6514)
\qline(542.6667,240.6514)(409.5,10)
\qline(9.999969,10)(209.75,125.3257)
\qline(209.75,125.3257)(276.3333,10)
\qline(276.3333,10)(9.999969,10)
\qline(409.5,240.6514)(409.5,10)
\qline(409.5,10)(276.3333,10)
\qline(276.3333,10)(409.5,240.6514)
\qline(409.5,240.6514)(209.75,125.3257)
\qline(209.75,125.3257)(276.3333,10)
\qline(276.3333,10)(409.5,240.6514)
\qline(409.5,240.6514)(609.25,355.9771)
\qline(609.25,355.9771)(675.8334,240.6514)
\qline(675.8334,240.6514)(409.5,240.6514)
\qline(809,471.3029)(809,240.6514)
\qline(809,240.6514)(675.8334,240.6514)
\qline(675.8334,240.6514)(809,471.3029)
\qline(809,471.3029)(609.25,355.9771)
\qline(609.25,355.9771)(675.8334,240.6514)
\qline(675.8334,240.6514)(809,471.3029)
\endpspicture}

\def\TwelveTilingB{
\pspicture(1,1)
\psset{unit=0.005cm} 
\qline(276.3333,10)(476.0833,125.3257)
\qline(476.0833,125.3257)(542.6667,10)
\qline(542.6667,10)(276.3333,10)
\qline(675.8333,240.6514)(675.8333,10)
\qline(675.8333,10)(542.6667,10)
\qline(542.6667,10)(675.8333,240.6514)
\qline(675.8333,240.6514)(476.0833,125.3257)
\qline(476.0833,125.3257)(542.6667,10)
\qline(542.6667,10)(675.8333,240.6514)
\qline(9.999969,10)(276.3333,10)
\qline(276.3333,10)(209.75,125.3257)
\qline(209.75,125.3257)(9.999969,10)
\qline(409.5,240.6514)(276.3333,10)
\qline(276.3333,10)(209.75,125.3257)
\qline(209.75,125.3257)(409.5,240.6514)
\qline(409.5,240.6514)(276.3333,10)
\qline(276.3333,10)(476.0833,125.3257)
\qline(476.0833,125.3257)(409.5,240.6514)
\qline(409.5,240.6514)(675.8334,240.6514)
\qline(675.8334,240.6514)(609.2501,355.9771)
\qline(609.2501,355.9771)(409.5,240.6514)
\qline(809,471.3029)(675.8334,240.6514)
\qline(675.8334,240.6514)(609.2501,355.9771)
\qline(609.2501,355.9771)(809,471.3029)
\qline(675.8333,10)(809,10)
\qline(809,10)(675.8333,240.6514)
\qline(675.8333,240.6514)(675.8333,10)
\qline(809,240.6514)(809,10)
\qline(809,10)(675.8333,240.6514)
\qline(675.8333,240.6514)(809,240.6514)
\qline(675.8333,240.6514)(809,240.6514)
\qline(809,240.6514)(809,471.3029)
\qline(809,471.3029)(675.8333,240.6514)
\endpspicture}
\def\TwentySevenTilingA{
\pspicture(1.0,1.0)
\psset{unit=0.005cm} 
\qline(542.6667,10)(675.8334,86.88382)
\qline(675.8334,86.88382)(720.2222,10)
\qline(720.2222,10)(542.6667,10)
\qline(809,163.7676)(809,10)
\qline(809,10)(720.2222,10)
\qline(720.2222,10)(809,163.7676)
\qline(809,163.7676)(675.8334,86.88382)
\qline(675.8334,86.88382)(720.2222,10)
\qline(720.2222,10)(809,163.7676)
\qline(809,163.7676)(675.8334,86.88382)
\qline(675.8334,86.88382)(631.4445,163.7676)
\qline(631.4445,163.7676)(809,163.7676)
\qline(542.6667,10)(542.6667,163.7676)
\qline(542.6667,163.7676)(631.4445,163.7676)
\qline(631.4445,163.7676)(542.6667,10)
\qline(542.6667,10)(675.8334,86.88382)
\qline(675.8334,86.88382)(631.4445,163.7676)
\qline(631.4445,163.7676)(542.6667,10)
\qline(276.3333,10)(409.5,86.88382)
\qline(409.5,86.88382)(453.8889,10)
\qline(453.8889,10)(276.3333,10)
\qline(542.6667,163.7676)(542.6667,10)
\qline(542.6667,10)(453.8889,10)
\qline(453.8889,10)(542.6667,163.7676)
\qline(542.6667,163.7676)(409.5,86.88382)
\qline(409.5,86.88382)(453.8889,10)
\qline(453.8889,10)(542.6667,163.7676)
\qline(542.6667,163.7676)(409.5,86.88382)
\qline(409.5,86.88382)(365.1111,163.7676)
\qline(365.1111,163.7676)(542.6667,163.7676)
\qline(276.3333,10)(276.3333,163.7676)
\qline(276.3333,163.7676)(365.1111,163.7676)
\qline(365.1111,163.7676)(276.3333,10)
\qline(276.3333,10)(409.5,86.88382)
\qline(409.5,86.88382)(365.1111,163.7676)
\qline(365.1111,163.7676)(276.3333,10)
\qline(9.999992,10)(143.1667,86.88382)
\qline(143.1667,86.88382)(187.5555,10)
\qline(187.5555,10)(9.999992,10)
\qline(276.3333,163.7676)(276.3333,10)
\qline(276.3333,10)(187.5555,10)
\qline(187.5555,10)(276.3333,163.7676)
\qline(276.3333,163.7676)(143.1667,86.88382)
\qline(143.1667,86.88382)(187.5555,10)
\qline(187.5555,10)(276.3333,163.7676)
\qline(542.6667,163.7676)(675.8334,240.6514)
\qline(675.8334,240.6514)(720.2222,163.7676)
\qline(720.2222,163.7676)(542.6667,163.7676)
\qline(809,317.5352)(809,163.7676)
\qline(809,163.7676)(720.2222,163.7676)
\qline(720.2222,163.7676)(809,317.5352)
\qline(809,317.5352)(675.8334,240.6514)
\qline(675.8334,240.6514)(720.2222,163.7676)
\qline(720.2222,163.7676)(809,317.5352)
\qline(809,317.5352)(675.8334,240.6514)
\qline(675.8334,240.6514)(631.4445,317.5352)
\qline(631.4445,317.5352)(809,317.5352)
\qline(542.6667,163.7676)(542.6667,317.5353)
\qline(542.6667,317.5353)(631.4445,317.5352)
\qline(631.4445,317.5352)(542.6667,163.7676)
\qline(542.6667,163.7676)(675.8334,240.6514)
\qline(675.8334,240.6514)(631.4445,317.5352)
\qline(631.4445,317.5352)(542.6667,163.7676)
\qline(276.3333,163.7676)(409.5,240.6514)
\qline(409.5,240.6514)(453.8889,163.7676)
\qline(453.8889,163.7676)(276.3333,163.7676)
\qline(542.6667,317.5352)(542.6667,163.7676)
\qline(542.6667,163.7676)(453.8889,163.7676)
\qline(453.8889,163.7676)(542.6667,317.5352)
\qline(542.6667,317.5352)(409.5,240.6514)
\qline(409.5,240.6514)(453.8889,163.7676)
\qline(453.8889,163.7676)(542.6667,317.5352)
\qline(542.6667,317.5352)(675.8334,394.4191)
\qline(675.8334,394.4191)(720.2222,317.5352)
\qline(720.2222,317.5352)(542.6667,317.5352)
\qline(809,471.3029)(809,317.5352)
\qline(809,317.5352)(720.2222,317.5352)
\qline(720.2222,317.5352)(809,471.3029)
\qline(809,471.3029)(675.8334,394.4191)
\qline(675.8334,394.4191)(720.2222,317.5352)
\qline(720.2222,317.5352)(809,471.3029)
\endpspicture}
\def\TwentySevenTilingB{
\pspicture(1.0,1.0)
\psset{unit=0.005cm} 
\qline(453.8889,10)(587.0555,86.88382)
\qline(587.0555,86.88382)(631.4445,10)
\qline(631.4445,10)(453.8889,10)
\qline(720.2222,163.7676)(720.2222,10)
\qline(720.2222,10)(631.4445,10)
\qline(631.4445,10)(720.2222,163.7676)
\qline(720.2222,163.7676)(587.0555,86.88382)
\qline(587.0555,86.88382)(631.4445,10)
\qline(631.4445,10)(720.2222,163.7676)
\qline(720.2222,163.7676)(587.0555,86.88382)
\qline(587.0555,86.88382)(542.6667,163.7676)
\qline(542.6667,163.7676)(720.2222,163.7676)
\qline(453.8889,10)(453.8889,163.7676)
\qline(453.8889,163.7676)(542.6667,163.7676)
\qline(542.6667,163.7676)(453.8889,10)
\qline(453.8889,10)(587.0555,86.88382)
\qline(587.0555,86.88382)(542.6667,163.7676)
\qline(542.6667,163.7676)(453.8889,10)
\qline(187.5555,10)(320.7222,86.88382)
\qline(320.7222,86.88382)(365.1111,10)
\qline(365.1111,10)(187.5555,10)
\qline(453.8889,163.7676)(453.8889,10)
\qline(453.8889,10)(365.1111,10)
\qline(365.1111,10)(453.8889,163.7676)
\qline(453.8889,163.7676)(320.7222,86.88382)
\qline(320.7222,86.88382)(365.1111,10)
\qline(365.1111,10)(453.8889,163.7676)
\qline(453.8889,163.7676)(587.0555,240.6514)
\qline(587.0555,240.6514)(631.4445,163.7676)
\qline(631.4445,163.7676)(453.8889,163.7676)
\qline(720.2222,317.5352)(720.2222,163.7676)
\qline(720.2222,163.7676)(631.4445,163.7676)
\qline(631.4445,163.7676)(720.2222,317.5352)
\qline(720.2222,317.5352)(587.0555,240.6514)
\qline(587.0555,240.6514)(631.4445,163.7676)
\qline(631.4445,163.7676)(720.2222,317.5352)
\qline(9.999992,10)(187.5555,10)
\qline(187.5555,10)(143.1667,86.88382)
\qline(143.1667,86.88382)(9.999992,10)
\qline(276.3333,163.7676)(187.5555,10)
\qline(187.5555,10)(143.1667,86.88382)
\qline(143.1667,86.88382)(276.3333,163.7676)
\qline(276.3333,163.7676)(187.5555,10)
\qline(187.5555,10)(320.7222,86.88382)
\qline(320.7222,86.88382)(276.3333,163.7676)
\qline(276.3333,163.7676)(453.8889,163.7676)
\qline(453.8889,163.7676)(409.5,240.6515)
\qline(409.5,240.6515)(276.3333,163.7676)
\qline(542.6666,317.5353)(453.8889,163.7676)
\qline(453.8889,163.7676)(409.5,240.6515)
\qline(409.5,240.6515)(542.6666,317.5353)
\qline(542.6666,317.5353)(453.8889,163.7676)
\qline(453.8889,163.7676)(587.0555,240.6515)
\qline(587.0555,240.6515)(542.6666,317.5353)
\qline(542.6667,317.5352)(720.2222,317.5352)
\qline(720.2222,317.5352)(675.8333,394.4191)
\qline(675.8333,394.4191)(542.6667,317.5352)
\qline(809,471.3029)(720.2222,317.5352)
\qline(720.2222,317.5352)(675.8333,394.4191)
\qline(675.8333,394.4191)(809,471.3029)
\qline(720.2222,10)(809,10)
\qline(809,10)(809,163.7676)
\qline(809,163.7676)(720.2222,10)
\qline(720.2222,163.7676)(720.2222,10)
\qline(720.2222,10)(809,163.7676)
\qline(809,163.7676)(720.2222,163.7676)
\qline(720.2222,163.7676)(809,163.7676)
\qline(809,163.7676)(720.2222,317.5353)
\qline(720.2222,317.5353)(720.2222,163.7676)
\qline(809,317.5353)(809,163.7676)
\qline(809,163.7676)(720.2222,317.5353)
\qline(720.2222,317.5353)(809,317.5353)
\qline(720.2222,317.5352)(809,317.5352)
\qline(809,317.5352)(809,471.3029)
\qline(809,471.3029)(720.2222,317.5352)
\endpspicture}

\def\FortyEightTiling{
\pspicture(2.0,2.0)
\psset{unit=0.01cm} 
\qline(10,10)(148.75,10)
\qline(148.75,10)(79.375,50.05365)
\qline(79.375,50.05365)(10,10)
\qline(148.75,10)(287.5,10)
\qline(287.5,10)(218.125,50.05365)
\qline(218.125,50.05365)(148.75,10)
\qline(287.5,10)(426.25,10)
\qline(426.25,10)(356.875,50.05365)
\qline(356.875,50.05365)(287.5,10)
\qline(426.25,10)(565,10)
\qline(565,10)(495.625,50.05365)
\qline(495.625,50.05365)(426.25,10)
\qline(10,10)(79.375,50.05365)
\qline(79.375,50.05365)(79.375,130.161)
\qline(79.375,130.161)(10,10)
\qline(79.375,130.161)(148.75,170.2147)
\qline(148.75,170.2147)(148.75,250.3221)
\qline(148.75,250.3221)(79.375,130.161)
\qline(148.75,250.3221)(218.125,290.3757)
\qline(218.125,290.3757)(218.125,370.4831)
\qline(218.125,370.4831)(148.75,250.3221)
\qline(218.125,370.4831)(287.5,410.5368)
\qline(287.5,410.5368)(287.5,490.6441)
\qline(287.5,490.6441)(218.125,370.4831)
\qline(565,10)(495.625,50.05365)
\qline(495.625,50.05365)(495.625,130.161)
\qline(495.625,130.161)(565,10)
\qline(495.625,130.161)(426.25,170.2147)
\qline(426.25,170.2147)(426.25,250.3221)
\qline(426.25,250.3221)(495.625,130.161)
\qline(426.25,250.3221)(356.875,290.3757)
\qline(356.875,290.3757)(356.875,370.4831)
\qline(356.875,370.4831)(426.25,250.3221)
\qline(356.875,370.4831)(287.5,410.5368)
\qline(287.5,410.5368)(287.5,490.6441)
\qline(287.5,490.6441)(356.875,370.4831)
\qline(287.5,410.5368)(218.125,290.3757)
\qline(218.125,290.3757)(218.125,370.4831)
\qline(218.125,370.4831)(287.5,410.5368)
\qline(287.5,410.5368)(218.125,290.3757)
\qline(218.125,290.3757)(287.5,330.4294)
\qline(287.5,330.4294)(287.5,410.5368)
\qline(218.125,290.3757)(356.875,290.3757)
\qline(356.875,290.3757)(287.5,250.3221)
\qline(287.5,250.3221)(218.125,290.3757)
\qline(218.125,290.3757)(356.875,290.3757)
\qline(356.875,290.3757)(287.5,330.4294)
\qline(287.5,330.4294)(218.125,290.3757)
\qline(356.875,290.3757)(287.5,410.5368)
\qline(287.5,410.5368)(356.875,370.4831)
\qline(356.875,370.4831)(356.875,290.3757)
\qline(356.875,290.3757)(287.5,410.5368)
\qline(287.5,410.5368)(287.5,330.4294)
\qline(287.5,330.4294)(356.875,290.3757)
\qline(218.125,290.3757)(148.75,170.2147)
\qline(148.75,170.2147)(148.75,250.3221)
\qline(148.75,250.3221)(218.125,290.3757)
\qline(218.125,290.3757)(148.75,170.2147)
\qline(148.75,170.2147)(218.125,210.2684)
\qline(218.125,210.2684)(218.125,290.3757)
\qline(148.75,170.2147)(287.5,170.2147)
\qline(287.5,170.2147)(218.125,130.161)
\qline(218.125,130.161)(148.75,170.2147)
\qline(148.75,170.2147)(287.5,170.2147)
\qline(287.5,170.2147)(218.125,210.2684)
\qline(218.125,210.2684)(148.75,170.2147)
\qline(287.5,170.2147)(218.125,290.3757)
\qline(218.125,290.3757)(287.5,250.3221)
\qline(287.5,250.3221)(287.5,170.2147)
\qline(287.5,170.2147)(218.125,290.3757)
\qline(218.125,290.3757)(218.125,210.2684)
\qline(218.125,210.2684)(287.5,170.2147)
\qline(356.875,290.3757)(287.5,170.2147)
\qline(287.5,170.2147)(287.5,250.3221)
\qline(287.5,250.3221)(356.875,290.3757)
\qline(356.875,290.3757)(287.5,170.2147)
\qline(287.5,170.2147)(356.875,210.2684)
\qline(356.875,210.2684)(356.875,290.3757)
\qline(287.5,170.2147)(426.25,170.2147)
\qline(426.25,170.2147)(356.875,130.161)
\qline(356.875,130.161)(287.5,170.2147)
\qline(287.5,170.2147)(426.25,170.2147)
\qline(426.25,170.2147)(356.875,210.2684)
\qline(356.875,210.2684)(287.5,170.2147)
\qline(426.25,170.2147)(356.875,290.3757)
\qline(356.875,290.3757)(426.25,250.3221)
\qline(426.25,250.3221)(426.25,170.2147)
\qline(426.25,170.2147)(356.875,290.3757)
\qline(356.875,290.3757)(356.875,210.2684)
\qline(356.875,210.2684)(426.25,170.2147)
\qline(148.75,170.2147)(79.375,50.05371)
\qline(79.375,50.05371)(79.375,130.161)
\qline(79.375,130.161)(148.75,170.2147)
\qline(148.75,170.2147)(79.375,50.05371)
\qline(79.375,50.05371)(148.75,90.10736)
\qline(148.75,90.10736)(148.75,170.2147)
\qline(79.375,50.05371)(218.125,50.05371)
\qline(218.125,50.05371)(148.75,10)
\qline(148.75,10)(79.375,50.05371)
\qline(79.375,50.05371)(218.125,50.05371)
\qline(218.125,50.05371)(148.75,90.10736)
\qline(148.75,90.10736)(79.375,50.05371)
\qline(218.125,50.05371)(148.75,170.2147)
\qline(148.75,170.2147)(218.125,130.161)
\qline(218.125,130.161)(218.125,50.05371)
\qline(218.125,50.05371)(148.75,170.2147)
\qline(148.75,170.2147)(148.75,90.10736)
\qline(148.75,90.10736)(218.125,50.05371)
\qline(287.5,170.2147)(218.125,50.05371)
\qline(218.125,50.05371)(218.125,130.161)
\qline(218.125,130.161)(287.5,170.2147)
\qline(287.5,170.2147)(218.125,50.05371)
\qline(218.125,50.05371)(287.5,90.10736)
\qline(287.5,90.10736)(287.5,170.2147)
\qline(218.125,50.05371)(356.875,50.05371)
\qline(356.875,50.05371)(287.5,10)
\qline(287.5,10)(218.125,50.05371)
\qline(218.125,50.05371)(356.875,50.05371)
\qline(356.875,50.05371)(287.5,90.10736)
\qline(287.5,90.10736)(218.125,50.05371)
\qline(356.875,50.05371)(287.5,170.2147)
\qline(287.5,170.2147)(356.875,130.161)
\qline(356.875,130.161)(356.875,50.05371)
\qline(356.875,50.05371)(287.5,170.2147)
\qline(287.5,170.2147)(287.5,90.10736)
\qline(287.5,90.10736)(356.875,50.05371)
\qline(426.25,170.2147)(356.875,50.05371)
\qline(356.875,50.05371)(356.875,130.161)
\qline(356.875,130.161)(426.25,170.2147)
\qline(426.25,170.2147)(356.875,50.05371)
\qline(356.875,50.05371)(426.25,90.10736)
\qline(426.25,90.10736)(426.25,170.2147)
\qline(356.875,50.05371)(495.625,50.05371)
\qline(495.625,50.05371)(426.25,10)
\qline(426.25,10)(356.875,50.05371)
\qline(356.875,50.05371)(495.625,50.05371)
\qline(495.625,50.05371)(426.25,90.10736)
\qline(426.25,90.10736)(356.875,50.05371)
\qline(495.625,50.05371)(426.25,170.2147)
\qline(426.25,170.2147)(495.625,130.161)
\qline(495.625,130.161)(495.625,50.05371)
\qline(495.625,50.05371)(426.25,170.2147)
\qline(426.25,170.2147)(426.25,90.10736)
\qline(426.25,90.10736)(495.625,50.05371)
\endpspicture}

\def\OneHundredTwentyFiveTiling{
\pspicture(2.0,2.0)
\psset{unit=0.01cm} 
\qline(10,10)(121,10)
\qline(121,10)(65.5,42.04297)
\qline(65.5,42.04297)(10,10)
\qline(121,10)(232,10)
\qline(232,10)(176.5,42.04297)
\qline(176.5,42.04297)(121,10)
\qline(232,10)(343,10)
\qline(343,10)(287.5,42.04297)
\qline(287.5,42.04297)(232,10)
\qline(343,10)(454,10)
\qline(454,10)(398.5,42.04297)
\qline(398.5,42.04297)(343,10)
\qline(454,10)(565,10)
\qline(565,10)(509.5,42.04297)
\qline(509.5,42.04297)(454,10)
\qline(10,10)(65.5,42.04297)
\qline(65.5,42.04297)(65.5,106.1288)
\qline(65.5,106.1288)(10,10)
\qline(65.5,106.1288)(121,138.1718)
\qline(121,138.1718)(121,202.2576)
\qline(121,202.2576)(65.5,106.1288)
\qline(121,202.2577)(176.5,234.3006)
\qline(176.5,234.3006)(176.5,298.3865)
\qline(176.5,298.3865)(121,202.2577)
\qline(176.5,298.3865)(232,330.4294)
\qline(232,330.4294)(232,394.5153)
\qline(232,394.5153)(176.5,298.3865)
\qline(232,394.5153)(287.5,426.5582)
\qline(287.5,426.5582)(287.5,490.6441)
\qline(287.5,490.6441)(232,394.5153)
\qline(565,10)(509.5,42.04297)
\qline(509.5,42.04297)(509.5,106.1288)
\qline(509.5,106.1288)(565,10)
\qline(509.5,106.1288)(454,138.1718)
\qline(454,138.1718)(454,202.2576)
\qline(454,202.2576)(509.5,106.1288)
\qline(454,202.2577)(398.5,234.3006)
\qline(398.5,234.3006)(398.5,298.3865)
\qline(398.5,298.3865)(454,202.2577)
\qline(398.5,298.3865)(343,330.4294)
\qline(343,330.4294)(343,394.5153)
\qline(343,394.5153)(398.5,298.3865)
\qline(343,394.5153)(287.5,426.5582)
\qline(287.5,426.5582)(287.5,490.6441)
\qline(287.5,490.6441)(343,394.5153)
\qline(287.5,426.5583)(232,330.4294)
\qline(232,330.4294)(232,394.5153)
\qline(232,394.5153)(287.5,426.5583)
\qline(287.5,426.5583)(232,330.4294)
\qline(232,330.4294)(287.5,362.4724)
\qline(287.5,362.4724)(287.5,426.5583)
\qline(232,330.4294)(343,330.4294)
\qline(343,330.4294)(287.5,298.3865)
\qline(287.5,298.3865)(232,330.4294)
\qline(232,330.4294)(343,330.4294)
\qline(343,330.4294)(287.5,362.4724)
\qline(287.5,362.4724)(232,330.4294)
\qline(343,330.4294)(287.5,426.5583)
\qline(287.5,426.5583)(343,394.5153)
\qline(343,394.5153)(343,330.4294)
\qline(343,330.4294)(287.5,426.5583)
\qline(287.5,426.5583)(287.5,362.4724)
\qline(287.5,362.4724)(343,330.4294)
\qline(232,330.4294)(176.5,234.3006)
\qline(176.5,234.3006)(176.5,298.3865)
\qline(176.5,298.3865)(232,330.4294)
\qline(232,330.4294)(176.5,234.3006)
\qline(176.5,234.3006)(232,266.3435)
\qline(232,266.3435)(232,330.4294)
\qline(176.5,234.3006)(287.5,234.3006)
\qline(287.5,234.3006)(232,202.2577)
\qline(232,202.2577)(176.5,234.3006)
\qline(176.5,234.3006)(287.5,234.3006)
\qline(287.5,234.3006)(232,266.3435)
\qline(232,266.3435)(176.5,234.3006)
\qline(287.5,234.3006)(232,330.4294)
\qline(232,330.4294)(287.5,298.3865)
\qline(287.5,298.3865)(287.5,234.3006)
\qline(287.5,234.3006)(232,330.4294)
\qline(232,330.4294)(232,266.3435)
\qline(232,266.3435)(287.5,234.3006)
\qline(343,330.4294)(287.5,234.3006)
\qline(287.5,234.3006)(287.5,298.3865)
\qline(287.5,298.3865)(343,330.4294)
\qline(343,330.4294)(287.5,234.3006)
\qline(287.5,234.3006)(343,266.3435)
\qline(343,266.3435)(343,330.4294)
\qline(287.5,234.3006)(398.5,234.3006)
\qline(398.5,234.3006)(343,202.2577)
\qline(343,202.2577)(287.5,234.3006)
\qline(287.5,234.3006)(398.5,234.3006)
\qline(398.5,234.3006)(343,266.3435)
\qline(343,266.3435)(287.5,234.3006)
\qline(398.5,234.3006)(343,330.4294)
\qline(343,330.4294)(398.5,298.3865)
\qline(398.5,298.3865)(398.5,234.3006)
\qline(398.5,234.3006)(343,330.4294)
\qline(343,330.4294)(343,266.3435)
\qline(343,266.3435)(398.5,234.3006)
\qline(176.5,234.3006)(121,138.1718)
\qline(121,138.1718)(121,202.2577)
\qline(121,202.2577)(176.5,234.3006)
\qline(176.5,234.3006)(121,138.1718)
\qline(121,138.1718)(176.5,170.2147)
\qline(176.5,170.2147)(176.5,234.3006)
\qline(121,138.1718)(232,138.1718)
\qline(232,138.1718)(176.5,106.1288)
\qline(176.5,106.1288)(121,138.1718)
\qline(121,138.1718)(232,138.1718)
\qline(232,138.1718)(176.5,170.2147)
\qline(176.5,170.2147)(121,138.1718)
\qline(232,138.1718)(176.5,234.3006)
\qline(176.5,234.3006)(232,202.2577)
\qline(232,202.2577)(232,138.1718)
\qline(232,138.1718)(176.5,234.3006)
\qline(176.5,234.3006)(176.5,170.2147)
\qline(176.5,170.2147)(232,138.1718)
\qline(287.5,234.3006)(232,138.1718)
\qline(232,138.1718)(232,202.2577)
\qline(232,202.2577)(287.5,234.3006)
\qline(287.5,234.3006)(232,138.1718)
\qline(232,138.1718)(287.5,170.2147)
\qline(287.5,170.2147)(287.5,234.3006)
\qline(232,138.1718)(343,138.1718)
\qline(343,138.1718)(287.5,106.1288)
\qline(287.5,106.1288)(232,138.1718)
\qline(232,138.1718)(343,138.1718)
\qline(343,138.1718)(287.5,170.2147)
\qline(287.5,170.2147)(232,138.1718)
\qline(343,138.1718)(287.5,234.3006)
\qline(287.5,234.3006)(343,202.2577)
\qline(343,202.2577)(343,138.1718)
\qline(343,138.1718)(287.5,234.3006)
\qline(287.5,234.3006)(287.5,170.2147)
\qline(287.5,170.2147)(343,138.1718)
\qline(398.5,234.3006)(343,138.1718)
\qline(343,138.1718)(343,202.2577)
\qline(343,202.2577)(398.5,234.3006)
\qline(398.5,234.3006)(343,138.1718)
\qline(343,138.1718)(398.5,170.2147)
\qline(398.5,170.2147)(398.5,234.3006)
\qline(343,138.1718)(454,138.1718)
\qline(454,138.1718)(398.5,106.1288)
\qline(398.5,106.1288)(343,138.1718)
\qline(343,138.1718)(454,138.1718)
\qline(454,138.1718)(398.5,170.2147)
\qline(398.5,170.2147)(343,138.1718)
\qline(454,138.1718)(398.5,234.3006)
\qline(398.5,234.3006)(454,202.2577)
\qline(454,202.2577)(454,138.1718)
\qline(454,138.1718)(398.5,234.3006)
\qline(398.5,234.3006)(398.5,170.2147)
\qline(398.5,170.2147)(454,138.1718)
\qline(121,138.1718)(65.5,42.04291)
\qline(65.5,42.04291)(65.5,106.1288)
\qline(65.5,106.1288)(121,138.1718)
\qline(121,138.1718)(65.5,42.04291)
\qline(65.5,42.04291)(121,74.08588)
\qline(121,74.08588)(121,138.1718)
\qline(65.5,42.04291)(176.5,42.04291)
\qline(176.5,42.04291)(121,10)
\qline(121,10)(65.5,42.04291)
\qline(65.5,42.04291)(176.5,42.04291)
\qline(176.5,42.04291)(121,74.08588)
\qline(121,74.08588)(65.5,42.04291)
\qline(176.5,42.04291)(121,138.1718)
\qline(121,138.1718)(176.5,106.1288)
\qline(176.5,106.1288)(176.5,42.04291)
\qline(176.5,42.04291)(121,138.1718)
\qline(121,138.1718)(121,74.08588)
\qline(121,74.08588)(176.5,42.04291)
\qline(232,138.1718)(176.5,42.04291)
\qline(176.5,42.04291)(176.5,106.1288)
\qline(176.5,106.1288)(232,138.1718)
\qline(232,138.1718)(176.5,42.04291)
\qline(176.5,42.04291)(232,74.08588)
\qline(232,74.08588)(232,138.1718)
\qline(176.5,42.04291)(287.5,42.04291)
\qline(287.5,42.04291)(232,10)
\qline(232,10)(176.5,42.04291)
\qline(176.5,42.04291)(287.5,42.04291)
\qline(287.5,42.04291)(232,74.08588)
\qline(232,74.08588)(176.5,42.04291)
\qline(287.5,42.04291)(232,138.1718)
\qline(232,138.1718)(287.5,106.1288)
\qline(287.5,106.1288)(287.5,42.04291)
\qline(287.5,42.04291)(232,138.1718)
\qline(232,138.1718)(232,74.08588)
\qline(232,74.08588)(287.5,42.04291)
\qline(343,138.1718)(287.5,42.04291)
\qline(287.5,42.04291)(287.5,106.1288)
\qline(287.5,106.1288)(343,138.1718)
\qline(343,138.1718)(287.5,42.04291)
\qline(287.5,42.04291)(343,74.08588)
\qline(343,74.08588)(343,138.1718)
\qline(287.5,42.04291)(398.5,42.04291)
\qline(398.5,42.04291)(343,10)
\qline(343,10)(287.5,42.04291)
\qline(287.5,42.04291)(398.5,42.04291)
\qline(398.5,42.04291)(343,74.08588)
\qline(343,74.08588)(287.5,42.04291)
\qline(398.5,42.04291)(343,138.1718)
\qline(343,138.1718)(398.5,106.1288)
\qline(398.5,106.1288)(398.5,42.04291)
\qline(398.5,42.04291)(343,138.1718)
\qline(343,138.1718)(343,74.08588)
\qline(343,74.08588)(398.5,42.04291)
\qline(454,138.1718)(398.5,42.04291)
\qline(398.5,42.04291)(398.5,106.1288)
\qline(398.5,106.1288)(454,138.1718)
\qline(454,138.1718)(398.5,42.04291)
\qline(398.5,42.04291)(454,74.08588)
\qline(454,74.08588)(454,138.1718)
\qline(398.5,42.04291)(509.5,42.04291)
\qline(509.5,42.04291)(454,10)
\qline(454,10)(398.5,42.04291)
\qline(398.5,42.04291)(509.5,42.04291)
\qline(509.5,42.04291)(454,74.08588)
\qline(454,74.08588)(398.5,42.04291)
\qline(509.5,42.04291)(454,138.1718)
\qline(454,138.1718)(509.5,106.1288)
\qline(509.5,106.1288)(509.5,42.04291)
\qline(509.5,42.04291)(454,138.1718)
\qline(454,138.1718)(454,74.08588)
\qline(454,74.08588)(509.5,42.04291)
\endpspicture}

\def\TwentyEightTiling{
\vskip 5in
\pspicture(2.0,2.0)
\psset{unit=1cm}
\qline(0,0)(12,0)\qline(0,0)(8.50,13.56)\qline(12,0)(8.50,13.56)
\psgrid(0,0)(12,14)\pspolygon[fillcolor=lightgray,fillstyle=solid](8.50,13.56)(7.44,11.86)(9.50,9.68)
\pspolygon[fillcolor=lightgray,fillstyle=solid](8.50,13.56)(6.38,10.17)(9.00,11.62)
\pspolygon[fillcolor=lightgray,fillstyle=solid](6.38,10.17)(9.00,11.62)(6.88,8.23)
\pspolygon[fillcolor=lightgray,fillstyle=solid](6.38,10.17)(6.88,8.23)(4.25,6.78)
\pspolygon[fillcolor=lightgray,fillstyle=solid](4.25,6.78)(6.88,8.23)(4.75,4.84)
\pspolygon[fillcolor=lightgray,fillstyle=solid](4.25,6.78)(4.75,4.84)(2.12,3.39)
\pspolygon[fillcolor=lightgray,fillstyle=solid](2.12,3.39)(4.75,4.84)(2.62,1.45)
\pspolygon[fillcolor=lightgray,fillstyle=solid](2.12,3.39)(2.62,1.45)(-0.00,0.00)
\pspolygon[fillcolor=lightgray,fillstyle=solid](-0.00,0.00)(2.62,1.45)(4.00,-0.00)
\pspolygon[fillcolor=lightgray,fillstyle=solid](4.00,-0.00)(2.62,1.45)(6.62,1.45)
\pspolygon[fillcolor=lightgray,fillstyle=solid](4.00,-0.00)(6.62,1.45)(8.00,-0.00)
\pspolygon[fillcolor=lightgray,fillstyle=solid](5.25,2.90)(9.25,2.90)(6.62,1.45)
\pspolygon[fillcolor=lightgray,fillstyle=solid](9.25,2.90)(6.62,1.45)(10.62,1.45)
\pspolygon[fillcolor=lightgray,fillstyle=solid](6.62,1.45)(10.62,1.45)(8.00,0.00)
\pspolygon[fillcolor=lightgray,fillstyle=solid](10.62,1.45)(12.00,0.00)(8.00,0.00)
\pspolygon[fillcolor=lightgray,fillstyle=solid](12.00,-0.00)(9.94,2.18)(11.00,3.87)
\pspolygon[fillcolor=lightgray,fillstyle=solid](11.00,3.87)(9.94,2.18)(8.94,6.05)
\pspolygon[fillcolor=lightgray,fillstyle=solid](11.00,3.87)(8.94,6.05)(10.00,7.75)
\pspolygon[fillcolor=lightgray,fillstyle=solid](10.00,7.75)(7.88,4.36)(7.38,6.29)
\pspolygon[fillcolor=lightgray,fillstyle=solid](10.00,7.75)(7.38,6.29)(9.50,9.68)
\pspolygon[fillcolor=lightgray,fillstyle=solid](9.50,9.68)(7.38,6.29)(6.88,8.23)
\pspolygon[fillcolor=lightgray,fillstyle=solid](9.00,11.62)(9.50,9.68)(6.88,8.23)
\pspolygon[fillcolor=lightgray,fillstyle=solid](6.88,8.23)(4.75,4.84)(7.38,6.29)
\pspolygon[fillcolor=lightgray,fillstyle=solid](4.75,4.84)(7.38,6.29)(5.25,2.90)
\pspolygon[fillcolor=lightgray,fillstyle=solid](5.25,2.90)(7.88,4.36)(7.38,6.29)
\pspolygon[fillcolor=lightgray,fillstyle=solid](4.75,4.84)(2.62,1.45)(5.25,2.90)
\pspolygon[fillcolor=lightgray,fillstyle=solid](2.62,1.45)(5.25,2.90)(6.62,1.45)
\pspolygon[fillcolor=lightgray,fillstyle=solid](5.25,2.90)(9.25,2.90)(7.88,4.36)
\pspolygon[fillcolor=lightgray,fillstyle=solid](7.88,4.36)(9.94,2.18)(8.94,6.05)
\endpspicture}

\def\ThreeTilingTwo{
\pspicture(1.5,1.15)
\psset{unit=2cm}
\newrgbcolor{lightblue}{0.8 0.8 1}
\newrgbcolor{pink}{1 0.8 0.8}
\newrgbcolor{lightgreen}{0.8 1 0.8}
\newrgbcolor{lightyellow}{1 1 0.8}
\pspolygon[fillstyle=solid,linewidth=1pt,fillcolor=lightblue](0.00,0.00)(0.87,1.50)(0.87,0.50)
\psline(0.87,1.50)(0.87,0.50)
\psline(0.00,0.00)(0.00,0.00)
\psline(0.87,0.50)(0.00,0.00)
\psline(0.87,1.50)(0.87,1.50)
\psline(0.00,0.00)(0.87,1.50)
\psline(0.87,0.50)(0.87,0.50)
\pspolygon[fillstyle=solid,linewidth=1pt,fillcolor=lightyellow](0.87,1.50)(0.87,0.50)(1.73,0.00)
\psline(0.87,0.50)(1.73,0.00)
\psline(0.87,1.50)(0.87,1.50)
\psline(1.73,0.00)(0.87,1.50)
\psline(0.87,0.50)(0.87,0.50)
\psline(0.87,1.50)(0.87,0.50)
\psline(1.73,0.00)(1.73,0.00)
\pspolygon[fillstyle=solid,linewidth=1pt,fillcolor=lightgreen](0.00,0.00)(0.87,0.50)(1.73,0.00)
\psline(0.87,0.50)(1.73,0.00)
\psline(0.00,0.00)(0.00,0.00)
\psline(1.73,0.00)(0.00,0.00)
\psline(0.87,0.50)(0.87,0.50)
\psline(0.00,0.00)(0.87,0.50)
\psline(1.73,0.00)(1.73,0.00)
\endpspicture}

\def\TwelveTilingTwo{
\pspicture(1.5,1.15)
\psset{unit=2cm}
\newrgbcolor{lightblue}{0.8 0.8 1}
\newrgbcolor{pink}{1 0.8 0.8}
\newrgbcolor{lightgreen}{0.8 1 0.8}
\newrgbcolor{lightyellow}{1 1 0.8}
\pspolygon[fillstyle=solid,linewidth=1pt,fillcolor=lightblue](0.00,0.00)(0.87,1.50)(0.87,0.50)
\psline(0.87,1.50)(0.87,0.50)
\psline(0.43,0.75)(0.43,0.25)
\psline(0.00,0.00)(0.00,0.00)
\psline(0.87,0.50)(0.00,0.00)
\psline(0.87,1.00)(0.43,0.75)
\psline(0.87,1.50)(0.87,1.50)
\psline(0.00,0.00)(0.87,1.50)
\psline(0.43,0.25)(0.87,1.00)
\psline(0.87,0.50)(0.87,0.50)
\pspolygon[fillstyle=solid,linewidth=1pt,fillcolor=lightyellow](0.87,1.50)(0.87,0.50)(1.73,0.00)
\psline(0.87,0.50)(1.73,0.00)
\psline(0.87,1.00)(1.30,0.75)
\psline(0.87,1.50)(0.87,1.50)
\psline(1.73,0.00)(0.87,1.50)
\psline(1.30,0.25)(0.87,1.00)
\psline(0.87,0.50)(0.87,0.50)
\psline(0.87,1.50)(0.87,0.50)
\psline(1.30,0.75)(1.30,0.25)
\psline(1.73,0.00)(1.73,0.00)
\pspolygon[fillstyle=solid,linewidth=1pt,fillcolor=lightgreen](0.00,0.00)(0.87,0.50)(1.73,0.00)
\psline(0.87,0.50)(1.73,0.00)
\psline(0.43,0.25)(0.87,0.00)
\psline(0.00,0.00)(0.00,0.00)
\psline(1.73,0.00)(0.00,0.00)
\psline(1.30,0.25)(0.43,0.25)
\psline(0.87,0.50)(0.87,0.50)
\psline(0.00,0.00)(0.87,0.50)
\psline(0.87,0.00)(1.30,0.25)
\psline(1.73,0.00)(1.73,0.00)
\endpspicture}

\def\TwelveTilingOne{
\pspicture(1.5,1.15)
\psset{unit=2cm}
\newrgbcolor{lightblue}{0.8 0.8 1}
\newrgbcolor{pink}{1 0.8 0.8}
\newrgbcolor{lightgreen}{0.8 1 0.8}
\newrgbcolor{lightyellow}{1 1 0.8}
\pspolygon[fillstyle=solid,linewidth=1pt,fillcolor=lightblue](0.00,0.00)(0.43,0.75)(0.43,0.25)
\psline(0.43,0.75)(0.43,0.25)
\psline(0.00,0.00)(0.00,0.00)
\psline(0.43,0.25)(0.00,0.00)
\psline(0.43,0.75)(0.43,0.75)
\psline(0.00,0.00)(0.43,0.75)
\psline(0.43,0.25)(0.43,0.25)
\pspolygon[fillstyle=solid,linewidth=1pt,fillcolor=lightyellow](0.00,0.00)(0.87,0.00)(0.43,0.25)
\psline(0.87,0.00)(0.43,0.25)
\psline(0.00,0.00)(0.00,0.00)
\psline(0.43,0.25)(0.00,0.00)
\psline(0.87,0.00)(0.87,0.00)
\psline(0.00,0.00)(0.87,0.00)
\psline(0.43,0.25)(0.43,0.25)
\pspolygon[fillstyle=solid,linewidth=1pt,fillcolor=lightgreen](0.87,0.00)(0.43,0.75)(0.43,0.25)
\psline(0.43,0.75)(0.43,0.25)
\psline(0.87,0.00)(0.87,0.00)
\psline(0.43,0.25)(0.87,0.00)
\psline(0.43,0.75)(0.43,0.75)
\psline(0.87,0.00)(0.43,0.75)
\psline(0.43,0.25)(0.43,0.25)
\pspolygon[fillstyle=solid,linewidth=1pt,fillcolor=lightblue](0.43,0.75)(1.30,0.75)(0.87,0.50)
\psline(1.30,0.75)(0.87,0.50)
\psline(0.43,0.75)(0.43,0.75)
\psline(0.87,0.50)(0.43,0.75)
\psline(1.30,0.75)(1.30,0.75)
\psline(0.43,0.75)(1.30,0.75)
\psline(0.87,0.50)(0.87,0.50)
\pspolygon[fillstyle=solid,linewidth=1pt,fillcolor=lightyellow](0.43,0.75)(0.87,0.00)(0.87,0.50)
\psline(0.87,0.00)(0.87,0.50)
\psline(0.43,0.75)(0.43,0.75)
\psline(0.87,0.50)(0.43,0.75)
\psline(0.87,0.00)(0.87,0.00)
\psline(0.43,0.75)(0.87,0.00)
\psline(0.87,0.50)(0.87,0.50)
\pspolygon[fillstyle=solid,linewidth=1pt,fillcolor=lightgreen](0.87,0.00)(1.30,0.75)(0.87,0.50)
\psline(1.30,0.75)(0.87,0.50)
\psline(0.87,0.00)(0.87,0.00)
\psline(0.87,0.50)(0.87,0.00)
\psline(1.30,0.75)(1.30,0.75)
\psline(0.87,0.00)(1.30,0.75)
\psline(0.87,0.50)(0.87,0.50)
\pspolygon[fillstyle=solid,linewidth=1pt,fillcolor=lightblue](0.87,0.00)(1.30,0.75)(1.30,0.25)
\psline(1.30,0.75)(1.30,0.25)
\psline(0.87,0.00)(0.87,0.00)
\psline(1.30,0.25)(0.87,0.00)
\psline(1.30,0.75)(1.30,0.75)
\psline(0.87,0.00)(1.30,0.75)
\psline(1.30,0.25)(1.30,0.25)
\pspolygon[fillstyle=solid,linewidth=1pt,fillcolor=lightyellow](0.87,0.00)(1.73,0.00)(1.30,0.25)
\psline(1.73,0.00)(1.30,0.25)
\psline(0.87,0.00)(0.87,0.00)
\psline(1.30,0.25)(0.87,0.00)
\psline(1.73,0.00)(1.73,0.00)
\psline(0.87,0.00)(1.73,0.00)
\psline(1.30,0.25)(1.30,0.25)
\pspolygon[fillstyle=solid,linewidth=1pt,fillcolor=lightgreen](1.73,0.00)(1.30,0.75)(1.30,0.25)
\psline(1.30,0.75)(1.30,0.25)
\psline(1.73,0.00)(1.73,0.00)
\psline(1.30,0.25)(1.73,0.00)
\psline(1.30,0.75)(1.30,0.75)
\psline(1.73,0.00)(1.30,0.75)
\psline(1.30,0.25)(1.30,0.25)
\pspolygon[fillstyle=solid,linewidth=1pt,fillcolor=lightblue](0.43,0.75)(0.87,1.50)(0.87,1.00)
\psline(0.87,1.50)(0.87,1.00)
\psline(0.43,0.75)(0.43,0.75)
\psline(0.87,1.00)(0.43,0.75)
\psline(0.87,1.50)(0.87,1.50)
\psline(0.43,0.75)(0.87,1.50)
\psline(0.87,1.00)(0.87,1.00)
\pspolygon[fillstyle=solid,linewidth=1pt,fillcolor=lightyellow](0.43,0.75)(1.30,0.75)(0.87,1.00)
\psline(1.30,0.75)(0.87,1.00)
\psline(0.43,0.75)(0.43,0.75)
\psline(0.87,1.00)(0.43,0.75)
\psline(1.30,0.75)(1.30,0.75)
\psline(0.43,0.75)(1.30,0.75)
\psline(0.87,1.00)(0.87,1.00)
\pspolygon[fillstyle=solid,linewidth=1pt,fillcolor=lightgreen](1.30,0.75)(0.87,1.50)(0.87,1.00)
\psline(0.87,1.50)(0.87,1.00)
\psline(1.30,0.75)(1.30,0.75)
\psline(0.87,1.00)(1.30,0.75)
\psline(0.87,1.50)(0.87,1.50)
\psline(1.30,0.75)(0.87,1.50)
\psline(0.87,1.00)(0.87,1.00)
\endpspicture}

\title{Triangle Tiling II: Nonexistence theorems}         
\author{Michael Beeson}        
\date{\today}          
\maketitle
\newtheorem{theorem}{Theorem}
\newtheorem{lemma}{Lemma}
\newtheorem{corollary}{Corollary}
\newtheorem{definition}{Definition}
\abstract{
An $N$-tiling of triangle $ABC$ by triangle $T$ is a way of writing $ABC$ as a union of $N$ triangles
congruent to $T$, overlapping only at their boundaries.   The triangle $T$ is the ``tile''.
  The tile may or may not be similar to $ABC$.  We wish to understand possible tilings by completely characterizing 
the triples $(ABC, T, N)$ such that $ABC$ can be $N$-tiled by $T$.  In particular, this understanding 
should enable us to specify for which $N$ there exists a tile $T$ and a triangle $ABC$ that is $N$-tiled by $T$;
or given $N$,  determine which tiles and triangles can be used for $N$-tilings; or given $ABC$, to determine 
which tiles and $N$ can be used to $N$-tile $ABC$.   This is the second of four papers on this 
subject.  In \cite{beeson1}, we dealt with the case when $ABC$ is similar to $T$, and the case when $T$
is a right triangle.  In this paper, we assume that $ABC$ is not similar to $T$, and 
$T$ is not a right triangle, and furthermore that if
$T$ has a $120^\circ$ angle then $T$ is  isosceles.  

The  main theorem is that under those hypotheses, the only 
$N$-tilings are of an equilateral triangle by an isosceles tile with base angles $\pi/6$, and $N$ three times a square.

Under those hypotheses,  
there are only two families of tilings.
There are tilings of an equilateral triangle $ABC$ by an isosceles tile with base angles $\pi/6$;  
and there are newly-discovered ``triquadratic tilings'', which are treated in the third paper 
\cite{beeson-triquadratics}
in this series.
These arise in the special case that $3\alpha + 2\beta = \pi$ or $3\beta+ 2\alpha = \pi$ with $\alpha$ not a rational multiple of $\pi$, where $\alpha$ and $\beta$ are the two smallest angles of the tile, and the sides of the 
tile have rational ratios.  
 
The case when the tile has an angle $2\pi/3$ and is not isosceles is taken up in \cite{beeson120}.
At the time of writing, our analysis of that case is incomplete, but there are no known   
 tilings of any $ABC$ in that case.  
  
We use techniques from linear algebra and elementary field theory,  and in one case we use some 
algebraic number theory.  We use some  counting arguments and some elementary geometry 
and trigonometry.    
}

\section{Introduction}

For a general introduction to the problem of triangle tiling, see \cite{beeson1}. This paper is 
largely concerned with the non-existence of tilings, rather than their existence.  In \cite{beeson1}
we enumerated the known tilings; it is our aim to prove that these families exhaust all the possible 
tilings, or at least, exhaust all the triples  $(ABC, N, T)$ such that $ABC$ can be $N$-tiled by 
tile $T$.  There will be few pictures of beautiful tilings here (though there will be at least {\em some});
most of the paper is filled with calculations dedicated to showing that other apparently possible tilings
actually are not possible.

We begin by studying a number of special cases that we will need later;  then we finish 
with a frontal attack on the general case, in which the tile   is neither a right triangle nor  similar to 
the tiled triangle $ABC$.   This methodical approach succeeds in all but two cases.  One
of those cases is treated
in the third paper in this series \cite{beeson-triquadratics}, and the other is 
taken up in \cite{beeson120}.

The following lemma identifies those relatively few rational multiples of $\pi$ that have rational tangents or 
whose sine and cosine satisfy a polynomial of low degree over $\Q$.  The lemma and its proof are 
of course well-known, but it is short and may help to make the paper more self-contained.

\begin{lemma}  \label{lemma:euler}
 Let $\zeta = e^{i\theta}$ be algebraic
of degree $d$ over $\Q$,  where
  $\theta$ is a rational multiple of $\pi$,  say $\theta = 2m \pi/n$, where $m$ and $n$ have no common factor.   
\smallskip

Then $d = \varphi(n)$, where $\varphi$ is the Euler totient function.
In particular if $d = 4$, which is the case when $\tan \theta$ is rational and $\sin \theta$ is not,
then $n$ is 5, 8, 10,  or 12; and if $d=8$ then $n$ is 15, 16, 20, 24, or 30.
\end{lemma}

\noindent{\em Remark}. For example, if $\theta = \pi/6$, we have $\sin \theta = 1 / 2$, which 
is of degree 1 over $\Q$.  Since $\cos \theta = \sqrt 3 / 2$, the number
$\zeta = e^{i\theta}$ is in $\Q(i,\sqrt{3})$, which is of degree 4 over $\Q$.   The number $\zeta$ is a 
12-th root of unity, i.e. $n$ in the theorem is 12 in this case;  so the minimal polynomial of $\zeta$ is
 of degree $\varphi(12) = 4$.   This example shows that the theorem is best possible.
\smallskip

\noindent{\em Remark}. The hypothesis that $\theta$ is a rational multiple of $\pi$ cannot be dropped.
For example, $x^4-2x^3+x^2-2x+1$ has two roots on the unit circle and two off the unit circle.
\smallskip

\noindent{\em Proof}.
Let $f$ be a polynomial with rational coefficients of degree $d$ satisfied by $\zeta$.
 Since $\zeta = e^{i 2m\pi/n}$,
$\zeta$ is an $n$-th root of unity, so its minimal polynomial has degree $d = \varphi(n)$,
where $\varphi$ is the Euler totient function.  Therefore $\varphi(n) \le d$.  If $\tan \theta$
is rational and $\sin \theta$ is not, then  $\sin \theta$ has degree 2 over $\Q$, so $\zeta$
has degree 2 over $Q(i)$, so $\zeta$ has degree 4 over $\Q$.  The stated values
of $n$ for the cases $d=4$ and $d=8$ follow from the well-known formula for $\varphi(n)$. That completes the 
proof of (ii) assuming (i).

\begin{corollary} \label{corollary:euler}
If $\sin \theta$ or $\cos \theta$ is rational, and $\theta < \pi$ is a rational multiple of $\pi$,
then $\theta$ is a multiple of $2\pi/n$ where $n$ is 5, 8, 10, or 12.
\end{corollary}

\noindent{\em Proof}. Let $\zeta = \cos \theta + i \sin \theta = e^{i\theta}$.  Under the stated hypotheses,
the degree of $\Q(\zeta)$ over $\Q$ is 2 or 4.  Hence, by the lemma, $\theta$ is a multiple of $2\pi/n$,
where $n = 5$, 8, 10, or 12 (if the degree is 4) or $n = 3$ or $6$ (if the degree is 3).  But the cases
3 and 6 are superfluous, since then $\theta$ is already a multiple of $2\pi/12$.

In \cite{beeson1}, we introduced the $\d$ matrix and the $\d$ matrix equation,
$$ \d \vector a b c = \vector X Y Z $$
where $a$, $b$, and $c$ are the sides of the tile, and $X$, $Y$, and $Z$ are the lengths of the 
sides of $ABC$,
in order of size.  The angles of $ABC$ are, in order of size, $A$, $B$, and $C$,  so $X = \overline{BC}$,
$Y = \overline{AC}$, and $Z = \overline{AB}$.  We keep this convention even if some the angles are equal.
The $\d$ matrix has nonnegative integer entries, describing how the sides of $ABC$ are composed of 
edges of tiles. 

The $d$ matrix is used in almost all our proofs.  To avoid having every page filled with cumbersome 
subscript notation $\d_{ij}$ for the entries of the matrix, we introduce letters for the entries. 
While this eliminates subscripts, it does require the reader to remember which element is denoted by 
which letter.  Here, for reference, we define 

$$  \d = \matrix p d e g m f h \ell r $$

\section{An isosceles tile $T$ with base angle $\pi/6$}

\begin{figure}[ht]
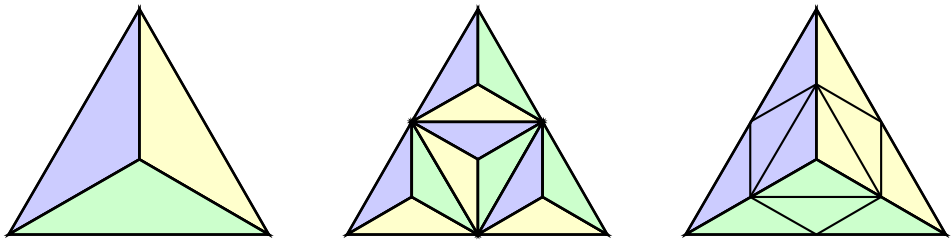
  
\caption{Tilings by the isosceles tile with base angle $\pi/6$}
\label{figure:12-tilings}
\hskip 1.3cm 
\ThreeTilingTwo
\TwelveTilingOne
\TwelveTilingTwo
\end{figure}

In this section, we investigate the tilings that can be constructed from the tile $T$ with two $\pi/6$ angles and 
one $2\pi/3$ angle, in which $T$ is not similar to $ABC$. One can find examples of such tilings by first quadratically 
tiling an equilateral triangle, and then 3-tiling each tile; the first 12-tiling shown Fig.~\ref{figure:12-tilings}
is of that form.
  
\begin{figure}[ht]
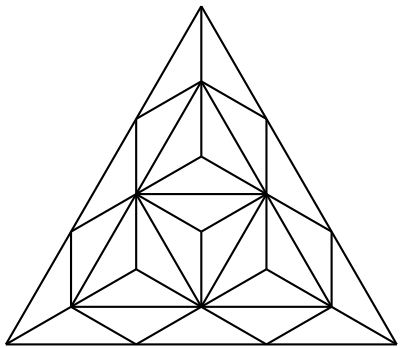
   
\caption{A prime 27-tiling by the isosceles tile with base angle $\pi/6$}
\label{figure:27-tiling}
\vskip 0.1in
\hskip 5cm
\TwentySevenTiling
\end{figure}

 But this does not exhaust the possibilities, as we can also 
3-tile an equilateral triangle $ABC$ and then quadratically tile the three resulting triangles. In this way we
can produce, for example, a (different) 12-tiling, illustrated in the second tiling in Fig.~\ref{figure:12-tilings}.
 These tilings are not prime; that is, they are made by first constructing one tiling and then tiling the tile.
 We once thought
that the equilateral 3-tiling might be the only prime tiling constructed from $T$.  That is incorrect, however, as 
shown by the prime 27-tiling in Fig.~\ref{figure:27-tiling}.

\begin{lemma} \label{lemma:piover6}
 Let triangle $T$ be isosceles with a vertex angle of $2\pi/3$.  Suppose triangle $ABC$ is $N$-tiled by $T$ and
not similar to $T$.  Then $ABC$ is equilateral, and $N$ has the form $3m^2$ for some integer $m$.
\end{lemma}

\noindent{\em Remark}.  For each integer $N$ of the form mentioned in the lemma, there does exist an $N$-tiling 
of an equilateral triangle by the tile $T$ of the lemma, as discussed above.
\smallskip

\noindent{\em Proof}.
Let $\alpha = \pi/6$ and $\gamma = 2\pi/3$ be the angles of $T$.  We note that no tile has its $\gamma$ angle at a 
vertex of $ABC$, since that would force the other two vertices of $ABC$ to have angle $\alpha$ and hence $ABC$ would 
be similar to $T$, contrary to hypothesis.
Hence the vertex angles of $ABC$ are composed (altogether) of six $\alpha$ angles.  The only two ways of 
writing 6 as a sum of smaller positive integers are $6=1+2+3$ and $6=2+2+2$.   Hence,   $ABC$ is either equilateral,
or it is a 30-60-90 right triangle. 

We are free to choose the size of $T$; 
let us suppose that the sides of $T$ are $a=1$ and $c =  \sqrt 3$, so the altitude of $T$ is $1/2$ and 
its area is $\sqrt 3 /4$  (because its base is $\sqrt 3$ and its height is $1/2$). 
If $X$, $Y$, and $Z$ are the sides of $ABC$, then the fundamental tiling equation is
  $$ \d \vector 1 1 {\sqrt 3} = \vector X  Y  Z$$
Our equations will be slightly simpler if we recall the convention that when the tile is isosceles, as it is in this proof,
so that $b=a$, that we take the middle column of the $\d$ matrix to be zero (counting short sides as $a$ rather than $b$,
since it is arbitrary).

We begin by proving that it is not the case that along one side of $ABC$,  all the tiles sharing that side of $ABC$ have an $a$ edge on the 
boundary of $ABC$.  Assume, for proof by contradiction,  that all the tiles on side $BC$ have an $a$ edge on $BC$. 
Then each of these tiles has a $\gamma$ angle at one of its vertices on $BC$.  But there cannot be a $\gamma$ 
angle at the endpoints $B$ and $C$.  Hence, there is one more $\gamma$ angle than vertices on (the interior of) $BC$,
so by the pigeonhole principle, some vertex has more than one $\gamma$ angle.  But this is impossible since
$\gamma > \pi/2$.  That contradiction shows that it is impossible for all the tiles with an edge on $BC$ to have an $a$ edge
on $BC$.  The same argument applies to each side of $ABC$, since all the angles of $ABC$ are less than $\gamma$.
To phrase this result in terms of the $\d$ matrix:  it is not the case that any entry in the third column is 
zero.  That is, $e \neq 0$, $f \neq 0$, and $r \neq 0$.

 We now take up the case when $ABC$ is equilateral. 
Let $X$ be the side of equilateral triangle $ABC$.  Then the altitude of $ABC$ is $X \sqrt 3/2$, and its area  is $X^2 \sqrt 3/4$.
But this area is also $N$ times the area of $T$.  Since the area of $T$ is $\sqrt 3/4$, we have $N \sqrt 3 /4= X^2 \sqrt 3/4$.  
Hence $N = X^2$, and $X =  \sqrt{N}$.  The tiling equations are
$$ \d \vector 1 1 {\sqrt 3} = \sqrt N \vector 1 1 1 .$$
We can assume that the components of this vector equation correspond to edges $AB$, $AC$, and $BC$ respectively.
Then, remembering that the middle column of $\d$ is zero because the tile is isosceles,  we have
\begin{eqnarray*}
\sqrt N &=& p +      e  \sqrt 3  \\
N &=& (p  + e  \sqrt 3  )^2 \\
&=& p^2 +  3e^2  +  2pe \sqrt 3 
\end{eqnarray*}
Since $\sqrt 3$ is irrational, we have $pe=0$.  Since we proved above that $e \neq 0$, we have $p = 0$.   
That is, no $a$ sides of the tile occur on the side $BC$  of $ABC$ corresponding to the first row of the tiling equation. 
Then $BC$ is entirely composed of $c$ edges (as it 
indeed is in the 27-tiling).   Since each side of $ABC$ has length $X$ and $c = \sqrt 3$, the number of tiles on each edge 
is $X/\sqrt 3$, which must be an integer, say $m$.  Then $X^2 = 3m^2$, so $N = X^2 = 3m^2$.  
That completes the proof in case $ABC$ is equilateral.   

We may therefore assume that $ABC$ is a 30-60-90 triangle.  
   Let the $\pi/6$ angle of $ABC$ be at $A$, the right angle   at $C$ and the $\pi/3$ angle at $B$. 
 Let the sides of $ABC$ be $X$, $X \sqrt 3 $, and $2X$.  
Then the area of $ABC$ is $X^2 \sqrt 3 /2$, which must be equal to the area of $N$ tiles, namely $N \sqrt 3 / 4$.
Hence $X^2 = \frac N 2$, or $X = \sqrt{N/2}$.  
Then we have 
$$ \d \vector 1 1 {\sqrt 3} =  \sqN \vector 1 {\sqrt 3} 2. $$
From the first row of this equation, we have 
$$ p + e \sqrt 3 = \sqN.$$
 Squaring this equation we have 
\begin{eqnarray*}
\frac N 2  &=& (p+ e\sqrt 3 )^2 \\
 &=& p^2 +  3e^2 +  2pe \sqrt 3
\end{eqnarray*}
Since the left side is rational, we have $pe=0$, and since $e \neq 0$, we have $p=0$.
Then 
\begin{eqnarray}
 \frac N 2 &=& 3e^2 \nonumber \\
 N &=& 6e^2.  \label{eq:1405}
 \end{eqnarray}
From the second row of the $\d$ matrix equation, we have 
$$g + f \sqrt 3 = \sqN \sqrt 3.$$
Squaring this equation, we have
\begin{eqnarray*}
g^2 + 3f^2 + 2gf \sqrt 3 &=& 3N
\end{eqnarray*}
Since $\sqrt 3$ is irrational, we have $gf = 0$.  Since $f \neq 0$ we have $g = 0$. 
Then 
\begin{eqnarray*}
3f^2 &=& 3N \\
N &=& f^2
\end{eqnarray*}
But by (\ref{eq:1405}), we have $N = 6e^2$.  Hence $f^2 = 6e^2$, so $\sqrt{6} = f/e$,  contradicting 
the irrationality of $\sqrt 6$. 
That completes the proof of the lemma.

\section{An isosceles tile $T$ with base angle  $\pi/5$}

\begin{lemma} \label{lemma:piover5}
Let $\alpha = \pi/5$.  Let $T$ be an isosceles triangle with two angles $\alpha$ and one angle $3\alpha$.
Let $ABC$ be an isosceles triangle with base angles $2\alpha$ and vertex angle $\alpha$.  Then there is 
no tiling of $ABC$ by tile $T$.
\end{lemma}

\noindent{\em Proof}. 
Let $a = \sin \alpha$.  We will work in the field $\Q(a)$ so we begin with some identities in that field.
We have 
\begin{eqnarray*}
 a &= &\sin \alpha\  = \ \frac 1 4 \sqrt{10 - 2 \sqrt 5}\\
a^2 &=& \frac 1 { 16} (10-2 \sqrt 5) \\
16a^2 - 10 &=& -2\sqrt 5 \\
\sqrt 5 &=& 5-8a^2 \\
(8a^2-5)^2 &=& 5 \\
64a^4 - 80a^2 + 25 &=& 5 \\
16a^4 - 20a^2 + 5 &=& 0
\end{eqnarray*}
This is the minimal polynomial of $a$ over $\Q$.  It is irreducible by Eisenstein's criterion, since 
20 and 5 are divisible by 5,  but 16 is not, and 5 is not divisible by $5^2$.   So $\Q(a)$ has degree 4 over $\Q$
and $\{1,a,a^2,a^3\}$ is a basis for $\Q(a)$ over $\Q$.   From the minimal polynomial we see 
\begin{equation} \label{eq:afourthpi5}
a^4 = 20a^2-5
\end{equation}
We have 
\begin{eqnarray*}
b &=& \cos \alpha \\
&& \frac {1 + \sqrt 5} 4 \\
&=& \frac {1 + (5-8a^2)} 4 \mbox{\qquad since $\sqrt 5 - 5-8a^2$} \\
&=& \frac{3-4a^2} 2 
\end{eqnarray*}
Hence
\begin{equation} \label{eq:bpi5}
b = \frac 3 2 - 2a^2 
\end{equation}
We have $\sin 2\alpha = 2ab = 2a(\frac 3 2 - 2a^2$.  Hence
\begin{equation} 
\sin 2\alpha = 3a-4a^3
\end{equation}
Recall (or prove) the trig identity (which is true for any $\alpha$)
$$ \frac {\sin 3\alpha} {\sin \alpha} = 3  - 4 \sin \alpha$$
Expressing this more briefly using $a$ we have 
\begin{equation} \label{eq:sin3alpha}
\frac {\sin 3 \alpha} {\sin \alpha} = 3-4a
\end{equation}

With these algebraic preliminaries in hand, we turn to the proof of the lemma.
The sides of $T$ are $\sin \alpha$ and $\sin(3\alpha)$.  The altitude of $T$ is $ \sin^2\alpha$,
so the area of $T$ is $\frac 1 2 \sin(3\alpha) \sin^2 \alpha$.  On the other hand, for some number $\lambda$,
the sides of $ABC$ are $\lambda \sin 2\alpha$ and $\lambda \sin \alpha$.  The altitude of $ABC$ is 
$\frac 1 2 \lambda \sin \alpha \sin(2\alpha)$, so its area is 
$\frac 1 2  \lambda^2 \sin \alpha \sin^2 2\alpha$.  Since the area of $ABC$ must be $N$ times the area of $T$, we have
\begin{eqnarray*}
N \frac 1 2 \sin(3\alpha) \sin^2 \alpha &=& \frac 1 2  \lambda^2 \sin \alpha \sin^2 2\alpha \\
N \sin(3\alpha) \sin \alpha &=& \lambda^2 \sin^2 2\alpha \\
\lambda = \frac { \sqrt{N \sin \alpha \sin(3\alpha)}} { \sin(2\alpha)}
\end{eqnarray*}
Let $\d$ be the $\d$ matrix of the tiling.   Since the tile is isosceles, we only have sides $a$ and $c$, not $b$, 
and the middle column of the $\d$ matrix is zero.    Then we have
\begin{eqnarray*}
\d \vector {\sin \alpha}{\sin \alpha} {\sin 3\alpha} &=& \lambda \vector {\sin \alpha}{ \sin 2\alpha} {\sin 2\alpha} \\
&=& \frac { \sqrt{N \sin \alpha \sin(3\alpha)}} { \sin(2\alpha)} \vector  {\sin \alpha}{ \sin 2\alpha} {\sin 2\alpha} 
\end{eqnarray*}
Writing the second of these three equations out we have 
\begin{eqnarray*}
g\sin \alpha + f {\sin 3\alpha} &=&   { \sqrt{N \sin \alpha \sin(3\alpha)}} \\
\end{eqnarray*}
Dividing both sides   by $\sin \alpha$ we have 
\begin{eqnarray*}
g  + f\frac {\sin 3\alpha} {\sin \alpha} &=& \sqrt{N \frac {\sin 3\alpha} {\sin \alpha}}\\
\end{eqnarray*}
Expressing this in terms of    $a=\sin \alpha$, we have 
\begin{eqnarray*}
g + f(3-4a) &=& \sqrt{N(3-4a)}
\end{eqnarray*}
 Squaring both sides we have 
\begin{eqnarray*}
(g  + f(3 -4a))^2 &=&  N (3-4a)\\
g^2 + 2 gf (3-4a) + f^2(3-4a)^2 &=& N(3-4a) \\
g^2 + 6g f + 9 f^2 + a(-8gf - 24f^2) + a^2(16f^2) &=& 3N-4aN
\end{eqnarray*}
Since $\{1,a,a^2,a^3\}$ is a basis for $\Q(a)$ over $\Q$, the coefficients of like powers of $a$ must be equal 
on both sides.  Therefore the coefficient of $a^2$ on the left is zero; that is, $f = 0$.  
Then the equation becomes
$$ g^2  = 3N-4aN$$
But now there is a nonzero coefficient of $a$ on the right, and no $a$ term on the left.  This 
contradiction completes the proof of the lemma.

\section{A non-isosceles tile $T$ with largest angle $2\pi/5$}
\begin{lemma}\label{lemma:twopioverfive}
Suppose $ABC$ is $N$-tiled by a triangle $T$.  Suppose $ABC$ is 
not similar to $T$, and $T$  is not isosceles.   Then the largest angle $\gamma$ of $T$ 
cannot be $2\pi/5$.  
\end{lemma}

\noindent{\em Proof}.  
Since $ABC$ is not similar to $T$, we have vertex splitting. 
As before, let $P$, $Q$, and $R$ be the total number of $\alpha$, $\beta$, and $\gamma$ angles, 
respectively, at vertices of $ABC$.   Since $T$ is not similar to $ABC$,
we have $P+Q+R \ge 5$.
Since $T$ is not isosceles,  $\alpha < \beta$.    We have $\beta > (\pi-\gamma)/2 = 3\pi/10 = (3/4) \gamma$.
The number $R$ of $\gamma$ angles 
is at most 2, since $3\gamma = 6\pi/5 > \pi$.   

Assume, for proof by contradiction, that $R=2$.  Then 
$P\alpha + Q\beta = \pi/5$, and since $\beta > 3\pi/10 > \pi/5$,
we must have $Q = 0$.  Hence $\alpha = \pi/(5P)$ and $\beta = \pi- \gamma - \alpha = 3\pi/5 - \pi/(5P)$.
Since $\beta < \gamma$, we have 
\begin{eqnarray*}
\frac 3 5 \pi - \frac \pi {5P} < \frac 2 5 \pi \\
\frac 3 5 - \frac 1 {5P} < \frac 2 5 \\
3P - 1 < 2P \\
P < 1
\end{eqnarray*}
This contradicts $P+Q \ge 5$.  That completes the proof by contradiction that $R \neq 2$.
Since $R \le 2$,  that implies $R = 1$ or $R=0$.  

Now assume, for proof by contradiction, that $R=1$.  Then the equation $P\alpha + Q\beta + R\gamma = \pi$
becomes
$P\alpha + Q \beta = 3\pi/5$, and since $\beta > 3\pi/10$, we must have $Q \le 1$.  We have 
two equations for $\alpha$ and $\beta$:
\begin{eqnarray*}
P\alpha + Q\beta &=& \frac {3\pi} 5  \qquad\mbox{since $R=1$} \\
\alpha + \beta &=& \frac{3\pi} 5  \qquad\mbox{since $\alpha + \beta + \gamma = \pi$}
\end{eqnarray*}
We cannot have $P=Q =1$, since $P + Q \ge 5$.  Hence $P \neq Q$, and the equations 
are uniquely solvable for $\alpha$ and $\beta$.  By Cramer's rule we have 
\begin{eqnarray*}
\alpha &=& \frac {3\pi}5 \frac{1-Q}{P-Q} \\
\beta &=&  \frac {3\pi}5\frac{P-1}{P-Q}
\end{eqnarray*}
Since $Q \le 1$, the equation for $\alpha$ shows that $Q=0$, since $\alpha > 0$.  Hence
\begin{eqnarray*}
\alpha &=& \frac {3\pi}{5 P} \\
\beta &=&  \frac {3\pi}5\frac{P-1}P
\end{eqnarray*}
We have $\beta < \gamma = 2\pi/5$.  Hence 
\begin{eqnarray*}
\frac {3\pi}5\frac{P-1}P&<& \frac {2\pi} 5 
\end{eqnarray*}
Dividing by $\pi$ and multiplying by $5P$, we have
\begin{eqnarray*}
3(P-1)  &<& 2P \\
3P-3  &<& 2P \\
P &<& 3 
\end{eqnarray*}
But this contradicts $P+Q+R \ge 5$, since $R=1$ and $Q=0$.  That completes the proof that $ R \neq 1$. 

Hence $R=0$.  The equation $P\alpha + Q\beta + R\beta = \pi$ then becomes $P\alpha + Q\beta = \pi$,
and we have the two equations
\begin{eqnarray*}
P\alpha + Q\beta &=& \pi \\
\alpha + \beta &=& \frac{3\pi} 5  \qquad\mbox{since $\alpha + \beta + \gamma = \pi$}
\end{eqnarray*}
If $P=Q$, then we have $\alpha + \beta = \pi/P = (3\pi)/5$.  Hence $P = 5/3$, which is 
impossible since $P$ is an integer.  Therefore $P\neq Q$, and the equations can be solved 
uniquely for $\alpha$ and $\beta$.  By Cramer's rule we have 
\begin{eqnarray*}
\alpha &=& \frac{ \pi - (3\pi/5)Q} {P-Q} \\
\beta &=& \frac { (3\pi/5)P - Q\pi}{P-Q}
\end{eqnarray*}
Assume, for proof by contradiction, that $Q > 1$.  Then the numerator of $\alpha$ is negative, 
so the denominator $P-Q$ is also negative, so $Q > P$.   Then 
\begin{eqnarray*}
\beta &=& \frac{ Q - (3/5)P}{Q-P}  \pi
\end{eqnarray*}
Since $\beta < \gamma = 2\pi/5$, we have 
\begin{eqnarray*}
\frac {Q - (3/5)P} {Q-P} &<& \frac 2 5 \\
\end{eqnarray*}
Multiplying both sides by $5(Q-P)$, which is positive since $Q > P$, we have
\begin{eqnarray*}
5Q - 3P &<& 2Q-2P \\
3Q < P 
\end{eqnarray*}
But we have $P < Q$, so $3Q < P < Q$, which is a contradiction.  Thus our assumption that $Q > 1$
is untenable.  Hence $Q \le 1$.  

Now assume $Q=0$.  Then the formula for $\beta$ above simplifies to $3\pi/5$,  which is more than
$\gamma$, contradiction.   Hence $Q=1$.   Then we have
\begin{eqnarray*}
\beta   &=& \frac {(3/5)P-1}{P-1} \pi \\
\end{eqnarray*}
The equation $\beta < \gamma$ then implies
\begin{eqnarray*}
\frac {(3/5)P-1}{P-1} &<& \frac 2 5
\end{eqnarray*}
Multiplying both sides by $5(P-1)$ we have 
\begin{eqnarray*}
3P-5 &<& 2(P-1) \\
P &<& 3
\end{eqnarray*}
But that contradicts $P+Q+R \ge 5$, since $R=0$ and $Q=1$.
This contradiction does not depend on any assumption, and hence completes the proof of the lemma.

\section{The case $\alpha = \pi/11$ and $\beta = 3\alpha$}
In this section and the next, we prove there are no tilings of any triangle $ABC$ in which the tile $T$ is 
not similar to $ABC$, for two special tiles $T$ that seem to slip through the cracks of our 
arguments in subsequent sections.   We deal with these two special cases using some   
field theory computations in the relevant cyclotomic fields, combined with simple geometric 
arguments.   Some of the computations were made using 
a computer algebra system, but they could have been made by hand with enough patience, and 
they can be checked (once made) with a calculator, as we shall now explain.

We used the computer algebra systems Sage and Mathematica.   
The results  we obtained are of two kinds.  Some are polynomial equations in a root of unity $\zeta$, with integer coefficients.  You 
can check these with a scientific calculator, and verify that the two sides agree to many decimal places.  While 
that is not a proof, it is good evidence that no mistake has been made in preparing input for Sage or typesetting
the output, or by Sage itself.   On the other hand, some of the results are equations containing not only $\zeta$ 
but also parameters for 
unknown integers $p$, $q$, $r$, etc.   Here a computer algebra system is used to divide one polynomial (with variables 
in the coefficients) by another and get the remainder.   The computations would be very difficult, but not impossible,
to perform by hand.   To check whether a mistake has been made in the input or 
the typesetting, you would have to give the problem to Sage or Mathematica or Maple yourself. 
Mathematica and Maple are also capable of doing these calculations,  and we checked 
some of them in both Mathematica and Sage.   
    
The most curious thing about this part of the work  
is that the two special cases are surprisingly different.  The two arguments have a method 
in common (field theory), but the fields are different, and the resulting equations are 
different, and the geometric arguments at the end are different.  They do not seem to be 
special cases of a more general proof.  

\begin{lemma} \label{lemma:special-11}  Let $ABC$ be tiled using as a tile the 
triangle $T$ with $\alpha = \pi/11$, $\beta = 3\pi/11$, and $\gamma = 7\pi/11$,
and suppose the angles of triangle $ABC$ are anything but $2\alpha$, $4\alpha$, and $5\alpha$.
Then $ABC$ is similar to $T$.
\end{lemma}

\noindent{\em Remark}.  We were not able to deal with the exceptional case mentioned by 
this method of proof.  
\medskip

\noindent
{\em Proof}. 
 First we note that $\sin \alpha$ and $\cos \alpha$ are not rational,
by Corollary \ref{corollary:euler}.  Let 
\begin{eqnarray*}
a &:=& \sin\alpha \\
\zeta &:=& e^{i \alpha}
\end{eqnarray*}
Then $\zeta^{22} = 1$.  Since $22$ is not divisible by 4,  $i$ does not belong to $\Q(\zeta)$.  The degree of $\Q(\zeta)$ 
over $\Q$ is $\varphi(22) = 10$, and $\Q(\zeta,i)$ then has degree 20.
There are no convenient automorphims of $\Q(\zeta)$ that fix two of $ia$, $ib$, and $ic$ and move the third, so 
we are forced to consider the area equation directly.   

The minimal polynomial of $\zeta$ is the cyclotomic polynomial 
$$x^{10} - x^9 + x^8 - x^7 + x^6 - x^5 + x^4 - x^3 + x^2 - x + 1.$$ 
We note that  $\sin^2 \alpha$ is not rational, since if it were, then $\sin \alpha$ and $\cos \alpha$
would be of degree 2, so  $\Q(\sin \alpha, \cos \alpha)$ would be of degree 1, 2, or 4, and  $\Q(\zeta,i) = \Q(\sin \alpha, \cos \alpha, i)$ would be of degree 1, 2, 4, or 8, but we have shown that it is 20.
  Hence $\Q(a^2) \neq \Q$.
Since $a^2 = -(\zeta-\zeta^{-1})^2$ belongs to $\Q(\zeta)$,  it is of degree dividing $10$.  Since it is fixed under the 
automorphism $\sigma_{21}$ that takes $\zeta$ to $\zeta^{21} = \zeta^{-1}$,  in fact $\Q(a^2)$ is of degree 5
over $\Q$.  


The minimal polynomial of $a = \sin(\pi/11)$ is 
$$
 x^{10} - \frac {11} 4 x^8 + \frac{11} 4 x^6 - \frac {77}{64}  x^4 + \frac{55}{256}x^2 - \frac{11}{1024}
 $$
according to Sage.  We will show that $\cos \alpha$ belongs to $\Q(a)$.  We have
\begin{eqnarray*}
\sin 6 \alpha &=& \sin 5 \alpha 
\end{eqnarray*}
Expanding both sides in $\sin \alpha$ and $\cos \alpha$ (using one's favorite computer algebra system, or pencil and paper),
subtracting one from the other, 
and then using $\cos^2 \alpha = 1-\sin^2 \alpha$, one finds
\begin{eqnarray*}
0 = (32 \sin^4 \alpha \cos \alpha -16 \sin^4 \alpha-32 \sin^2 \alpha \cos \alpha +20 \sin^2 \alpha +6 \cos \alpha -5) \sin \alpha
\end{eqnarray*}
Dividing by $\sin \alpha$ we have
\begin{eqnarray*}
0 &=& 32 \sin^4 \alpha \cos \alpha -16 \sin^4 \alpha-32 \sin^2 \alpha \cos \alpha +20 \sin^2 \alpha +6 \cos \alpha -5 \\
  &=& \cos \alpha ( 32 \sin^4 \alpha -32 \sin^2 \alpha + 6)  -16 \sin^4 \alpha+20 \sin^2 \alpha  -5 
\end{eqnarray*}
Solving for $\cos \alpha$ we find
\begin{eqnarray*}
\cos \alpha &=& \frac { 16 \sin^4 \alpha -20 \sin^2 \alpha  +5 } { 32 \sin^4 \alpha -32 \sin^2 \alpha + 6}\\
&=& \frac {16 a^4 - 20 a^2+5}{32 a^4 - 32 a^2 + 6}   
\end{eqnarray*}

Squaring and using $\cos^2 \alpha = 1-\sin^2\alpha$, one soon can confirm the minimal polynomial reported by Sage;
we stop here, having proved that $\cos \alpha$ belongs to $\Q(a^2)$.
It is awkward to have $\cos \alpha$ expressed as a fraction; we want to have it expressed in terms of the basis elements of 
$\Q(a^2)$.   Sage can do that for us, and we can verify the result numerically to check that there was no human or machine error.%
\footnote{If we wanted to do it by hand, we would use the extended Euclidean algorithm to find polynomials $\lambda$ and 
$\mu$ such that $\lambda(32x^4-32x^2 + 6) + \mu \psi(x) = 1$, where $\psi$ is the the minimal polynomial of $a$; then 
the answer is $(16x^4 - 20x^2 + 5) \lambda \mod \psi$.   This is tedious by hand, easy with Mathematica or Sage, but then 
you might as well just let Sage do the whole calculation.}
\begin{equation} \label{eq:cospiover11}
 \cos \alpha = 128 a^8 - 288 a^6 + 216 a^4 - 60a^2 + \frac 9 2   
\end{equation}
We next wish to express $b = \sin \beta$ and $c = \sin \gamma$ in terms of $a = \sin \alpha$ and $\cos \alpha$. 
We have 
\begin{eqnarray*}
b &=& \sin 3\alpha\\
&=& 3a-4a^3
\end{eqnarray*}
 and 
\begin{eqnarray*}
c &=& \sin 4\alpha  \\
&=& \sin(2 \cdot 2 \alpha) \\
&=& 2 \sin 2\alpha \cos 2\alpha \\
&=& 4 \sin \alpha \cos \alpha (1-2\sin^2 \alpha) \\
&=& 4a(1-2a^2)\cos \alpha  \\
&=& 4a(1-2a^2)\big(128 a^8 - 288 a^6 + 216 a^4 - 60a^2 + \frac 9 2\big)  
\end{eqnarray*}
Now we ask Sage to evaluate this expression in $\Q(a)$.  It tells us
$$c =  -64 a^7 + 112a^5 - 56a^3 + 7a  $$
As a check, this value is numerically equal to $\sin(4\alpha)$ to many decimal places.

Let $\theta$ be one of the angles of triangle $ABC$ and let $U$ and $V$ be the adjacent sides.
Let $p$, $q$, $r$, $m$, $n$, and $\ell$ be the rows of the $\d$ matrix associated with $U$ and $V$,
i.e. the numbers of tiles with $a$, $b$, and $c$ sides on $U$ and $V$, respectively.  Explicitly,
\begin{eqnarray*}
U &=& pa + qb + rc \\
V &=& ma + nb + \ell c
\end{eqnarray*}
The area equation $2N\A_T = 2\A_{ABC}$ becomes
\begin{eqnarray*}
N abc &=& UV \sin\theta   \\
  &=& (pa +qb + rc)(ma+nb + \ell c) \sin \theta  \\
 &=& (pa + q(3a - 4a^3) + r(-64 a^7 + 112a^5 - 56a^3 + 7a))\\
&&\big(ma + n(3a-4a^3) + \ell(-64 a^7 + 112a^5 - 56a^3 + 7a)\big) \sin \theta \\
\end{eqnarray*}
Dividing both sides by $a^2 \sin \theta$ we have
\begin{eqnarray*}
\frac{Nbc} {a \sin \theta} &=& \big(p + q(3-4a^2) + r(-64 a^6 + 112a^4 - 56a^2 + 7)\big)  \\
&&\big(m + n(3a-4a^2) + \ell(-64 a^6 + 112a^4 - 56a^2 + 7)\big)  \nonumber
\end{eqnarray*}
Multiplying out, and reducing by the minimal polynomial of $a$, we have (according to Sage)
\begin{eqnarray*}
\frac{Nbc} {a \sin \theta}&=& 
(256(q + 14r)\ell + 256(14\ell + n)r - 7168\ell r)a^8  \\
&& + (-448(q + 14r)\ell - 448(14\ell + n)r - 64(p + 3q + 7r)\ell - 64(7\ell + m + 3n)r + 13376\ell r)a^6   \\
&& + (16(q + 14r)(14\ell + n) + 112(p + 3q + 7r)\ell + 112(7\ell + m + 3n)r - 4576\ell r)a^4  \\
&& + (-4(q + 14r)(7\ell + m + 3n) - 4(14\ell + n)(p + 3q + 7r) + 704\ell r)a^2 \\
&& + (p + 3q + 7r)(7\ell + m + 3n) - 33\ell r
\end{eqnarray*}
Expanding and simplifying the coefficients of $a^8$ and $a^6$ we find  
\begin{eqnarray}\label{eq:areaeqnforpiover11}
\frac{Nbc} {a \sin \theta}&=& 
256(q\ell + rn)a^8  \nonumber\\
&& -( 640 (q\ell + rn) + 64(\ell p + \ell r + mr) )a^6 \nonumber \\
&& + (16(q + 14r)(14\ell + n) + 112(p + 3q + 7r)\ell + 112(7\ell + m + 3n)r - 4576\ell r)a^4 \nonumber\\
&& + (-4(q + 14r)(7\ell + m + 3n) - 4(14\ell + n)(p + 3q + 7r) + 704\ell r)a^2 \nonumber \\
&& + (p + 3q + 7r)(7\ell + m + 3n) - 33\ell r
\end{eqnarray}
 Now suppose, for proof by contradiction, that $\theta = \gamma$.  Then $\sin \theta = c$ and we can 
cancel $c$ from numerator and denominator on the left.  
 Substituting $b = 3a-4a^3$ on the left, a factor of $a$ cancels and we obtain $3N-4Na^2$ on the left.
 Comparing coefficients of $a^8$ we have zero on the left and $q\ell + rn$ on the right; hence 
 $q\ell = 0$ and
  $rn = 0$ (since both quantities are nonnegative).  Then the coefficient of $a^6$ is zero on the 
  left, and on the right it is $64(\ell p + \ell r + mr)$. Hence $\ell p = 0 = \ell r = mr$.  Now
  we have proved that $nr = mr = \ell r = 0$, so $r(n+m+\ell) = 0$.  But $n+m+\ell > 0$, so $r = 0$.
  Then the coefficient of $a^4$ simplifies to $16nq$;  hence $nq = 0$.  We have also proved
  $\ell p = 0 = \ell q = \ell r$,  so $\ell(p+q+r) = 0$. 
    But $p+q+r > 0$, so $\ell = 0$.
The coefficient of $a^2$ simplifies to $-4(mq + np)$ on the right and $-4N$ on the left.
Therefore $mq +np = N$.  The constant term is $3N$ on the left and on the right it simplifies
to $mp + 3mp+3np$.  Hence $mp + 3mp + 3np = 3N$.   Subtracting three times the equation $mq + np= N$ we find $mp = 0$.
Assume for the moment that $m \neq 0$. Then $p = 0$; since $\ell = 0$ we must have $q \neq 0$.  Then since $nq = 0$ we have $n = 0$,
  and side $V$ is composed entirely of $a$ sides of tiles, while side $U$ is composed entirely of $b$ sides of tiles.
On the other hand if $m = 0$ then side $V$ is composed entirely of $b$ sides of tiles, so $n \neq 0$; then since $nq = 0$
we have $q =0$ and side $U$ is composed entirely of $a$ sides of tiles.  Note that these deductions are consistent
with the case when $ABC$ is similar to $T$ and the tiling is quadratic.  Interchanging the arbitrary labels $U$ and $V$
if necessary, we may assume that side $U$ is composed entirely of $a$ sides of tiles.  The other end of side $U$
is at either $A$ or $B$.  The angle of $ABC$ at $A$ cannot be $\alpha$, for then $ABC$ would be similar to $T$,
since angle $C = \gamma$; but it must be less than $\beta$, since $2\beta + \gamma > \pi$.    Also the 
angle at $B$ must be less than $\beta$,  since $A + B = \alpha + \beta$, and angle $A > \alpha$.   Hence the angle 
at the other endpoint of $U$ is less than $\beta$, whether that endpoint is $A$ or $B$.   Let $T_1$ be the 
tile along $U$ at that vertex of $ABC$.  Then $T_1$ has its $\alpha$ angle in the corner of $ABC$. But
then it cannot have its $a$ side along $U$, contradiction.  That completes the proof that  
  $\theta \neq \gamma$.   Since $\theta$ was originally any angle of $ABC$, we can now assume 
that none of the angles of $ABC$ are equal to $\gamma$. 

Returning to (\ref{eq:areaeqnforpiover11}), the last equation before we assumed $\theta = \gamma$, we now 
assume, for proof by contradiction, that $\theta = \alpha$.  Then $\sin \theta = a$ and on the left 
we have $Nbc/a^2$.   Substituting the expressions derived above for $b/a$ and $c/a$ on the left, we have 
\begin{eqnarray*} 
&&N(3-4a^2)( -64 a^6 + 112a^3 - 56a^2 + 7 ) \\
&=& 256(q\ell + rn)a^8   \\
&& -( 640 (q\ell + rn) + 64(\ell p + \ell r + mr) )a^6   \\
&& + (16(q + 14r)(14\ell + n) + 112(p + 3q + 7r)\ell + 112(7\ell + m + 3n)r - 4576\ell r)a^4  \\
&& + (-4(q + 14r)(7\ell + m + 3n) - 4(14\ell + n)(p + 3q + 7r) + 704\ell r)a^2  \\
&& + (p + 3q + 7r)(7\ell + m + 3n) - 33\ell r
\end{eqnarray*}
Expanding the left hand side we have 
\begin{eqnarray*}
&&	N(256a^8 - 640a^6 + 560a^4 - 196a^2 + 21) \\
&=& 256(q\ell + rn)a^8   \\
&& -( 640 (q\ell + rn) + 64(\ell p + \ell r + mr) )a^6   \\
&& + (16(q + 14r)(14\ell + n) + 112(p + 3q + 7r)\ell + 112(7\ell + m + 3n)r - 4576\ell r)a^4  \\
&& + (-4(q + 14r)(7\ell + m + 3n) - 4(14\ell + n)(p + 3q + 7r) + 704\ell r)a^2  \\
&& + (p + 3q + 7r)(7\ell + m + 3n) - 33\ell r
\end{eqnarray*}
Now we can equate coefficients of like powers of $a$.  From the coefficient of $a^8$ we see $N = q\ell + rn$.
Then from the coefficient of $a^6$ we have 
\begin{eqnarray*}
640 N &=& 640(q\ell + rn) + 64 (\ell p + \ell r + mr) \\
&=& 640N + 64 (\ell p + \ell r + mr) 
\end{eqnarray*}
Subtracting the left side from the right we see $\ell p = \ell r = mr = 0$.
Turning to the coefficient of $a^4$, we have 
\begin{eqnarray*}
560 N &=& 16 \cdot 14 (q \ell + rn) + 16 qn + 3 \cdot 112(q\ell + rn) \\
&=& 560(q\ell + rn) + 16qn \\
&=& 560N + 16 qn
\end{eqnarray*}
Subtracting the left side from the right, we see $qn = 0$.
Finally, from the constant term we have
\begin{eqnarray*}
21N &=& pm + 3pn + 3qm + 21(q\ell + rn) \\
&=& pm + 3pn + 3qm + 21N
\end{eqnarray*}
Subtracting the left side from the right, we see $pm = 0 = pn = qm$.  Since we
already deduced $p\ell = 0$ we now have $p(m+n+\ell) = 0$, which implies $p=0$. 
Similarly, having deduced $mp = mq = rm = 0$, we have $m(p+q+r) = 0$, which 
implies $m=0$.  Now, if $q \neq 0$ then since $nq = 0$ we would have $n=0$;
and then since $m = 0$ we would have $\ell \neq 0$, and since $\ell r = 0$ we 
would have $\ell = 0$.  Thus if $q \neq 0$ then $U$ is composed entirely of $b$ sides
of tiles (since $p=r=0$) and $V$ is composed entirely of $c$ sides of tiles (since $m=n=0$).
On the other hand, if $q = 0$ then since $p=0$, we would have $r\neq 0$, so $\ell = 0$, 
so $n \neq 0$, and then $U$ would be composed entirely of $c$ sides of tiles and 
$V$ entirely of $b$ sides.   Hence in either case,  one of $U$ and $V$ is 
composed entirely of $b$ sides; we may assume without loss of generality that it is $U$.
Let $T_1$ be the tile at vertex $A$ (there is only one since the angle there is $\alpha$).
Let $V_2$ be the next vertex on side $U$, at the other end of the $b$ side of $T_1$.  
Then $T_1$ has its $\gamma$ angle at $V_2$, since it has its $\beta$ angle opposite its $b$ 
side on $U$, and its $\alpha$ angle at $A$.  Since $\gamma > \pi/2$, there cannot be 
two $\gamma$ angles on the same side of $U$ at any one vertex.  Then by the pigeonhole
principle, the last tile  
$T_k$ must have its $\gamma$ angle 
at the vertex $C$ of $ABC$, and $U$ is the side $AC$.

  But angle $C$ cannot be exactly $\gamma$, since then $ABC$
would be similar to $T$ (and besides, we already showed none of the angles of $ABC$ are equal
to $\gamma$).   Assume, for proof by contradiction, 
 that angle $C = \gamma + \alpha$.     Let $P$ be the point on side $AB$ such that angle $ACP$ is $\gamma$; 
since angle $C$ is more  than $\gamma$, such a point $P$ exists.  Then triangle $ACP$ is similar to $T$,
and the similarity factor is $q$, since the side opposite the $\beta$ angle of $ACP$ (which is $AC$)
has length $qb$.  Hence 
the length of $AP$ is $qc$ and the area of triangle $ACP$ is $q^2 \A_T$.   Since the length of $AP$ is $qc$
and that is less than the length of $AB$, which is $\ell c$,  we have $q < \ell$. 

Observe that triangle $CPB$ is similar to triangle $ABC$, since angle $PCB$ is $\alpha$ and the two triangles share 
angle $B = 2\alpha$.  The length of $PC$
is $qa$, because $q$ is the similarity
factor between $T$ and triangle $ACP$.   Side $PC$ is  opposite angle $B$ in $PCB$; and $AC$, which has 
length $qb$, is opposite angle $B$ in triangle $ABC$.  The similarity factor between $CPB$ and $ABC$ is 
thus $a/b$.  We now decompose triangle $ABC$ into triangle $ABP$ and triangle $CPB$ (although we do not claim 
that segment $CP$ belongs to the tiling).  We have 
\begin{eqnarray*}
\A_{ABC} &=& \A_{ACP} + \A_{CPB} \\
&=& q^2 \A_T + \frac a b \A_{ABC}
\end{eqnarray*}
Solving for $\A_{ABC}$ we find 
\begin{eqnarray*}
\A_{ABC} &=& \frac {q^2}{1-\frac a b} \A_T \\
\end{eqnarray*}
On the other hand $\A_{ABC} = N \A_T$ because $ABC$ is $N$-tiled by $T$.  Hence
\begin{eqnarray*}
N \A_T &=& \frac{q^2}{1-\frac a b} \A_T 
\end{eqnarray*}
Canceling $\A_T$ and simplifying, we find 
\begin{eqnarray*}
q^2 &=& N\big(1-\frac a b\big) \\
q^2 b &=& N(b-a)
\end{eqnarray*}
Substituting $b = 3a-4a^3$ and canceling $a$ from both sides we have
\begin{eqnarray*}
q^2(3-4a^2) &=& N(2-4a^2)\\
-4q^2 a^2  + 3q^2 &=& -4Na^2 + 2N
\end{eqnarray*}
Comparing the coefficients of $a^2$ we have $N = q^2$. 
Comparing the constant coefficients we have $q^2 = 2N/3$.  Since $N > 0$ this is 
a contradiction.   Hence,  if angle $A = \alpha$, we do not have angle $C = \gamma + \alpha$.

Now suppose, for proof by contradiction, that angle $A = \alpha$ and angle $C = \gamma + 2\alpha$.
Then angle $B = \alpha$, so $ABC$ is isosceles.  Let $P$ and $Q$ be points on $AB$ 
such that angle $ACP = \gamma$ and angle $ACQ = \gamma + \alpha$.  Then angle $QCB = \alpha$.
Since angle $B = \alpha$, triangle $ABC$ is similar to triangle $CAB$, as both have two $\alpha$
angles.   Angle $PCQ = \alpha$ and angle $CPQ = \pi-\beta = \alpha +\gamma$.
Therefore triangle $PCQ$ is similar to triangle $CQA$.   Triangle $ACP$ is similar to the 
tile, with similarity factor $q$, so $PC = qa$.  In the similarity between $PCQ$ and $CQA$,
side $PC$ corresponds to $AC$, so the similarity factor from $CQA$ to $PCQ$ is 
$PC/AC = a/b$.   Then 
\begin{eqnarray*}
AC &=& qb \\
AP &=& qc \\
PC &=& qa \\
QC &=& (a/b) AQ \\
&=& (a/b)(AP + PQ) \\
&=& (a/b) qc + (a/b) PQ\\
PQ &=& (a/b) QC
\end{eqnarray*}
Substituting this value of $PQ$ into $QC = (a/b)qc + (a/b)PQ$ we have 
\begin{eqnarray*}
QC &=& (a/b) qc + (a/b)^2 QC 
\end{eqnarray*}
Solving for $QC$ we have
\begin{eqnarray*}
QC &=& \frac { abqc}{b^2-a^2}
\end{eqnarray*}
Since $PQ = (a/b) QC$ this give us 
\begin{eqnarray*}
PQ &=& \frac {a^2 qc}{b^2-a^2}
\end{eqnarray*}
Since triangle $ABC$ is isosceles we have $BC = AC = qb$.  
Since triangle $QBC$ is isosceles we have $QB = QC$.  We compute
\begin{eqnarray*}
AB &=& AP + PQ + QB \\
&=& qc + PQ + QB \\
&=& qc + \frac {a^2 qc}{b^2-a^2} + QB \\
&=& qc + \frac {a^2 qc}{b^2-a^2} + QC \mbox{since $QB=QC$} \\
&=&qc + \frac {a^2 qc}{b^2-a^2} +\frac { abqc}{b^2-a^2} \\
&=& qc\big( 1 + \frac {a^2 + ab}{b^2-a^2}\big) \\
&=& qc \big( 1 + \frac{a(b+a)}{b^2-a^2}\big)\\
&=& gc \big( 1 + \frac a{b-a} \big) \\
&=& gc \big(\frac{b-a + a}{b-a} \big) \\
&=& gc \big(\frac b {b-a}\big)
\end{eqnarray*}
But $AB$, which was the side $V$,  is composed entirely of $c$ edges.  Hence 
$gb(b-a)$ is an integer, namely $\ell$ where $V = \ell c$.  In particular $b(b-a) = \ell/g$ 
is rational. But 
$b = 3a-4a^3$, so we have 
\begin{eqnarray*}
\frac {\ell} g &=& (3a-4a^3)((3a-4a^3)-a) \\
   0     &=& (3a-4a^3)(2a-4a^3) - \frac {\ell} g
\end{eqnarray*}
But then $a$ satisfies a rational polynomial of degree 6,  contradicting the 
fact that its minimal polynomial has degree 10.   This contradiction shows that 
$ABC$ cannot have an $\alpha$ angle and a $\gamma + 2\alpha$ angle.  
But we have now checked all the possibilities for angle $C$ if angle $A = \alpha$.
Hence, angle $A$ cannot be $\alpha$.

Now assume, for proof by contradiction, that triangle $ABC$ is isosceles, with two angles 
equal to $2\alpha$.  We may suppose the two angles are angles $A$ and $B$. 
Then angle $C = \gamma$, since $\gamma + 4\alpha = \pi$.  But we have already shown 
that triangle $ABC$ cannot have a $\gamma$ angle.  Hence triangle $ABC$ cannot be 
isosceles with two $2\alpha$ angles.

The smallest angle of triangle $ABC$ cannot exceed $2\alpha$, since if it is $3\alpha$ or more then the second smallest 
angle is at least $4\alpha$, and the largest angle would exceed $5 \alpha$, making the total at least $12\alpha$,
while it must be only $11\alpha = \pi$.  

The  smallest angle of triangle $ABC$ must therefore be exactly  $2\alpha$.  Then the second
 smallest angle is at least $3\alpha = \beta$.  If the next smallest angle 
is $\beta = 3\alpha$ then the largest angle is $\gamma$, which we have shown to be impossible.   Hence the 
second smallest angle must be at least $4\alpha$.  If it is $4\alpha$ then the third angle is $5\alpha$.  The second
smallest cannot be $5\alpha$ or more because then the sum of all three angles would exceed $\pi = 11 \alpha$. 
Hence there is only one possible shape of $ABC$ remaining  when angle $A = 2\alpha$, angle $B = 4\alpha$,
and angle $C = 5\alpha$.  But that shape is ruled out by hypothesis.  That completes the proof of the lemma.

\section{The case $\alpha = \pi/14$ and $\beta = 4\alpha$}
Here is a another special case that requires a separate treatment, which we provide by using a cyclotomic field.

\begin{lemma} \label{lemma:special-14}  Let $ABC$ be tiled using as a tile the 
triangle $T$ with $\alpha = \pi/14$, $\beta = 4\pi/14$, and $\gamma = 9\pi/14$.
Suppose angle $C$ is equal to $\gamma$ or $\gamma + \alpha$.  
Then $ABC$ is similar to $T$.
\end{lemma}

\noindent{\em Proof}. Let $\zeta = e^{2\pi i/28}$.  The degree of $\Q(\zeta)$ over $\Q$ is $\varphi(28) = 12$.
Since 28 is divisible by 4, $i$ belongs to $\Q(\zeta)$.  Let $\sigma = \sigma_{15}$ be the automorphism of 
$\Q(\zeta)$ that takes $\zeta$ to $\zeta^{15}$.   Then $\sigma$ takes $i$ to $-i$.   Since 
$2\sin(j\alpha) = -i(\zeta^j - \zeta^{-j})$,  $\sigma$ takes $\sin(j \alpha)$ to $-\sin (15 j  \alpha)$,
so $\sigma$ fixes $a = \sin \alpha$ and also fixes $c$ since 
\begin{eqnarray*}
c \sigma &=& -\sin( 15 \cdot 9 \alpha)  \\
&=& -\sin 23 \alpha  \\
&=& \sin 5 \alpha  \\
&=& \sin 9 \alpha \\
&=& \sin \gamma \ = \ c
\end{eqnarray*}
On the other hand $b\sigma = -\sin 60 \alpha = -\sin 4\alpha = -b$.  
 
Let $U$ and $V$ be (the  lengths of) two sides of $ABC$, and let $\theta$ be the 
angle between sides $U$ and $V$.
 Then let the number of $a$, $b$, and $c$ sides of tiles 
on side $U$ 
be $p$, $q$, and $r$ respectively, and let $m$, $n$, and $\ell$ be the number of $a$, $b$, and $c$ sides on $V$.
Then we have 
\begin{eqnarray*}
U &=& pa + qb + rc \\
V &=& ma + nb + \ell c \\
U\sigma &=& pa - qb + rc \\
V\sigma &=& ma - nb + \ell c
\end{eqnarray*}
All the angles of $ABC$ are multiples of $\alpha$.  For which $J$ do we have $\sin J\alpha$ 
  fixed by $\sigma$?  Exactly for $J$ even, because 
 \begin{eqnarray*}
(2i \sin J\alpha ) \sigma &=&  (\zeta^J-\zeta^{-J})\sigma \\
&=& (\zeta \sigma)^J - (\zeta\sigma)^{-J} \\
&=& \zeta^{15J} - \zeta^{-15J} \\
&=& (-\zeta)^J - (-\zeta)^{-J} \\
&=& (-1)^J (\zeta^J - \zeta^{-J}) \\
&=& (-1)^J 2i \sin J\alpha 
\end{eqnarray*}
But $\sigma$ changes the sign of $i$, since $i\sigma = \zeta^7 \sigma = \zeta^{7 \cdot 15} = -i$.
Hence the left side is equal to $-2i(\sin J\alpha) \sigma$, and we can divide by $-2i$, obtaining
\begin{eqnarray*}
(\sin J \alpha) \sigma &=& (-1)^{J+1} \sin J \alpha
\end{eqnarray*}
So for $J$ odd, $\sigma$ fixes $\sin J \alpha$,  and for $J$ even, $\sigma$ changes
the sign of $ \sin J \alpha$.
   
For some integer $J$, angle $\theta =    J \alpha$.   
  The area equation is 
\begin{eqnarray}
2N \A_T &=& 2\A_{ABC} \nonumber \\
Nabc &=& UV \sin  \theta  \nonumber \\
  &=& (pa + qb + rc)(ma + nb + \ell c) \sin \theta  \label{eq:2034}
  \end{eqnarray}
Applying $\sigma$ we find
$$ -Nabc = (-1)^{J+1}(pa -qb + rc)(ma-nb+\ell c) (-1)^{J+1}\sin \theta$$
since $\sigma$ fixes $a$ and $c$ but changes the sign of $b$.
Dividing the last two equations, and multiplying by $-1$, we have 
\begin{eqnarray*}
1 &=& (-1)^{J} \frac {(pa + qb + rc)(ma + nb + \ell c)}
      {(pa -qb + rc)(ma-nb+\ell c)}
\end{eqnarray*}
The fraction on the right can be 1 only if the $b$ terms do not appear, 
i.e. $q = n = 0$, and it can be $-1$ only if the $a$ and $c$ terms do not 
appear, i.e. $p=r=m=\ell = 0$.  That is, it can be $1$ only if $U$ and $V$
both have no $b$ edges,  and it can be $-1$ only if $U$ and $V$ are composed
{\em only} of $b$ edges.  (If that conclusion is not obvious, it can be 
be reached by multiplying by the denominator and simplifying, remember that all 
the letters denote nonnegative quantities.)  Hence: if $J$ is even, 
$U$ and $V$ have no $b$ edges, and if $J$ is odd, they have only $b$ edges.

Consider the case in which $J$ is odd, and  $U$ and $V$ both are composed of
 only $b$ sides of tiles.   In that case the area equation becomes
 $ Nabc = qnb^2 \sin \theta$.
 Canceling $b$ we have $$Nac = qnb \sin \theta.$$ 
Now the plan is to express
  $a$, $b$, and $\sin \theta$ in terms of $\zeta$.   The minimal polynomial $\psi(x)$ 
  of $\zeta$
  is the 28-th cyclotomic polynomial, so it is known.  The degree of $\psi$ 
  is $\varphi(28) = 12$, where $\varphi$ is the Euler totient function.  
  We will get a polynomial $f(x)$   that is satisfied by $\zeta$,
  so it must be zero mod $\psi$.  But it will turn out not to be zero mod $\psi$.
 Here are the details.
First, the minimal polynomial of $\zeta$ is  
\begin{equation}
 \psi(x) =  x^{12}-x^{10} +x^8 - x^6 + x^4 - x^2 + 1 \label{eq:psi28}
 \end{equation}
as can be found using the {\tt cyclotomic\_polynomial} or the {\tt monopoly} function of Sage,
 or by factoring $x^{28} -1 $ in Mathematica or Maple, or by using techniques found in any number theory book.

Next, we express $a$, $b$, and $\sin \theta$ in terms of $\zeta$:
\begin{eqnarray*}
a &=& \sin \frac {\pi}{14} \ = \ - \frac i 2 (\zeta - \zeta^{-1} )\\
b &=& \sin \frac {4\pi}{14} \ = \  - \frac i 2 (\zeta^4 - \zeta^{-4}) \\
c &=& \sin \frac {9\pi}{14} \ = \  - \frac i 2 (\zeta^9 - \zeta^{-9}) \\
\theta &=& \sin \frac {J\pi}{14} \ = \  - \frac i 2 (\zeta^{J} - \zeta^{-J})\\
4Nac &=& - N(\zeta-\zeta^{-1})(\zeta^9-\zeta^{-9}) \\
4qnb \sin \theta &=& - qn(\zeta^4 - \zeta^{-4})(\zeta^{J} - \zeta^{-J}) \\
\end{eqnarray*}
Since $Nac = qnb \sin C$ we can equate the last two right hand sides:
\begin{eqnarray*}
  N(\zeta-\zeta^{-1})(\zeta^9-\zeta^{-9}) &=& qn(\zeta^4 - \zeta^{-4})(\zeta^{J} - \zeta^{-J}) \\
 N(\zeta^{10} -\zeta^8 - \zeta^{-8} + \zeta^{-10}) &=& qn(\zeta^{J+4} - \zeta^{J-4} - \zeta^{4-J} + \zeta^{-J-4} )
\end{eqnarray*}
Now we want to multiply by a sufficient power of $\zeta$ to kill the negative exponents;  that power 
has to be at least 10 and at least $J+4$.   Since $14\alpha = \pi$,  if $ABC$ has any angles that 
are odd multiples of $\alpha$,  then it has at least two such, and the smallest one is less than 
$\pi/2 = 7\alpha$,  so we only need to consider $J = 1$, 3, and 5 among odd $J$.  Hence $\zeta^{10}$
is a sufficient power to kill the negative exponents.  Multiplying by $\zeta^{10}$ we have
 \begin{eqnarray*}
N(\zeta^{20} - \zeta^{18} - \zeta^2 + 1) &=& qn(\zeta^{J+14} - \zeta^{J+6} -\zeta^{14-J} 
+ \zeta^{6-J}) 
\end{eqnarray*}
Since $\zeta^{14} = -1$ we can simplify this a bit:
\begin{eqnarray*}
N(-\zeta^{6} + \zeta^{4} - \zeta^2 + 1) &=& qn(-\zeta^{J} - \zeta^{J+6} -\zeta^{14-J} 
+ \zeta^{6-J})
\end{eqnarray*}
Define 
\begin{eqnarray*}
f(x) &:=& Nx^{6} - Nx^{4} - qnx^{J+6} - qnx^{J} - qnx^{14-J} + qnx^{6-J} + Nx^2 - N
\end{eqnarray*}
Then $f$ is a[ polynomial satisfied by $\zeta$.  Therefore $f(x)$ must 
be zero mod the cyclotomic polynomial $\psi$ exhibited in (\ref{eq:psi28}), which has degree $12$.
The degree of $f$ can only exceed 12 if $J+6 \ge 12$ or $14-J \ge 12$. We only have 
to consider $J = 1$, 3, or 5,  and when $J=3$ or $J=5$ we are done already, because
the degree of $f$ is less than 12 and the constant term $-N$ is not zero, so $f$ is not 
the zero polynomial.  That leaves only the case $N=1$, in which we must compute $f$ mod $\psi$.
 We did this computation in Sage as follows for $J=1$:
\begin{verbatim}
sage: R.<x,N,q,n> = PolynomialRing(QQ,4)
sage: f = N*x^6 - N*x^4 - q*n* x^7 - q*n*x^13 + q*n*x^5 + N*x^2 - N
sage: psi = x^12 - x^10 + x^8 - x^6 + x^4 - x^2 + 1
sage: f.quo_rem(psi)
(-q*n*x, -q*n*x^11 + q*n*x^9 - 2*q*n*x^7 + 
      2*q*n*x^5 + x^6*N - q*n*x^3 - x^4*N + q*n*x + x^2*N - N)
\end{verbatim}
This tells us that $f(x)$ mod $\psi(x)$ is 
\begin{eqnarray*}
-qnx^{11} + qnx^9 - 2qnx^7 + Nx^6 + 2qn x^5 - Nx^4 - qnx^3  + Nx^2  + qnx - N
\end{eqnarray*}
This is not the zero polynomial since its constant term $N$ is not zero.   Hence 
$ABC$ has no angle that is an odd multiple of $\alpha$.

Now suppose $\theta$ is an even multiple of $\alpha$.  
Since $ABC$
has no angles that are odd multiples of $\alpha$, all its angles are even multiples
of $\alpha$, so the smallest one is less than or equal to $4 \alpha$; hence we 
only have to rule out angles $\theta = 2\alpha$ and $4\alpha$.   If $\theta$ is 
an even multiple of $\alpha$, then (as shown above)  
 $U$ and $V$ have no 
$b$ edges, and the area equation (\ref{eq:2034}) becomes 
\begin{eqnarray}
Nabc    &=& (pa  + rc)(ma +   \ell c) \sin \theta   \label{eq:2157}
\end{eqnarray}
First we take up the case $J=4$.  Then $\theta = 4\alpha = \beta$,
so $\sin \theta = b$.  Then we can divide the equation by $b$, obtaining
\begin{eqnarray*}
Nac &=& (pa + rc)(ma + \ell c)
\end{eqnarray*}
The plan is the express everything in terms of $\zeta$ and then take the 
equation mod $\psi$.  We have
\begin{eqnarray*}
N(\zeta-\zeta^{-1})(\zeta^9-\zeta^{-9}) &=& (p(\zeta-\zeta^{-1}) + r(\zeta^9-\zeta^{-9})
               (m(\zeta-\zeta^{-1}) + \ell(\zeta^9-\zeta^{-9})) \\
N\zeta^{10}-N\zeta^8 - N\zeta^{-8}+ N\zeta^{-10}
&=& \zeta^{18} r\ell + \zeta^{10}(rm + \ell p) + \zeta^2pm + \zeta^{-2} pm  \\
&& + \zeta^{-10}(rm + \ell p) + \zeta^{-18} r\ell
\end{eqnarray*}
Since $\zeta^{14} = -1$ we have 
\begin{eqnarray*}
N\zeta^{10}-N\zeta^8 - N\zeta^{-8}+ N\zeta^{-10} &=& \zeta^{10}(rm + \ell p) -\zeta^4 r\ell + \zeta^2 pm \\
&& +  \zeta^{-2}pm  - \zeta^{-4} rl + \zeta^{-10}(rm + \ell p) 
\end{eqnarray*}
Bringing everything to one side of the equation we have 
\begin{eqnarray*}
0 &=& \zeta^{10}(rm + \ell p - N) + N \zeta^8  - \zeta^4 r\ell + \zeta^2 pm \\
&&+ \zeta^{-2} pm - \zeta^{-4} r\ell + N \zeta^{-8} + \zeta^{-10}(rm  + \ell p - N)
\end{eqnarray*}
Multiplying by $\zeta^{10}$ we have 
\begin{eqnarray*}
0 &=& \zeta^{20}(rm + \ell p - N) + N \zeta^{18}  - \zeta^{14} r\ell + \zeta^{12} pm + \zeta^{8} pm - \zeta^{6} r\ell + N \zeta^{2} +  (rm  + \ell p - N)
\end{eqnarray*}
Again using $\zeta^{14} = -1$ we have 
\begin{eqnarray*}
0 &=& -\zeta^{6}(rm + \ell p - N) - N \zeta^{4}  + r\ell + \zeta^{12} pm + \zeta^{8} pm - \zeta^{6} r\ell + N \zeta^{2} +  (rm  + \ell p - N) \\
&=& pm \zeta^{12} + pm \zeta^8 -\zeta^6(rm + \ell p - N + r\ell) -N\zeta^4 + N\zeta^2 + (rm + \ell p - N)
\end{eqnarray*}
This polynomial in $\zeta$ has the same degree as the minimal polynomial $\psi$, so it must
be either identically zero, or a constant multiple of $\psi$.  It is not a constant multiple, since 
it lacks a term in $\zeta^{10}$, while $\psi$ has such a term.  It is not identically zero,
since the coefficient of $\zeta^2$ is $N \neq 0$.  This contradiction shows that $ABC$
has no angle $\theta = 4\alpha$.  

 Therefore $ABC$ must have an angle $\theta = 2\alpha$.  There is no longer a canceling 
 factor in the area equation (\ref{eq:2157}), so the algebra is more tedious.  We 
reach for computer support:
\begin{verbatim}
sage: R.<x,N,p,r,m,l> = PolynomialRing(QQ,6)
sage: a = x - x^-1
sage: b = x^4-x^-4
sage: c = x^9 - x^-9
sage: f = N*a*b*c - (p*a+r*c)*(m*a + l *c)
sage: t = x^2  - x^-2
sage: f = x^20*(N*a*b*c - (p*a+r*c)*(m*a + l *c)*t)
sage: f 
-x^40*r*l + x^36*r*l + x^34*N - x^32*r*m - x^32*p*l - x^32*N 
+ x^30*r*m + x^30*p*l + x^28*r*m + x^28*p*l - x^26*r*m - x^26*p*l
 - x^26*N - x^24*p*m + x^24*N + 2*x^22*p*m + 2*x^22*r*l - 2*x^18*p*m 
 - 2*x^18*r*l + x^16*p*m - x^16*N + x^14*r*m + x^14*p*l + x^14*N 
 - x^12*r*m - x^12*p*l - x^10*r*m - x^10*p*l + x^8*r*m + x^8*p*l 
 + x^8*N - x^6*N - x^4*r*l + r*l
\end{verbatim}
That is, 
\begin{eqnarray*}
f(x) &=& -r\ell x^{40} + rl x^{36} + x^{34} N - x^{32} (rm-p\ell- N) +x^{30}(rm+p\ell) + x^{28}(rm+p\ell) \\
&&+ x^{26} (rm-p\ell-N) + x^{24}(N-pm) +2x^{22}(pm + r\ell) - 2x^{18}(pm + r\ell) + x^{16}(pm-N) \\
&&+ x^{14}(rm + p\ell+N) - x^{12}(rm + p\ell) - x^{10}(rm + p\ell) + x^8(rm + p\ell + N) \\
&& - Nx^6 - r\ell x^4 + r\ell
\end{eqnarray*}
Since $x^{14} = -1$ and $x^{28} = 1$ we have
\begin{eqnarray*}
f(x) &=& =rl x^{12} + rlx^8 + x^6 N - x^4(rm-p\ell-N) + x^2(rm + p\ell) + (rm+p\ell) \\
&& -x^{12}(rm-p\ell-N) - x^{10}(N-pm) - 2x^8(pm + r\ell) + 2x^4(pm + r\ell) - x^2(pm-N) \\
&& +(rm + p\ell-N) - x^{12} (rm + pl)  - x^{10}(rm + p\ell) + x^8(rm + p\ell + N) \\
&& - Nx^6 - r\ell x^4 + r\ell
\end{eqnarray*}
Collecting like terms we have
\begin{eqnarray*}
f(x) &=& x^{12}(r\ell - 2rm   + N) + x^{10}(pm-N -rm - p\ell) + x^8(-2pm -r\ell + rm + p\ell +N) \\
&& x^4(-rm+p\ell + N+2pm+r\ell) + x^2(rm+p\ell-pm + N) + (rm+p\ell- N + r\ell)
\end{eqnarray*}
This polynomial must vanish at $\zeta$, and it has degree at most 12, so either it is 
identically zero or it is a multiple of $\psi$.  Since $f$ has a zero coefficient of $x^6$,
it is not a multiple of $\psi$; hence it must be identically zero.   From the constant 
term we have $N = rm + p\ell + r\ell$;  substituting that value for $N$ into the 
coefficient of $x^4$ we have 
\begin{eqnarray*}
0 &=& -rm  + p\ell + N + 2pm + r\ell \\
&=& -rm + p\ell + (rm + p\ell + r\ell) + 2pm + r\ell\\
&=& 2p\ell +r\ell + 2pm + r\ell \\
&=& 2\ell(p+r) + 2pm 
\end{eqnarray*}
Hence $pm = 0$ and $\ell(p+r) = 0$.
Hence either $\ell =0 $ or $p=r=0$.
   If $p=r=0$ then the constant term is $-N \neq 0$,
so we are finished.  Hence we can assume $\ell = 0$.   Then the constant term is $rm-N$,
so $N= rm$.  Hence $m \neq 0$; but since $pm = 0$ we have $p=0$.  Then the coefficient 
of $x^2$ is $rm + p\ell -\pm + N = rm +N$, but since $rm = N$ that is $2N$, which is 
not zero.  That contradiction completes the proof of the lemma.
\medskip

\section{The case when $3\alpha + 2\beta = \pi$ or $2\alpha + 3\beta = \pi$}

For this section only, we drop the assumption that $\alpha \le \beta$, assuming only that $\alpha \le \gamma$ and 
$\beta \le \gamma$.
We consider tilings of a triangle $ABC$ that is not isosceles,  by a tile with angles
$\alpha$, $\beta$, and $\gamma$, where $3\alpha + 2\beta = \pi$ and $\alpha$ is not a rational multiple of $\pi$.
 Note that we must have $\gamma = \beta + 2\alpha$,
since $\gamma = \pi-\alpha-\beta = 3\alpha+2\beta - \alpha - \beta = 2\alpha + \beta$. The possible shapes of $ABC$ are quite
limited:  

\begin{lemma} \label{lemma:twoshapes}
Let $ABC$ be $N$-tiled by a tile with angles $\alpha$, $\beta$, and $\gamma$, and suppose $3\alpha + 2\beta = \pi$.
Suppose $ABC$ is not similar to the tile and not isosceles, and suppose that $\alpha$ is not a rational multiple of $\pi$.
Then either $ABC$ has angles $2\alpha$, $\beta$, and $\beta + \alpha$,
or $ABC$ has angles $\alpha$, $2\alpha$, and $2 \beta$.
\end{lemma}

\noindent{\em Proof}.  Note that $\beta$ is not a rational multiple of $\alpha$, since then the relation $3\alpha + 2\beta = \pi$
would make $\alpha$ a rational multiple of $\pi$; similarly $\beta$ is not a rational multiple of $\pi$.
Suppose first that $ABC$ does not have an $\alpha$ angle.   If one angle (say $A$) is $2\alpha$,  then the tiling 
must split that angle into two $\alpha$ angles, since $\beta \neq 2\alpha$.   Suppose, for proof by contradiction,
that the tiing split one of the other two vertex angles of $ABC$ into several $\alpha$ angles.  It must be 
least 3 of them since $ABC$ is not isosceles.  Hence at least five $\alpha$ angles are involved at the vertices of $ABC$.
Hence the equation that says the sum of the vertex angles is $\pi$ is an integral relation of the form 
$n \alpha + m \beta + r \gamma = \pi$, with $n \ge 5$.   Since $\gamma = \beta + 2\alpha$ we have 
$(n+2r) \alpha + (m+r) \beta = \pi$.  Since $n+2r \ge 5$ this relation is not a multiple of $3\alpha + 2\beta = \pi$,
and hence the two equations can be solved for $\alpha$ and $\beta$, which will be rational multiples of $\pi$.  This 
contradicts the hypothesis, and the contradiction shows that $ABC$ does not have an angle that splits into $\alpha$ angles.
As a consequence of this observation, $ABC$ does not contain an angle larger than $\gamma - \alpha = \beta + \alpha$,
since if it did, then the third angle would be less than $\beta$ and so would have to split into $\alpha$ angles.  
Hence the other two angles of $ABC$ must be at least $\beta$ and at most $\beta + \alpha$.   There is only only 
possibility, namely the first triangle mentioned in the lemma.

Now suppose that $ABC$ does have an $\alpha$ angle.  If the second smallest angle is $2\alpha$, then the remaining 
angle must be $2\beta$;  that is the second possibility mentioned in the lemma.  If the second smallest angle is 
greater than $2\alpha$, it must be at least $\beta$, since if not it will split into at least three $\alpha$ angles,
giving rise to an integral relation that is not a multiple of $3\alpha + 2\beta = \pi$.  It cannot be exactly $\beta$
as that would make $ABC$ similar to the tile.  Therefore the second angle is more than $\beta$, and the third 
angle is therefore less than $\gamma$.  If neither of these two angles splits into $\alpha$ angles,  then each one 
contains a $\beta$ angle plus one or more $\alpha$ angles. But since $ABC$ is not isosceles, they cannot both 
be $\beta + \alpha$; hence at least four  $\alpha$ angles are involved in the tiling at the vertices of $ABC$.  That 
is contradictory, since $4\alpha + 2\beta > \pi$.  Hence one of the second two angles does split into $\alpha$ angles.
Since the second smallest angle is more than $2\alpha$, at least 3 $\alpha$ angles are required, and 
that gives rise to an integer relation that is 
not a multiple of $3\alpha + 2\beta$, contradiction.  That completes the proof of the lemma.

We next introduce notation and definitions that will be used throughout this section, i.e.,  regardless of the 
shape of $ABC$.
As always,  $a$, $b$, and $c$ are the sides of the tile, and we assume the tile and the triangle $ABC$ are 
scaled so that 
\begin{eqnarray*}
a &=& \sin \alpha \\
b &=& \sin \beta \\
c &=& \sin \gamma
\end{eqnarray*}
For use in this entire section, we introduce the following notation:
\begin{eqnarray*}
z &=& e^{i\pi\alpha/2} \\
\zeta &=& e^{i\pi\beta}
\end{eqnarray*}
We note that 
$$z^6\zeta^2 = -1$$
since $3\alpha + 2\beta = \pi$. Hence $\zeta = iz^{-3}$.  We have
\begin{eqnarray*}
2ia &=& z^2 - z^{-2} \\
2ib &=& \zeta- \zeta^{-1} \\
&=& i(z^3 + z^{-3}) \\
2b &=& z^3 + z^{-3} \\
2\sin 2\alpha &=& -i(z^4-z^{-4})\\
2ic &=& e^{i(2\alpha + \beta)} - e^{-i   (2\alpha + \beta)} \\
&=& z^4 \zeta - z^{-4}\zeta^{-1} \\
&=& iz^4 z^{-3} + iz^{-4}z^3  \mbox{\quad since $\zeta = iz^{-3}$}\\
&=& i(z + z^{-1})\\
2c &=& z + z^{-1}
\end{eqnarray*}
In other words,  $\sin \gamma = \cos (\alpha/2)$.   We have
\begin{eqnarray*} 
\sin(\beta + \alpha) &=& \sin(\pi-\gamma) \\
&=& \sin \gamma \\ 
&=&  \frac 1 2 (z + z^{-1})
\end{eqnarray*}
We note the factorizations
\begin{eqnarray*}
z^2 - z^{-2} &=& (z+z^{-1})(z-z^{-1}) \\
z^3 + z^{-3} &=& (z+z^{-1})(z^2 - 1 + z^{-2}) 
\end{eqnarray*}
In view of these factorizations we have
\begin{eqnarray*}
\frac a c &=& =-i(z-z^{-1}) \\
\frac b c &=& z^2 - 1 + z^{-2} 
\end{eqnarray*}
We introduce 
\begin{eqnarray*}
t &:=& 2\cos \alpha  \\
&=& z^2 + z^{-2} \\
s &:=& 2\sin(\alpha/2) \\
  &=& -i(z-z^{-1})
\end{eqnarray*}
Then $t = 2-s^2$, and 
\begin{eqnarray*}
\frac a c &=& =i(z-z^{-1}) \\
&=& s \\
\frac b c &=& z^2 - 1 + z^{-2} \\
&=& t-1 \\
&=& 1-s^2 \\
\sin 2\alpha &=& 2\sin\alpha\cos \alpha \\
&=& at \\
&=& a(2-s^2)  
\end{eqnarray*}
The $\d$ matrix equation is
\begin{eqnarray*}
\d \vector a b c &=& \vector X Y Z
\end{eqnarray*}
 We introduce individual names for the components of the $\d$ matrix to avoid so many 
subscripts:
\begin{eqnarray*}
\d &=& \matrix p d e g m f h \ell r 
\end{eqnarray*}
Therefore we have for the three sides  $X$, $Y$, and $Z$ of triangle $ABC$ 
\begin{eqnarray*}
X &=& pa + db + ec \\
Y &=& ga + mb + fc \\
Z &=& ha + \ell b + rc
\end{eqnarray*}

Generally our convention is that $X$, $Y$, and $Z$ are in order of size, so $X$ is opposite
the smallest angle $A$, and $Y$ is opposite the middle angle $B$, and $Z$ is opposite the 
largest angle $C$.  In this section, however, we are not assuming $\alpha < \beta$, 
and that means that we may not have $Y < Z$ either.  In case triangle $ABC$ has angles $\alpha$, $2\alpha$, and $2\beta$, 
we will assume that $X$ is opposite the $\alpha$ angle at $A$, and the $2\alpha$ angle is at $B$, opposite $Y$, and 
that $Z$ is opposite the $2\beta$ angle.    In case $ABC$ has angles $2\alpha$, $\beta$, and $\beta + \alpha$, we will
assume $X$ is opposite the $2\alpha$ angle at $A$, and $Y$ is opposite the $\beta$ angle at $B$, and $Z$ is 
opposite the $\beta + \alpha$ angle at $C$.   Thus in case $\beta < \alpha$,  we may not have $X$, $Y$, $Z$ in 
order of size.

\begin{lemma} \label{lemma:zerolimits1} Suppose triangle $ABC$ is $N$-tiled by a tile in which $\gamma > \pi/2$
(as is the case when $3\alpha + 2\beta = \pi$).  Suppose all the tiles along one side of $ABC$ do not have their $c$ sides
along that side of $ABC$.  Then there is a tile with a $\gamma$ angle at one of the endpoints of that side of $ABC$.
\end{lemma}

\noindent{\em Proof}.  Let $PQ$ be the side of $ABC$ with no $c$ sides of tiles along it.  Then the 
$\gamma$ angle of each of those tiles   occurs at a vertex on $PQ$, since the angle opposite the side 
of the tile on $PQ$ must be $\alpha$ or $\beta$.  Let $n$ be the number of tiles along $PQ$;  then there 
are $n-1$ vertices of these tiles on the interior of $PQ$.  Since $\gamma > \pi/2$,  no one vertex has more 
than one $\gamma$ angle.  By the pigeonhole principle, there is at least one tile whose $\gamma$ angle is 
not at one of those $n-1$ interior vertices;  that angle must be at $P$ or $Q$.  That completes the proof of the lemma.

\begin{lemma} \label{lemma:zerolimits2}
 Suppose $3\alpha + 2\beta = \pi$, and suppose triangle $ABC$  has angles $2\alpha$, $\beta$,
and $\beta + \alpha$.  Let $T$ be a triangle with angles $\alpha$, $\beta$, and $\gamma$ and suppose that $ABC$ 
is $N$-tiled by $T$ in such a way that just five tiles share vertices of $ABC$.  Then the following restrictions
on the elements of the $\d$ matrix apply:  We have $e \neq 0$, $f \neq 0$, and $r \neq 0$. In other words,
the third element in each row is nonzero.  
\end{lemma}

\noindent{\em Proof}.  None of the angles of $ABC$ is large enough to accommodate the 
$\gamma$ vertex of a tile, since those angles are $2 \alpha$, $\beta$, and $\beta + \alpha < \beta + 2\alpha = \gamma$.
The lemma  then follows from Lemma \ref{lemma:zerolimits1}, since the three numbers $e$, $f$, and $r$ are the 
numbers of $c$ sides of tiles on the three sides of $ABC$.   That completes the proof of the lemma.

 \begin{lemma} \label{lemma:zerolimits}
Let $3\alpha + 2\beta = \pi$, and let $ABC$ have angles $\alpha$, $2\alpha$, and $2\beta$, 
and suppose there is an $N$-tiling of $ABC$ by a tile with angles $\alpha$, $\beta$,
and $\gamma$ such that exactly two tiles meet at the $2\beta$ vertex of $ABC$.
  Then  we cannot have  $m=\ell = 0$, 
i.e. the bottom two entries in the second  column cannot both be zero. 
Also $e \neq 0$ and $f \neq 0$ and $r \neq 0$.
\end{lemma}

\noindent{Proof}.  At vertex $A$, where the angle is $\alpha$, there is just one tile.
It must have its $b$ side on $AB$ and its $c$ side on $AC$, or vice-versa.  The number of 
$b$ sides on $AC$ is $m$ and the number of $b$ sides on $AB$ is $\ell$, so they cannot 
both be zero. 

In triangle $ABC$, no $\gamma$ angles of tiles can occur at any vertex, since the $\alpha$ and $2\alpha$ angles 
are too small and the $2\beta$ angle splits into two $\beta$ angles.  Hence the third column of the $\d$ matrix
cannot contain any zero entries, by Lemma \ref{lemma:zerolimits1}.  Hence $e$, $f$, and $r$ are each nonzero.
  That completes the proof of the lemma.

\subsection{The case when $ABC$ has angles $\alpha$, $2\alpha$, and $2\beta$}

\begin{lemma} \label{lemma:3alpha2betaCase2}
Let $3\alpha + 2\beta = \pi$, and let $ABC$ have angles $\alpha$, $2\alpha$, and $2\beta$,
and let $N$ be arbitrary.  Then there is no $N$-tiling of $ABC$ by a tile with angles $\alpha$, $\beta$,
and $\gamma$ in which exactly two tiles meet at the $2\beta$ angle of $ABC$.  (Here we do not assume $\alpha < \beta$.)
\end{lemma}

\noindent{\em Proof}.
Let $X$, $Y$, and $Z$ be the sides of triangle $ABC$ opposite $A$, $B$, and $C$ respectively.
Let $\lambda$ be the ratio of the sides to the opposite angle, which (by the law of sines) is 
the same for all three sides:
\begin{eqnarray*}
X &=& \lambda \sin \alpha \\
Y &=& \lambda \sin 2\alpha \\
Z &=& \lambda \sin2 \beta 
\end{eqnarray*}
We note that
\begin{eqnarray*}
\sin \alpha &=& a \\
\sin 2\alpha &=& 2 \sin \alpha \cos \alpha \\
&=& t a \\
&=& (2-s^2)a \\
\sin 2\beta &=& \sin 3 \alpha \\
&=& 4\sin^3 - 3\sin^2 \alpha \\
&=& 4a^3-3a^2
\end{eqnarray*}
The $\d$ matrix equation is 
\begin{equation} \label{eq:dmatrix3alpha2betacase2}
 \matrix p d e g m f h \ell r \vector a b c =  \vector X Y Z = \lambda \vector a {\sin 2\alpha} {\sin 2 \beta} 
 \end{equation}
As noted above, we have
$\sin 2 \beta = 4a^3-3a^2$ and $\sin 2\alpha = a(2-s^2)$
we have
$$ \matrix p d e g m f h \ell r \vector a b c = \lambda \vector a {a(2-s^2)} {4a^3-3a^2} $$
Taking the ratio of the second row to the first row, we have 
\begin{eqnarray*}
\frac {ga + mb + fc} {pa + db + ec} &=& 2-s^2
\end{eqnarray*}
Multiplying by the denominator, we have 
\begin{eqnarray*}
ga + mb + fc &=& (2-s^2)(pa + db + ec)
\end{eqnarray*}
Dividing both sides by $c$ and expressing $a/c$ and $b/c$ in terms of $s$ we have
\begin{eqnarray*}
g \frac a c  + m \frac b c  + f &=& (2-s^2)(p \frac a c + d \frac b c + e) \\
 gs + m(1-s^2) + f &=& (2-s^2)(ps + d(1-s^2) + e)
\end{eqnarray*}
This is a fourth-degree polynomial equation for $s$ over $\Q$.   
Bringing the equation to polynomial form we have 
\begin{eqnarray}
0 &=& \psi(s)\  :=  \ ds^4 - ps^3 + (m-3d-e)s^2 + (2p-g)s + (2d+2e-f-m) \label{eq:psidef32}
\end{eqnarray} 
If $\psi$ is identically zero,
   then from the first two 
coefficients of $\psi$ in (\ref{eq:psidef32}) we have $d = p = 0$.  From the coefficient of $s^2$ we have $m-3d-e = m-e = 0$, so $m=e$.
From the coefficient of $s$ we have $2p-g =0$, and since $p=0$ we have $g=0$.  From the constant 
coefficient we have $2d + 2e-f - m = 2e-f-m = e-f$, so $f=e$.  The $\d$ matrix then has the form 
$$ \matrix 0 0 e 0 e e h \ell r $$ 
This is not immediately contradictory, and we shall return to proving that $\psi$ is not identically zero below.
 
From the first row of the $\d$ matrix equation we have 
$$ pa + d b + ec = \lambda a$$
Dividing by $c$ and using $a/c = s$ and $b/c = 1- s^2$ we have 
$$ \lambda s = ps + d(1-s^2) + e  
$$
Solving for $\lambda$ we have 
\begin{equation}
 \lambda = p -ds +   \frac {e+d} s \label{eq:lambdaofs}
\end{equation}

The $d$ matrix equation (\ref{eq:dmatrix3alpha2betacase2}) can be written as an eigenvalue problem 
this way:
$$ \matrix {p} d e 
           {\frac {gb} {\sin 2\alpha}} {\frac {mb} {\sin 2\alpha}} {\frac {fb} {\sin 2\alpha}}
           {\frac {hc} {\sin 3\alpha}} {\frac {\ell c} {\sin 3\alpha}} {\frac {rc} {\sin 3\alpha}}
    \vector a b c = \lambda \vector a b c 
$$
To find an eigenvector by the cofactor method, we need to compute the cofactors of the matrix
$$ \matrix {p-\lambda} d e 
           {\frac {gb} {\sin 2\alpha}} {\frac {mb} {\sin 2\alpha} - \lambda} {\frac {fb} {\sin 2\alpha}}
           {\frac {hc} {\sin 3\alpha}} {\frac {\ell c}{\sin 3\alpha}} {\frac {rc} {\sin 3\alpha}-\lambda}
$$
Taking the cofactors of the third row, and multiplying by $\sin 2\alpha$, we find a candidate for an eigenvector
(it is only a candidate until we prove that its components are all nonzero):
\begin{eqnarray} 
\vector u v w  &=& \vector {(df-em)b+ e \lambda \sin 2\alpha} 
                           {(eg-pf)b+ f\lambda b}
                           {(pm - dg)b - (mb + p \sin 2 \alpha)\lambda  +  \lambda^2\sin 2\alpha}  \label{eq:74}
\end{eqnarray}
Now assume, for proof by contradiction, that $\psi$ is identically zero.
As shown above, the $\d$ matrix then has the form 
$$ \matrix 0 0 e 0 e e h \ell r $$ 
Then we have
$$ \vector u v w = \vector {-be^2  +  e \lambda \sin 2\alpha}   { eb \lambda} {eb \lambda + \lambda^2 \sin 2\alpha}
$$
Since $e \neq 0$ (otherwise the whole first row of the $\d$ matrix is zero), the second two components are nonzero.
Suppose, for proof by contradiction, that $u = 0$.  Then since $e \neq 0$ we have
\begin{eqnarray*}
eb &=& \lambda \sin 2\alpha
\end{eqnarray*}
Since $\lambda \sin 2\alpha =  ga + mb + fc = eb + ec$, we have $eb = eb + ec$, and hence $ec = 0$.
Since $e \neq 0$ we have $c=0$, a contradiction. 
This contradiction shows that $u \neq 0$ (still under the assumption that $\psi$ is identically zero).  
  Having proved $u \neq 0$, we conclude that  the eigenspace has dimension 1, and 
  $(u,v,w)$ is a multiple of $(a,b,c)$.   Therefore $v/w = b/c$.   Cross multiplying, we have $vc = bw$.
Putting in $v = eb\lambda$ and $w = eb\lambda + \lambda^2 \sin 2\alpha$ we obtain
\begin{eqnarray*}
eb\lambda c &=& eb^2\lambda + b\lambda^2 \sin 2\alpha \\
            &=& eb^2\lambda + b\lambda^2 at \mbox{\qquad since $\sin 2\alpha = at)$}
\end{eqnarray*}
Dividing by $bc\lambda $ we have
\begin{eqnarray*}
e &=& e\frac b c + \lambda t \frac a c 
\end{eqnarray*}
Since $p=d=0$, the first row of the $\d$ matrix equation is $ec = \lambda a$.  Hence $\lambda = ec/a$.  Putting
this value in for $\lambda$ we have
\begin{eqnarray*}
e &=& e \frac b c +\frac {ec} a t \frac a c \\
&=& e \frac b c + et \\
&=& e(t-1) + et \mbox{\qquad since $b/c = t-1$} \\
\end{eqnarray*}
Since $e \neq 0$ we can divide by $e$:
\begin{eqnarray*}
1 &=& t-1 + t \\
&=& 2t-1 \\
\end{eqnarray*}
Solving for $t$ we find $t=1$.  Since $t=2\cos \alpha$ we have $\cos \alpha = 1/2$, so $\alpha = \pi/3$. Then
$\beta = \pi/2- 3\alpha/2 =0$,   contradiction.
Note that this is a contradiction even if we do not assume $\alpha \le \beta$.
This contradiction shows that $\psi$ is not identically zero.

We next want to prove that the components of $(u,v,w)$ are not zero,  without the assumption that $\psi$ is identically 
zero.  First assume, for proof by contradiction, that $v=0$.   Then 
\begin{eqnarray*}
(eg-pf)b + f \lambda b &=& 0 
\end{eqnarray*}
Dividing both sides by $b$ we have
\begin{eqnarray*}
eg-pf + f \lambda &=& 0
\end{eqnarray*}
Putting in the value of $\lambda$ from (\ref{eq:lambdaofs}) we have
\begin{eqnarray*}
0 &=& eg-pf + f\big(p-ds + \frac{e+d} s\big) \\
&=& (eg-pf)s + f(ps-ds^2 + e+d) 
\end{eqnarray*}
Collecting terms and changing the sign, we have 
\begin{eqnarray*}
0&=& dfs^2 -eg s - (e+d)f
\end{eqnarray*}
Let $H(s) = dfs^2 - egs - (e+d)f$.  Then $H(0) = (e+d)f \le 0$, and $H(1) = df-eg-(e+d)f = -eg-ef \le 0$.
Since $H$ is a quadratic polynomial with positive leading coefficient, it cannot have a zero between 0 and 1.
But $s$ is such a zero, contradiction.  Hence $v \neq 0$.

Next assume, for proof by contradiction, that $u=0$.   Then 
\begin{eqnarray*}
0 &=& (df-em) b + e \lambda \sin 2\alpha  \\
&=& (df-em) b + e \lambda at  \mbox{\qquad since $\sin 2\alpha = at$}
\end{eqnarray*}
Dividing both sides by $c$ and using (\ref{eq:lambdaofs}) we have
\begin{eqnarray*}
0 &=& (df - em)\frac b c + e \big(p-ds + \frac {e+d} s\big)t \frac a c 
\end{eqnarray*}
Using $b/c = 1-s^2$ and $t = 2-s^2$  and $a/c = s$, we have
\begin{eqnarray*}
0 &=& (df-em)(1-s^2) + e\big(p-ds + \frac {e+d} s\big) s(2-s^2) \\
&=&  (df-em)(1-s^2) + e(ps-ds^2 + e+d)(2-s^2) \\
&=& eds^4 - eps^3 + s^2(em-df -3ed -e^2 ) +2eps + (df-em + 2e^2 + 2de)
\end{eqnarray*}
By Lemma \ref{lemma:zerolimits1}, we have $e \neq 0$, since there cannot be a tile with a $\gamma$ angle at $B$ (where the 
$2\alpha$ angle of $ABC$ is),  and there cannot be a tile with a $\gamma$ angle at $C$, where the $2\beta$ angle of $ABC$ is, 
since by hypothesis, there are two tiles each with a $\beta$ angle at that vertex.

Define
$$H(s): = ds^4  - ps^3 + s^2\big(m- \frac{df} e -3d -e \big) +2ps + \big(\frac{df} e- m + 2e + 2d\big)$$
Since $e \neq 0$, we can divide the previous equation by $e$, obtaining  $H(s) = 0$:
\begin{equation} \label{eq:May1-2}
0= ds^4  - ps^3 + s^2\big(m- \frac{df} e -3d -e \big) +2ps + \big(\frac{df} e- m + 2e + 2d\big)
\end{equation}
 The first two terms of $H$ are equal to the first two terms of 
$\psi$, so $H-\psi$ is a quadratic in $s$:
\begin{eqnarray*}
0 &=& H(s) - \psi(s) \\
&=& s^2\big(m-\frac{df} e - 3d-e -(m-3d-e)\big) + s(2p-(2p-g))\\
&& + \big(\frac {df} e - m + 2e+2d\big) -(2d+2e-f-m)  
\end{eqnarray*}
\begin{eqnarray}
&=& -s^2\big(\frac{df} e\big) +gs + \frac{df} e  \label{eq:May1-1}
\end{eqnarray}
We have  $f \neq 0$, by Lemma \ref{lemma:zerolimits}.
 Assume, for proof by contradiction, that $d \neq 0$.  Then we can divide by $-df/e$, obtaining
\begin{eqnarray*}
0 &=& s^2  - \frac {eg}{df} s - 1 
\end{eqnarray*}
This quadratic function is negative when $s=0$, and negative when $s=1$, and its leading coefficient 
is positive. Therefore it has no zero between 0 and 1, contradiction.   This contradiction shows that $d=0$.

Since $d=0$,  equation (\ref{eq:May1-1}) becomes $gs = 0$.  Hence $g=0$, since $s\neq 0$.
The equation (\ref{eq:psidef32}) ($\psi(s) = 0$)  then becomes
\begin{eqnarray*}
0 &=& - ps^3 + (m-e)s^2 + 2ps + (2e-f-m) 
\end{eqnarray*}
On the other hand, the equation (\ref{eq:May1-2}) ($H(s) = 0$)  now becomes (with $d=g=0$)
$$ -ps^3 + (m-e)s^2 + 2ps + (2e-m) $$
Subtracting that from the previous equation we find $f=0$; but we already proved $f \neq 0$ using Lemma \ref{lemma:zerolimits1}.
That contradiction depended on the assumption $u=0$, and that completes the proof by contradiction that $u\neq 0$.

Now assume, for proof by contradiction, that $w=0$.   That is, 
\begin{eqnarray*}
0 &=& (pm-dg)b - (mb + p \sin 2\alpha) \lambda + \lambda^2 \sin 2\alpha
\end{eqnarray*}
To write this as a function of $s$, use  $\sin 2\alpha = at = cs(2-s^2)$ and $b = c(1-s^2)$
and $\lambda = p-ds + (e+d)/s$, and then multiply by $s/c$.  We find (with the aid of a computer algebra system)
the following polynomial equation:
\begin{eqnarray*}
0 &=& -d^2 s^7 + dp s^5  +(4d^2 + 2de-dm)s^4 + (dg-3dp-ep)s^3 \\
&& (-5d^2 - 6de - e^2 + 2dm + em) s^2 \\
&&+ (-dg + 2dp + 2ep) s + (2d^2+4de+2e^2 - dm - em)
\end{eqnarray*}
Computing the polynomial remainder of this on division by $\psi$, and changing the sign, we find
\begin{eqnarray*}
0 &=&  G(s) \ := \  dfs^2 - egs - f(d+e)
\end{eqnarray*}
Assume, for proof by contradiction, that $d=0$.  Then $0 = -egs - fe$.  
By Lemma \ref{lemma:zerolimits}, $e \neq 0$ and $f \neq 0$, so $fe > 0$; since $egs \ge 0$ this is a contradiction.
Hence $d \neq 0$.  Since $f \neq 0$ as just noted,  the leading coefficient $df$ is not zero.
Thus $G$ is a quadratic function such that $G(0) = -f(d+e) < 0$.   We have $G(1) = df-eg-f(d+e) = -eg-ef < 0$.
Since $G(0)$ and $G(1)$ are both negative, and $G^{\prime\prime}$ is positive,  $G$ has no zero between 0 and 1.
This is a contradiction, since $G(s) = 0$ for $s = 2 \sin(\alpha/2)$, which is between 0 and 1 since $\alpha < \pi/3$.
That completes the proof by contradiction that $w \neq 0$.   Hence, the eigenspace is one-dimensional, and 
$(u,v,w)$ is a multiple of the eigenvector $(a,b,c)$.

Returning to equation (\ref{eq:74}), now that we know $(u,v,w)$ is a multiple of $(a,b,c)$, we have 
we have $a/c = u/w$.  That is,
\begin{eqnarray*}
\frac a c &=& \frac { (df-em)g + e \lambda \sin 2\alpha}{(pm-dg)b - (mb + p \sin 2\alpha)\lambda + \lambda^2 \sin 2 \alpha}
\end{eqnarray*}
We have $\sin 2\alpha = at = a(2-s^2)$.   Putting that in, and dividing numerator and denominator on the right by $c$,
and then using $b/c = 1-s^2$ and $a/c = s$, we have 
\begin{eqnarray*}
s &=& \frac {(df-em)(1-s^2) + es\lambda(2-s^2)}
         {(pm-dg)(1-s^2) - (m(1-s^2) + ps(1-s^2))\lambda + \lambda^2 s(1-s^2)} \\
  &=& \frac {(df-em)(1-s^2) + es\lambda(2-s^2)}
         {(1-s^2)(pm-dg- (m + ps)\lambda + \lambda^2 s} 
\end{eqnarray*}
Multiplying by the denominator we have 
\begin{eqnarray*}
s(1-s^2)(pm-dg-  (m + ps)\lambda + \lambda^2 s) &=& (df-em)(1-s^2) + es\lambda(2-s^2) 
\end{eqnarray*}
Putting in $\lambda = p - es + (e+d)/s$ and bringing this to polynomial form (using Mathematica, for example,
or by hand if you wish) we have
\begin{eqnarray*}
0 &=& d^2 s^6 - dps^5 +(dm-de-3d^2  )s^4 +(2dp-dg)s^3 \\
&&+(3d^2+de-df-2dm)s^2 + (dg-dp+ep)s +(d^2-e^2 -df-dm) 
\end{eqnarray*}
Now taking this polynomial mod $\psi$, using the {\em PolynomialRemainder} command of Mathematica  (you won't want 
to do that polynomial division by hand) we find a simple quadratic equation: 
\begin{eqnarray}
0 &=& d^2 - e^2 - df - dm - (dg - dp + ep)s + (dm + de-d^2)s^2  \label{eq:simplequadratic}
\end{eqnarray}

By Lemma \ref{lemma:zerolimits},  $e \neq 0$. 
 Suppose, for proof by contradiction, that $d=0$.  Then 
 (\ref{eq:simplequadratic})  becomes 
 $ eps + e^2 = 0$.  Since $e \neq 0$ then $ps + e = 0$.   Since $e$, $p$, and $s$ are all $\ge 0$ and 
 $e > 0$, this is a contradiction.   Hence $d \neq 0$.
 
Since $s = -i(z+z^{-1}) = 2\sin(\alpha/2)$, we have
  \begin{eqnarray*}
 \sin \frac \alpha 2 &<& \frac 1 2 \mbox{\qquad since $\alpha < \pi/3$}\\
 s &<& 1  \mbox{\qquad since $s = 2\sin(\alpha/2)$}
 \end{eqnarray*}
  Then, an $N$-tiling of $ABC$ as in the lemma gives rise to
 a  solution $s= 2\sin(\alpha/2)$ of the following two equations (which are  (\ref{eq:psidef32}) and (\ref{eq:simplequadratic})):
  \begin{eqnarray*}
 0 &=& \psi(s) \ = \    \ ds^4 - ps^3 + (m-3d-e)s^2 + (2p-g)s + (2d+2e-f-m) \\
 0 &=& \chi(s) \ := \ (dm+de-d^2)s^2 - (dg-dp+ep)s + d^2-e^2 - df-dm  
 \end{eqnarray*}
 We will finish the proof 
 by showing that these two equations have no simultaneous solution $s$ between 0 and 1.

 We compute a somewhat simpler quartic than $\psi$ satisfied by $s$, as follows:
 \begin{eqnarray*}
 F(s)\ := \ d\psi-\chi &=& d^2s^4-dps^3 -2d(d+e)s^2  + p(d+e)s + (d+e)^2
 \end{eqnarray*}
$F$ has  a positive constant coefficient (as well as positive leading coefficient) since $d > 0$.
Assume, for proof by contradiction, that $p=0$.  Then 
\begin{eqnarray*}
F(s) &=& d^2 s^4 - 2d(d+e)s^2 + (d+e))2 \\
&=& (ds-(d+e))^2
\end{eqnarray*}
Hence the only zero of $F$ is $s = (d+e)/d \ge 1$.  But (if there is an $N$-tiling),  $F$ has a zero 
between 0 and 1, contradiction.  Hence $p \neq 0$.

Now we compute the value of $F(1)$.
\begin{eqnarray*}
F(1) &=& d^2 - dp -2d(d+e) + p(d+e) + (d+e)^2 \\
&=& d^2 - dp - 2d^2 -2de + pd + pe + d^2 + 2de + e^2 \\
&=&   pe + e^2 \\
F(1) &>& 0 \mbox{\qquad since $p >0$ and $e>0$ }
\end{eqnarray*}
  We have $F(0) = (d+e)^2$. This is positive since $d > 0$.  
Now $F$ is a quartic that is positive at  0 and positive at 1, and its derivative at 0 is $F^\prime(0) = p(d+e)$, which is 
  positive since $p \neq 0$ and $d \neq 0$.  The derivative $F^\prime(s)$ is a cubic;  one of 
its zeroes is negative, since for large negative $s$ we have $F^\prime(s) < 0$.  It therefore has 
at most two zeroes between 0 and 1.  The only way that $F(s)$ can have a zero strictly between 0 and 1 is 
if $F$ first increases (as $s$ increases from 0), then reaches a maximum and decreases, 
crossing the $s$-axis (or possibly just touching the $s$-axis, if $F$ has a double zero),
then reaches a minimum, then increases (crossing the $s$-axis again,
unless it has a double zero) 
and reaches its positive value at $s=1$.  Then $F^\prime$ has two zeroes between 0 and 1, 
and in particular $F^\prime(1)$ must be positive. 
We compute 
\begin{eqnarray*}
F^\prime(s) &=& 4d^2s^3 - 3dps^2 - 4d(d+e)s + p(d+e)
\end{eqnarray*}

The product of its three roots is thus $-p(d+e)/4d^2$, which is negative.  Because $F^\prime(0)$ is 
positive, and $F^\prime(s)$ is negative for large negative $s$,  one of the three roots is negative.
But (if there is an $N$-tiling),  the other two are between 0 and 1, so the product of the three 
roots is negative.   

The sum of the three roots of $F^\prime$ is $4d(d+e)/4d^2 = (d+e)/d$.  Since one of these roots is 
negative and the other two are between 0 and 1, the sum is less than 2.  Hence $(d+e)/d  < 2$.  Hence $e < d$.
We proved above that $F^\prime(1) \ge 0$.  Therefore 
\begin{eqnarray*}
0 &\le& F^\prime(1) \\
&=& 4d^2-3dp-4d(d+e)+p(d+e) \\
&=& -2dp-4de+pe \\
&=& p(e-2d) - 4de
\end{eqnarray*}
But since $e<d$,  and $e> 0$ and $d > 0$, the right side is negative, contradiction.  That completes the proof of the lemma.

\subsection{The case when $ABC$ has angles $2\alpha$, $\beta$, and $\beta + \alpha$ and $\sin(\alpha/2)$ is irrational}
Recall that $\sin(\beta + \alpha) = \sin \gamma = c$ in this case, as follows from $3\alpha + 2\beta = \pi$.
The $2\alpha$ angle is at vertex $A$, the $\beta$ angle at vertex $B$, and the $\beta + \alpha$ angle 
at vertex $C$.  We do not assume $\alpha < \beta$.  
The $\d$ matrix equation is 
$$ \matrix p d e g m f h \ell r \vector a b c = \vector X Y Z .$$
We suppose the tile is scaled so that each side is equal to the sine of the opposite angle.  Then 
\begin{eqnarray*}
a &=& \sin \alpha \\
b &=& \sin \beta \\
c &=& \sin (\alpha + \beta) \ = \ \sin \gamma  \mbox{\qquad since $\alpha + \beta = \pi - \gamma$}
\end{eqnarray*}
Let $\lambda = X/ \sin(2\alpha)$.  By the law of sines for triangle $ABC$ we then have 
\begin{eqnarray*}
X &=& \lambda \sin 2\alpha \\
Y &=& \lambda \sin \beta \\
Z &=& \lambda \sin(\beta + \alpha) 
\end{eqnarray*}
Define $t = (\sin 2 \alpha)/\sin \alpha = 2\cos \alpha$.  Then 
$\lambda at = \lambda (\sin \alpha) t = \lambda \sin 2 \alpha = X$.
Hence the $\d$ matrix equation can be written 
\begin{equation} \label{eq:MatrixEquation32}
 \matrix p d e g m f h \ell r \vector a b c =  \lambda \vector {at} b c 
\end{equation}

\noindent
Dividing the top row  
by $t$ we write it as an eigenvalue equation:
\begin{eqnarray*}
\matrix {\frac p t} {\frac d t } {\frac e t} g m f h \ell r \vector a b c =& \lambda \vector a b c 
\end{eqnarray*}
Our next goal is to prove that the eigenspace of the eigenvalue $\lambda$ is one-dimensional.
To do that, we use the method of computing a candidate eigenvector by the cofactor method, 
and showing that its three components are each nonzero.

 The matrix whose cofactors we need is 
\begin{eqnarray*}
\matrix  {\frac p t -\lambda} {\frac d t } {\frac e t} g {m - \lambda} f h \ell {r-\lambda}
\end{eqnarray*}

We expand in cofactors of the elements on the first row.  We find a candidate eigenvector $(u,v,w)$ where
\begin{eqnarray*}
 \vector u v w &=& \vector{ (m-\lambda)(r-\lambda) - f \ell} 
                            { hf - g(r-\lambda)}
                             {g\ell  -  h(m-\lambda)} 
 \end{eqnarray*} 
 \smallskip  
 
\begin{lemma} Suppose that $ABC$ is $N$-tiled by tile $T$ with $3\alpha + 2\beta = \pi$ and
$\sin(\alpha/2)$ irrational, and triangle $ABC$ has angles $2\alpha$, $\beta$, and $\beta + \alpha$.
Suppose also that not both $g$ and $h$ are zero.   Then the eigenspace of $\lambda$ is one-dimensional
and $(u,v,w)$ defined above is a multiple of $(a,b,c)$.                                                     
\end{lemma}

\noindent{\em Proof}.
 Our first aim is to prove $v \neq 0$.
  Since $f \neq 0$, 
 if $g=0$ then also $h=0$, which contradicts the hypothesis that not both $g$ and $h$ are zero, 
 and similarly, if $h=0$ then $g=0$. 
 Hence $g\neq 0$ and $h \neq 0$.  Then $hf > 0$.  But we have $grc \le \lambda c$ from the third row of the
 $\d$ matrix equation; hence $gr \le \lambda$.  Then
 \begin{eqnarray*}
 v &=& hf - g(r-\lambda) \\
 &=& hf + g(\lambda - r) \\
 &>& 0   \mbox{\qquad since not both $g=0$ and $h=0$}
 \end{eqnarray*}

Next we will prove $w \neq 0$. 
From the second row of the $\d$ matrix equation, we have $ga + mb \le \lambda b$.   Since $g > 0$,
we have $m < \lambda$, not just $m \le \lambda$.  Then
 \begin{eqnarray*}
  w &=& g\ell - h(m-\lambda)\\
  &=& g\ell + h(\lambda - m)\\
  &>& 0 \mbox{\qquad since not both $g=0$ and $h=0$}
 \end{eqnarray*}
 
 Finally we will prove $u \neq 0$.  We have from the second and third rows of the $\d$ matrix equation 
 \begin{eqnarray*}
 mb + fc  &\le& \lambda b \mbox{\qquad with equality only if $g=0$}\\
 \ell b + rc &\le & \lambda c \mbox{\qquad with equality only if $h=0$}
 \end{eqnarray*}
 Rearranging the terms of these equations we have 
 \begin{eqnarray*}
 fc &\le& (\lambda - m)b \\
 \ell b &\le& (\lambda - r)c
 \end{eqnarray*}
 Multiplying these two equations and dividing by $bc$ we have 
 \begin{eqnarray*}
 \ell f &\le& (\lambda - m)(\lambda - r)
 \end{eqnarray*}
 The left hand side is $u$.  We have proved $u \ge 0$ with equality only if $g=0$ and $h=0$.
But at present we have assumed that not both $g$ and $h$ are zero; so we have $u > 0$.

We have proved all three components of the candidate eigenvector are nonzero.  It follows that the 
eigenspace is one-dimensional; hence $(u,v,w)$ is a multiple of $(a,b,c)$.   That completes the proof of the lemma.

\begin{lemma} \label{lemma:singamma}
  Suppose $3\alpha + 2\beta = \pi$, and $s = 2\sin \alpha/2$.  Then we have (whether or not $s$ 
is rational)
\begin{eqnarray*}
\sin \gamma &=& \cos \frac \alpha 2 \\
\frac a c &=& s \\
\frac b c &=& 1-s^2 
\end{eqnarray*}
\end{lemma} 

\noindent{\em Proof}.   
Since $\gamma = \pi-(\alpha + \beta)$, we have 
\begin{eqnarray*}
\sin \gamma &=& \sin(\pi-(\alpha + \beta))  \\
 &=& \sin(\alpha + \beta) \nonumber \\
 &=& \cos(\pi/2- (\alpha + \beta))   \\
 &=& \cos \frac \alpha 2  \mbox{\qquad since $\pi/2 - \beta = 3\alpha/2$} 
 \end{eqnarray*}
Then $c = \sin \gamma = \cos \alpha/2$, and $a = \sin \alpha = 2 \sin(\alpha/2) \cos(\alpha/2)$.
Hence 
$$\frac a c = 2\sin \alpha/2.$$
Since $3 \alpha + 2\beta = \pi$, we have 
\begin{eqnarray*}
\sin \beta &=& \sin (\pi/2 - 3 \alpha/2) \\
&=& \cos(3 \alpha/2) \\
&=& 4 \cos^3 \frac \alpha 2 - 3 \cos \frac \alpha 2 
\end{eqnarray*}
Hence 
\begin{eqnarray*}
b/c &=& 4 \cos^2(\alpha/2) - 3 \\
&=& 4(1-\sin^2 \alpha/2) - 3 \\
&=& 1-4\sin^2 \alpha/2
\end{eqnarray*}
 Then we have 
\begin{eqnarray*}
\frac a c &=& s \\
\frac b c &=& 1-s^2
\end{eqnarray*}
That completes the proof of the lemma.
\smallskip

\begin{lemma} \label{lemma:PQlinearind}
Assume $ABC$ has angles $2\alpha$, $\beta$, and $\beta + \alpha$,
and $\sin(\alpha/2)$ is irrational. Then  $a$, $b$, and $c$ are 
linearly independent over $\Q$.   
\end{lemma}

\noindent{\em Proof}. 
Suppose, for proof by contradiction, that $a$, $b$, and $c$ are linearly dependent over $\Q$.  Then there are
(positive or negative) integers $n$, $j$, and $k$ such that $na + jb + kc = 0$.  Dividing by $c$
and using $a/c = s$ and $b/c = 1-s^2$ we find
$$0 = ns + j(1-s^2) + k.$$ 
If $j=0$ then $s = -k/n = 2\sin(\alpha/2)$ so $\sin(\alpha/2)$  is rational, contradicting the hypothesis that it is irrational.
Hence $j \neq 0$.
Solving for $1-s^2$ we have 
\begin{equation}  \label{eq:tsquaresub}
1-s^2 = -\frac {ns} j - \frac k {j} 
\end{equation}
There is an angle $\beta + \alpha$
at vertex $C$;  the adjacent sides are $X$ and $Y$ and 
 the area equation can be written as 
$$ Nab = XY  = (pa + db + ec)(ga + mb + fc)$$
 Dividing by $c^2$ we have 
$$ N \frac a c \frac b c = (p \frac a c + d \frac b c + e)(g \frac a c + m \frac b c + f) $$
Putting everything in terms of $s$ we have
\begin{eqnarray*}
 Ns(1-s^2) &=& ps + d(1-s^2) + e)(gs + m(1-s^2) + f  \\
 Ns \big(-\frac {ns} j - \frac k {j} \big)   &=& ps + d\big(-\frac {ns} j - \frac k {j} \big) + e)(gs + m \big(-\frac {ns} j - \frac k {j}  \big) f) 
 \end{eqnarray*}
 This is quadratic in $s$.
 We can then use (\ref{eq:tsquaresub})
to replace $s^2$ by a   term linear in $s$.  Hence $s$ satisfies a linear equation over $\Q$.  Hence $s$ is rational.
But $s = 2\sin(\alpha/2)$, so $\sin(\alpha/2)$ is rational.  It is not necessary to write the equation out explicitly.
 That completes the proof of the linear independence of $a$, $b$, and $c$.   
 
 \begin{lemma} \label{lemma:siscubic}  Suppose 
 $ABC$ has angles $2\alpha$, $\beta$, and $\beta + \alpha$ and $\sin(\alpha/2)$ is irrational. Then
 $s = 2\sin(\alpha/2)$ satisfies a cubic equation over $\Q$; hence its degree over $\Q$ is 3.
 \end{lemma}
 
 \noindent{\em Proof}.
 Taking the ratio of the first  row of (\ref{eq:MatrixEquation32}) to the third, we have 
\begin{eqnarray*}
\frac { pa + db + ec} {ha + \ell b + rc} &=& \frac {at} c
\end{eqnarray*}
Multiplying by the denominator and then dividing both sides by $c^2$ we have 
\begin{eqnarray*}
p \frac a c + d \frac b c + e &=& \frac {at} c \big(h \frac a c + \ell  \frac b c + r\big)
\end{eqnarray*}
We now express this equation in terms of $s$.  We have already derived  $a/c = s$ and $b/c = 1-s^2$;
we also have   $t = 2-s^2$, since $t$ is by definition $2\cos\alpha$, so we have 
$t = 2(1-2\sin^2(\alpha/2)) = 2-4\sin^2(\alpha/2) = 2-s^2$.  Putting these values into the 
previous equation, we have
\begin{eqnarray*}
ps + d(1-s^2) + e &=& s(2-s^2)(hs + \ell(1-s^2) + r)
\end{eqnarray*}
This is a polynomial equation of degree 5.  Bringing it to polynomial form we have 
\begin{eqnarray*}
0 &=& F(s) \ := \  \ell s^5 - hs^4 - (3\ell + r )s^3 + (d+2h)s^2 + (2\ell + 2r -p)s -(e+d)  
\end{eqnarray*}
If $F$ were identically zero, then we would have $\ell = 0$ from the leading coefficient, 
 and  then $r=0$ from the coefficient of $s^3$, contradicting Lemma \ref{lemma:zerolimits2}.
 Therefore $F$ is not identically zero.   
 
 Taking the ratio of the second row to the third, we have
 \begin{eqnarray*}
 \frac {ga + mb + fc} {ha + \ell b + rc} &=& \frac b c \\
 (ga + mb + fc) &=& (ha + \ell b + rc)\frac b c \\
 (g \frac a c + m \frac b c + f) &=& h \frac a c + \ell \frac b c + r) \frac b c \\ 
 (g s + m (1-s^2) + f) &=& ( hs + \ell (1-s^2) + r)(1-s^2)
 \end{eqnarray*}
 Subtracting the left side from the right and bringing the equation to polynomial form, we have 
 \begin{eqnarray*}
 0&=& H(s) \ := \ \ell s^4 - hs^3 + (m-2\ell - r)s^2 + (h-g)s + \ell + r-m-f
 \end{eqnarray*}
 If $H$ were identically zero, we would have $\ell =0$ from the leading coefficient, and $h=0$ from 
 the coefficient of $s^3$, then $r=m$ from the coefficient of $s^2$, then  $g=0$ from the coefficient of $s$,
 then $f=0$ from the constant coefficient; but by Lemma \ref{lemma:zerolimits2}, we do not have $g=f=0$.
  Hence $H$ is not identically zero.  We compute $K = -(F \mod H)$ (using the {\em PolynomialRemainder}
 command in Mathematica, for example), and we find 
 \begin{eqnarray}
 0 &=& K(s) \ = \ 
  (\ell + m) s^3 - (d+g +h)s^2 - (f+m + \ell + r -p)s  + (d+e)  \label{eq:April27-1}
 \end{eqnarray}
 If $K$ were identically zero, then we would have $\ell = m = 0$ from the coefficient of $s^3$, and 
 $d = g = h = 0$ from the coefficient of $s^2$, and then $e=0$ from the constant term, but $e=0$
 contradicts Lemma \ref{lemma:zerolimits1}.   Hence $K$ is not identically zero.
 Assume, for proof by contradiction, that the leading term of $K$ is zero. Then $m=\ell =0$ 
 and $H(s) = -hs^3 -rs^2 + (h-g)s -(m+f) = 0$ 
 That completes the proof of the lemma.
 
 \begin{lemma} Under the hypotheses in the section title,   $1$, $\sin(\alpha/2)$, and $\cos(\alpha)$ are linearly independent.
 \end{lemma}
 
 \noindent{\em Proof.}
  We have $\cos \alpha = 1-2\sin^2(\alpha/2) =
 1-(1/2)s^2$, and $\sin (\alpha/2) = s/2$, so the claim is equivalent to the claim that $1$, $s$, and $s^2$ 
 are linearly independent over $\Q$;  that is, that $s$ does not have degree 2 over $\Q$.    By Lemma 
 \ref{lemma:siscubic}, $s$ satisfies a cubic equation; hence its degree over $\Q$ is either 1 or 3, not 2. 
\medskip

\medskip
\begin{lemma} \label{lemma:notghzero}
Let $3\alpha + 2\beta = \pi$ and assume triangle $ABC$ is $N$-tiled by
a tile with angles $\alpha$ and $\beta$.  Suppose $ABC$  has angles $2\alpha$, $\beta$, and $\beta + \alpha$,
and suppose $s = 2\sin(\alpha/2)$ is not rational. 
Then the $\d$ matrix entries $g$ and $h$ are not both zero.
\end{lemma}

\noindent{\em Proof}.  Suppose, for proof by contradiction, that $g=h=0$.  The $\d$ matrix equation is, by (\ref{eq:MatrixEquation32}),
\begin{eqnarray*}
\matrix  p d e 0 m f 0 \ell r \vector a b c &=& \lambda \vector {at} b c
\end{eqnarray*}
Writing this in eigenvalue form we have
\begin{eqnarray*}
\matrix  {\frac p t - \lambda} {\frac d t} {\frac e t} 0 {m - \lambda} f 0 \ell {r - \lambda} \vector a b c &=& 0
\end{eqnarray*}
Hence the determinant of the matrix on the left is zero.  Expanding it by cofactors we have
\begin{eqnarray*}
\big(\frac p t - \lambda\big)\big( (m-\lambda)(r-\lambda) - f \ell)\big) &=& 0
\end{eqnarray*}
There are two cases: either $p = \lambda t$ or not.  In case $p = \lambda t$ then the last row of the 
$\d$ matrix equation is
\begin{eqnarray*}
\ell b + r c &=& \lambda c 
\end{eqnarray*}
Dividing by $c$ and using $b/c = 1-s^2$ and $\lambda = p/t$ we have
\begin{eqnarray*}
\ell (1-s^2) + r &=& \frac p t \\
&=& \frac p {2-s^2} \\
\ell (1-s^2)(2-s^2) &=& p 
\end{eqnarray*}
Thus $s^2$ satisfies a quadratic equation over $\Q$.  Hence the degree of $\Q(s)$ over $\Q$ is either 2 or 4,
since $s$ is not rational by hypothesis.  But this contradicts  Lemma \ref{lemma:siscubic},
  which says $s$ satisfies a cubic equation
and hence has degree 3 or 1.  That completes the case $p=\lambda t$.  

Now assume $p \neq \lambda t$. Then the second factor in the characteristic equation is zero:
\begin{eqnarray*}
(m-\lambda)(r-\lambda) - f \ell) &=& 0
\end{eqnarray*}
This is a quadratic equation for $\lambda$;  hence $\Q(\lambda)$ has degree 2 or 1 over $\Q$.  
But from the third row of the $\d$ matrix equation we have 
$$ \lambda = \ell (1-s^2) + r$$
which implies that the degree of $\Q(s)$ over $\Q(\lambda)$ is 1 or 2.   Hence the degree of $\Q(s)$ over
$\Q$ is 1, 2, or 4.   Degree 1 contradicts the hypothesis that $s$ is irrational;  degrees 2 or 4 contradict
(\ref{eq:April27-1}).   That completes the proof of the lemma.
\medskip 
 
\begin{lemma} Let $3\alpha + 2\beta = \pi$ and assume triangle $ABC$ has angles $2\alpha$ and $\beta$ and is $N$-tiled by
a tile with angles $\alpha$ and $\beta$.  Suppose   $\sin(\alpha/2)$ is irrational.
Then $\ell \neq 0$, $p \neq 0$, $h \neq 0$, $g\ell - hm \neq 0$, and 
 we have the following equations between the elements of the $\d$ matrix, where $\Delta$ is the determinant of 
 the $\d$ matrix:
\begin{eqnarray*}
\frac p {g+h}  &=& \frac {mr-\ell f} {g\ell -hm} \\
 \frac p {g+h} &=& \frac{-p(m+r) +eh +dg -N}{g\ell-hm + hf-gr}  \\
\frac p {g+h} &=& \frac {\Delta + N(m+r)} {-Nh} 
\end{eqnarray*}
and in the second equation, if the denominator on the right is zero, so is the numerator.
\end{lemma}

\noindent{\em Proof}.  
Since $t = N/\lambda^2$,  $t$ belongs to $\Q(\lambda)$.
The third 
row of the $\d$ matrix equation, $h a + \ell b + r c = \lambda c$,  shows that $\lambda$ belongs
to $\Q(a,b,c) = \Q(t)$.  Since $t$ belongs to $\Q(\lambda)$ and $\lambda$ belongs to $\Q(t)$, 
we have $\Q(\lambda) = \Q(t) = \Q(a,b,c)$.   This field has degree at least 3 over $\Q$,
 by Lemma \ref{lemma:PQlinearind}.  We will soon see that the degree is exactly 3. 

We are in a position to find three different cubic equations for $\lambda$.  
The first equation 
come from the characteristic equation for $\lambda$:
$$ \matrix {\frac p t - \lambda} {\frac d t} {\frac e t} g {m-\lambda} f h \ell {r-\lambda} = 0$$
Multiplying the top row by $t$, the determinant is still zero:
$$ \matrix {p - \lambda t} d e g {m-\lambda} f h \ell {r-\lambda} = 0$$
Since $\lambda^2 = N/t$, we have $\lambda t = N/\lambda$. Multiplying by $\lambda$ we have
$$
 \lambda \matrix {p - \frac N \lambda} d e g {m-\lambda} f h \ell {r-\lambda} = 0
$$
Expanding and collecting like powers of $\lambda$ we find the following equation, in which 
$\Delta$ is the determinant of the $\d$ matrix:
\begin{equation} \label{eq:cubic1}
p \lambda^3 + ( -p(m+r) +eh +dg -N)\lambda^2 + (\Delta + N(m+r)) \lambda - N(mr-\ell f) = 0
 \end{equation}
This is a cubic polynomial equation for $\lambda$.   Hence $\Q(\lambda)$ has degree 3 over $\Q$, and this 
is the minimal polynomial of $\lambda$.  Since we have already proved that the eigenspace of $\lambda$ is 
one-dimensional, it is not the case the all the coefficients are zero.  Since the degree of $\Q(\lambda)$
is 3, the highest and lowest coefficients are not zero:  $p \neq 0$ and $mr \neq 0$.

Next we compute $b/c$ in terms of $\lambda$:
\begin{eqnarray*}
 \frac b c &=& ( z^2 - 1 + z^2)\\
&=&  t-1 \mbox{\qquad since $t = z^2 + z^{-2}$} 
\end{eqnarray*}
Since $\lambda^2 = N/t$, we have $t  = N/\lambda^2$, which gives us 
\begin{equation} \label{eq:lambdaboverc}
 \frac b c = \frac  N {\lambda^2} - 1
\end{equation}
Since the eigenvector $(u,v,w)$ is a multiple of $(a, b,c)$, we have 
\begin{eqnarray*}
\frac b c &=& \frac v w  \\
\frac N {\lambda^2} - 1 &=& \frac {hf - g(r-\lambda)}{g \ell - h(m-\lambda)}
\end{eqnarray*}
Cross multiplying we have 
\begin{eqnarray*}
(N - \lambda^2)(g\ell - hm + h \lambda) &=& (hf-gr)\lambda^2 + g\lambda^3
\end{eqnarray*}
Subtracting the left side from the right and expanding and collecting, we have 
\begin{eqnarray}
0 &=& (g+h)\lambda^3 + (g\ell -hm + hf-gr) \lambda^2  -Nh \lambda  - N(g\ell  - hm) \label{eq:cubic2}
\end{eqnarray}
By Lemma \ref{lemma:notghzero},   not both $g$ and $h$ are zero.    Hence (since both are nonnegative) $g+h > 0$,
and this is a non-trivial equation. Since $\Q(\lambda) = \Q(a,b,c)$ (from the third row of the $\d$ matrix
equation), and this field has degree 3,  the constant term  $N(g\ell-hm)\neq 0$ (or else the equation could be 
divided by $\lambda$ to yield a quadratic equation for $\lambda$).  Hence $g\ell - hm \neq 0$.
Similarly, the constant term of (\ref{eq:cubic1}), namely $N(mr-\ell f)$, is also not zero, so $mr -\ell f\neq 0$.
 
Now (\ref{eq:cubic2}) and (\ref{eq:cubic1}) are two cubic equations for $\lambda$.  
Therefore (\ref{eq:cubic2}) is a multiple of (\ref{eq:cubic1}).  The ratio of the coefficients of $\lambda^3$ is 
$p/(g+h)$.  Therefore the other nonzero coefficients are also in that ratio.  
From the constant coefficients (neither of which is zero) we have
\begin{eqnarray}
\frac p {g+h}  &=& \frac {mr-\ell f} {g\ell -hm} \label{eq:64}
\end{eqnarray}
    It follows that
\begin{equation} \label{eq:ellnonzero}
\ell \neq 0
\end{equation}
since otherwise the numerator and denominator on the right of (\ref{eq:64}) have opposite signs, but the 
left side is positive.
From the ratio of the quadratic terms we have (unless the numerator and denominator are both zero) 
\begin{eqnarray}
\frac p {g+h} &=& \frac{-p(m+r) +eh +dg -N}{g\ell-hm + hf-gr}    \label{eq:66} 
\end{eqnarray}
From the linear terms we have (unless both numerator and denominator are zero) 
\begin{eqnarray}
\frac p {g+h} &=& \frac {\Delta + N(m+r)} {-Nh} \label{eq:67}
\end{eqnarray}
These are the three equations mentioned in the lemma.  It remains to prove $h\neq 0$.
Assume, for proof by contradiction, that $h = 0$.  Then $ g \neq 0$, since the constant term $g+h$ 
of (\ref{eq:cubic2}) is not zero.  Then from (\ref{eq:64}) we have $p\ell = mr$,
or $m = p\ell/r$.  When $h=0$ we have 
\begin{eqnarray*}
 \Delta &=& \determinant p d e g m f 0 \ell r \\
 &=& pmr  + eg\ell - gdr - pf\ell
\end{eqnarray*}
Since both numerator and denominator on the right of (\ref{eq:67}) are zero we have 
\begin{eqnarray*}
0 &=& \Delta + N(m+r) \\
&=& pmr  + eg\ell - gdr - pf\ell + N(m+r) 
\end{eqnarray*}
Substituting $m = p\ell/r$ and multiplying by $r$ we have 
\begin{eqnarray*}
0 &=& p^2 \ell r + eg\ell r - gdr^2 -pf\ell r+ N(p\ell + r^2) 
\end{eqnarray*}
Bringing the last term to the left side and changing signs we have
\begin{eqnarray*}
N(p\ell+r^2) &=& pf\ell r + gd r^2 - eg\ell r - p^2 \ell r \\
&\le& pf\ell r + gdr^2 
\end{eqnarray*}
Note that $gd \le \lambda^2$, since $ga \le \lambda b$ and $db \le \lambda a$, from the $\d$ matrix.
We also have $r \le \lambda$ and $fc \le \lambda b < \lambda c$, so $fr < \lambda^2$.  Then 
\begin{eqnarray*}
N(p\ell + r^2) &<& p\ell \lambda^2 + r^2 \lambda^2 \\
  &=& (p\ell + r^2) \lambda^2
\end{eqnarray*}
Since $p\ell + r^2 \neq 0$ we have 
\begin{eqnarray*}
N &<& \lambda^2
\end{eqnarray*}
which is a contradiction, since $\lambda^2 = N/t$ and $t = 2 \cos \alpha > 1$ since $\alpha > \pi/5$.  
That contradiction completes the proof that $h \neq 0$.  Hence equation (\ref{eq:67}) is valid, i.e.
its denominator is not zero.   That completes the proof of the lemma.%
\footnote{After finding these equations in March 2010, it was still some time before I could show that
they had no solutions; in the meantime, I wrote a short C program that checked there are no solutions 
for $N \le 400$.  That put an end to my attempts to find a tiling of this kind by placing paper triangles
on a card table.}

\begin{lemma} \label{lemma:3alpha2betaCase1}
 Suppose $3\alpha + 2\beta = \pi$, and triangle $ABC$ has angles $2\alpha$ and $\beta$ and is
not isosceles.  Suppose  $\sin(\alpha/2)$ is irrational.  Let $N$ be arbitrary.  Then there is no $N$-tiling of 
triangle   $ABC$ by a tile with angles $\alpha$ and  $\beta$.
\end{lemma}

\noindent{\em Proof}.   Suppose, for proof by contradiction, that there is a triangle $ABC$ and an $N$-tiling
as in the lemma.   Then the equations of the previous lemma hold.  From (\ref{eq:67}) we have 
\begin{eqnarray*}
N(m+r) &\le& -\Delta - N h \frac p {g+h} \\
N(m+r) + Nh \frac p {g+h} &\le& - \Delta \\
N\big(m+r+  \frac {hp}{g+h}\big) &\le -\Delta \\
&\le& dgr + pf \ell + hme - pmr - dfh - eg\ell \\
&\le& mhe + rdg +  p(\ell f - mr) - dfh - eg\ell
\end{eqnarray*}
From the third row of the $\d$ matrix we have $ha < Z = \lambda c$. Equality cannot hold since $\ell \neq 0$ (and we have $ha + \ell b < Z$).
From the first row we have $ec \le \lambda at$.  Multiplying these two inequalities (which is legal since all
these quantities are positive) and dividing by $ac$ we have
$he \le \lambda^2 t < N$.  Similarly, we have $dg \le \lambda^2 t  = N$.   Putting these results into 
the inequality above, we have 
\begin{eqnarray*}
N(m+r) + N \frac {hp} {g+h} &<& N(m+r) + p( \ell f - mr) - dfh - eg\ell
\end{eqnarray*}
Note that the inequality has become strict.  Subtracting $N(m+r)$ from both sides we have
\begin{eqnarray*}
N \frac {hp} {g+h} &<& p(\ell f-mr) - dfh - eg\ell
\end{eqnarray*}
We have $p \neq 0$ and $mr -\ell f < 0$,
since otherwise the right side would be negative, which would be a 
contradiction. It follows from (\ref{eq:64}) that 
\begin{equation} \label{eq:Rneg}
g\ell -hm < 0
\end{equation}
From the second row of the $\d$ matrix equation we have  $fc \le \lambda b$.
From the third row we have $\ell b \le \lambda c$. Hence $\ell f \le \lambda^2$.  Hence
\begin{eqnarray*}
N \frac {hp}{g+h} &<& p\lambda^2 - pmr - dfh - eg\ell
\end{eqnarray*}
Now  add $Np g(g+h)$ to both sides of this inequality:
\begin{eqnarray*}
N \big( \frac {hp}{g+h} + \frac {gp}{g+h} \big) &<& p\lambda^2 + N \frac {pg}{g+h} - pmr -dfh - eg\ell \\
Np &<&  p \lambda^2 + N\frac {pg}{g+h} - pmr -dfh - eg\ell 
\end{eqnarray*}
Subtracting $p\lambda^2$ from both sides we have
\begin{eqnarray*}
p(N-\lambda^2) &<&  N\frac {pg}{g+h}  - pmr -dfh - eg\ell  
\end{eqnarray*}
Dividing by $p$ (which we proved above is not zero), we have
\begin{eqnarray}
N-\lambda^2 &<&  N\frac {g}{g+h}  - mr - \frac {dfh} p - \frac{eg\ell} p \label{eq:69}
\end{eqnarray}
The left side is positive, and only 
the first term on the right can be positive, and it is small when $h$ is large.  This observation motivates
us to investigate the ratio $h/g$. 
 From (\ref{eq:Rneg}) we have $g\ell < mh$. Then $\ell/m < h/g$.   
 From (\ref{eq:64}) and (\ref{eq:Rneg}) we see that the numerator of (\ref{eq:64}) is negative, i.e. $mr-\ell f < 0$.
 Then $r < (\ell/m) f < (h/g)f$.  Hence $r/f < h/g$.  Now we have $\ell /m $ and $r/f$ each less than $h/g$.
 From the bottom row of the $\d$ matrix equation we then have 
 \begin{eqnarray*}
 Z &=& ha + \ell b + rc \\
 &=& ga \frac h g + m \frac \ell m  b + f  \frac r f c \\
 &<& \big(\frac h g \big)(ga + m b + f c ) \\
 &=& \big( \frac h g \big) Y \\
 \frac Z Y &<& \frac h g \\
 \end{eqnarray*}
 By the law of sines we have $Z/Y = \sin(\gamma) / \sin(\beta) = c/b$.  Hence
 \begin{eqnarray*}
 \frac c b &<& \frac h g
 \end{eqnarray*}
 Applying (\ref{eq:lambdaboverc}) to write $c/b = \lambda^2 /(N-\lambda^2)$, we have
 \begin{eqnarray}
 \frac h g &>& \frac {\lambda^2} {N-\lambda^2} \label{eq:hgbound}
 \end{eqnarray}
 From (\ref{eq:69}) we have 
 \begin{eqnarray*}
 N-\lambda^2 &<& N \frac g {g + h} \\
 &=& N \frac 1 { 1 + h/g} 
 \end{eqnarray*}
 Now we replace $h/g$ in the denominator on the right by the smaller quantity from (\ref{eq:hgbound}). 
 Then the right hand side increases.
 \begin{eqnarray*}
  N-\lambda^2 &<& N \frac 1 { 1 + \frac {\lambda^2} {N-\lambda^2}} \\
 &=& N \frac {N-\lambda^2}{N- \lambda^2  +   \lambda^2} \\
 &=& N \frac {N-\lambda^2} N \\
 &=& N-\lambda^2
 \end{eqnarray*}
 We have derived $N-\lambda^2 < N-\lambda^2$, which is a contradiction.
 That completes the proof of the lemma.

 \subsection{The case when $ABC$ has angles $2\alpha$, $\beta$, and $\beta + \alpha$, and $\alpha$ is a rational multiple of $\pi$}
 
  Note that our work above has assumed that $\alpha$ is not a rational multiple of $\pi$ only 
to limit the possible shapes of $ABC$.  The following lemma shows that there are only two possible rational multiples of $\pi$
that we might need to consider. 
\begin{lemma} \label{lemma:46} Let $\alpha$ be a rational multiple of $\pi$, and suppose $3\alpha + 2\beta = \pi$.  Let triangle $ABC$ have one angle $\beta$, one angle $2\alpha$,
and one angle $\beta + \alpha$.  Suppose that there is an $N$-tiling of $ABC$ by the tile with angles $\alpha$ and $\beta$.
Then $\alpha = 2\pi/7$ and $\beta  = \pi/14$,  or $\alpha = 2\pi/9$ and $\beta = \pi/6$.   In either case, $\sin(\alpha/2)$ 
is irrational.
\end{lemma}

\noindent{\em Proof}.  Suppose $\alpha = 2m\pi/n$.  Then the degree $d$ of $\Q(e^{i\alpha})$ over $\Q$  is finite and 
by Lemma \ref{lemma:euler}, we have $d = \varphi(n)$.  If both $\sin(\alpha)$ and $\cos(\alpha)$ are rational
we have $d = 2$, so $\varphi(n) = 2$, so $n=3$ or $4$.  Remember we are not assuming $\alpha < \beta$;  so the 
smaller one of $\alpha$ and $\beta$ is less than $\pi/5$, and $\alpha < \pi/3$ and $\beta < \pi/2$.   Hence $\alpha$
cannot be $2\pi/3$ or $2\pi/4$.  Hence not both $\sin \alpha$ and $\cos \alpha$ are rational. 

Assume first that $\sin(\alpha/2)$ is not rational.
By (\ref{eq:April27-1}),  the degree of $\Q(s)$ over $\Q$ is at most 3.  By Lemma \ref{lemma:PQlinearind}, 
under the assumption that $\sin(\alpha/2)$ is irrational, that degree is at least 3, and hence it is exactly 3.   Since $s = 2\sin(\alpha/2)$,
$\Q(s)$ is the real subfield of $\Q(e^{i\alpha/2})$, so the degree of that field is 6.   Since $\alpha = 2m\pi/n$, we have
$\alpha/2 = 2m\pi/(2n)$.
Hence by Lemma \ref{lemma:euler}, the degree of $\Q(e^{i\alpha/2})$ over $\Q$ is $\varphi(2n)$.   The possibilities are exactly these: 
$2n = 7, 9, 14$, or $18$.   Since $n$ is an integer, we cannot have $2n=7$ or 9; therefore $2n=14$ or $18$, i.e. $n=7$ or $n=9$.
Hence $\alpha = 2m\pi/7$ or $2m\pi/9$.  Since $\alpha < \pi/3$, we must have $m=1$ and $\alpha = 2\pi/7 $ or $\alpha = 2\pi/9$.
In case $\alpha = 2\pi/7$, we have $\beta = (\pi-3\alpha)/2 = \pi/14$, and in case $\alpha = 2\pi/9$ we have $\beta = \pi/6$.
That completes the proof in case $\sin(\alpha/2)$ is irrational.  

Suppose then that $\sin(\alpha/2)$ is rational.  If $\cos(\alpha/2)$ is also rational, then the degree of $\Q(e^{i\alpha/2})$ 
over $\Q$ is 2, so by Lemma \ref{lemma:euler}, we have $2 = \varphi(2n)$, so $2n = 3$ or $4$;  so $n=2$,  so $\alpha = 2\pi/2 = \pi$,
which is impossible.  Hence $\cos(\alpha/2)$ has degree 2 over $\Q$, so $\Q(e^{i\alpha/2})$ has degree 4 over $\Q$.  Hence
we have $4 = \varphi(2n)$, so by Lemma \ref{lemma:euler}, $2n = 5, 8, 10$, or $12$.  Hence $n = 4, 5$, or $6$.  Hence 
$\alpha = \pi/2$, $2\pi/5$, or $\pi/3$.   Since $\alpha < \pi/3$, none of these is possible.  
That completes the proof of the lemma.

\section{The general case when $T$ is not similar to $ABC$}

We have exhibited a few tilings in which the tile $T$ is not similar to $ABC$.  In this section we shall 
prove that if there are such tilings,  then the number of tiles and the shape of the tile and triangle 
correspond to the exhibited families.  In previous sections we have dealt with the cases when $T$
is isosceles (and isosceles includes equilateral), and when $T$ is similar to $ABC$.  Therefore we now assume that $T$ is 
not isosceles and not similar to $ABC$, and $ABC$ is $N$-tiled by $T$.   These assumptions will be in force
for this entire section.

Let $P$, $Q$, and $R$ be the total number of $\alpha$, $\beta$, and $\gamma$ angles (respectively) occurring at the vertices 
of triangle $ABC$.  These numbers control the ``vertex splitting'', i.e., the way the vertices of $ABC$ are divided into 
the angles of $T$.  We have 
$$P \alpha + Q \beta + R \gamma = \pi.$$
Together with $\alpha +\beta + \gamma = \pi$,  that makes two equations for $\alpha$, $\beta$, and $\gamma$.
We seek another equation that would enable us to solve for the three angles.   At each vertex $V$ of the tiling,
which is not a vertex of $ABC$, 
there meet several copies of the tile with vertices at $V$.  If $V$ is a non-strict or a boundary vertex, then the angle sum 
of the vertex angles is $\pi$, and if it is a strict interior vertex, the angle sum is $2\pi$.  We set $k=1$ for a non-strict
or boundary vertex, and $k=2$ for a strict interior vertex.  Then we have 
$$ p \alpha + q \beta + r \gamma = k\pi$$
where $p$, $q$, and $r$ are the numbers of copies of $T$ with their $\alpha$, $\beta$, or $\gamma$ angles at $V$.

The numbers $P$, $Q$, and $R$ depend only on $ABC$ and the tiling, but $p$, $q$, and $r$ depend on the particular vertex $V$.
Our plan is to analyze the possibilities for $P$, $Q$, and $R$,  and given those possibilities, to analyze the possibilities
for $p$, $q$, and $r$,  in such a way as to eliminate all but a manageable (small finite) number of special cases for the 
tile $T$.

\begin{lemma} \label{lemma:R0}
$P + Q + R \ge 5$.
\end{lemma}

\noindent{\em Proof}.  Suppose $P+Q+R < 5$.  Then there must be two vertices of  $ABC$ that are not split, since 
each vertex contributes at least one to the sum $P+Q+R$, and split vertices contribute at least 2.  
Each vertex that does not split is equal to one of the angles of $T$. 
Since $T$ is not isosceles, that means that $T$ has two angles equal to angles of $ABC$.  Hence $T$ is similar to $ABC$, contrary to hypothesis.  That completes the proof of the lemma.

\begin{lemma} \label{lemma:R2}
 $R < 2$, i.e. $R=0$ or $R=1$.
\end{lemma}

\noindent{\em Proof}. 
If $\gamma = \pi/3$,
then $T$ is equilateral, contrary to hypothesis. Since $\gamma \ge \pi/3$ for any triangle, we have $\gamma > \pi/3$.
Hence $R \le 3$.  But if $R=3$, then each angle of $ABC$ is at least $\gamma$.  Since the sum of the angles of $ABC$ is 
$\pi$, this is possible only if $\gamma = \pi/3$, which we have shown is false.  Hence $R < 3$;  but $R$ is an integer,
so $R \le 2$. Hence, it suffices to show $R \neq 2$.
Assume, for proof by contradiction, that $R = 2$.
Then $P\alpha + Q \beta + 2\gamma = \pi$.  Subtracting $\alpha + \beta + \gamma = \pi$,
we find $(P-1)\alpha + (Q-1) \beta + \gamma = 0$.  We must have $Q=0$ or $P=0$, since if not the left side is positive.
If $P=0$ then $\alpha = (Q-1) \beta  + \gamma \ge \beta $, contradiction.  Hence $Q=0$.  Then 
our equation simplifies to $\beta = (P-1)\alpha + \gamma \ge \gamma$, contradicting $\beta < \gamma$.    That completes the proof of the lemma.
 
\begin{lemma} \label{lemma:R1}  If $R=1$ then $Q=0$ and $\beta = (P-1)\alpha$.
\end{lemma}

\noindent{\em Proof}. Assume $R=1$.  Then $P\alpha + Q \beta + \gamma = \pi$.  Subtracting $\alpha + \beta + \gamma = \pi$,
we find $(P-1)\alpha + (Q-1) \beta = 0$.   Since there is vertex splitting and $T$ is not similar to $ABC$,
we have $P+Q \ge 4$ by Lemma \ref{lemma:R0}.   Hence we must have $Q=0$ or $P=0$.
If $P=0$ then $\alpha = (Q-1) \beta \ge \beta$, contradiction.  Hence $Q=0$ and $P \ge 4$,
and $\beta = (P-1)\alpha$.   That completes the proof of the lemma.

\begin{theorem} \label{theorem:NotSimilar2}
Suppose that triangle $ABC$ is $N$-tiled by triangle $T$, and $T$ is not similar to $ABC$. 
Suppose also that if $\gamma = 2\pi/3$ then $T$ is isosceles.   Then
one of the following holds:
\smallskip

(i)  $N$ has the form $3m^2$  or $6m^2$ for some integer $m$ and $ABC$ is equilateral; and the tile is (in case $N = 3m^2$) the tile used 
in the equilateral 3-tiling, that is, $\gamma = 2\pi/3$,  $\alpha = \beta = \pi/6$, or (in case $N = 6m^2$) it is half of that tile.
\smallskip

(ii) $N$ is twice a square, or six times a square, or twice a sum of two squares;  $\gamma = \pi/2$; and  $ABC$ is isosceles, having angles $\alpha$, $\alpha$, and $2\beta$,  or  angles
$\beta$, $\beta$, and $2\alpha$, so that $T$ is a right triangle similar to half of $ABC$.
\smallskip

(iii) $3\alpha + 2\beta = \pi$,  triangle $ABC$ has angles $2\alpha$, $\beta$, and $\beta + \alpha$,  
$s = \sin(\alpha/2)$ is rational, and $N \ge 28$.   Here we do not assume $\alpha < \beta$. 

\end{theorem}
\medskip

\noindent{\em Remark}. Case (iii) is taken up in \cite{beeson-triquadratics}, where a new family of tilings 
is presented, and a necessary and sufficient condition on $N$ is given for such tilings to exist.
\medskip

\begin{corollary}
If $N$ is not divisible by 2 or by 3 (or even if it is, and the 
quotient is not a square),
there are no $N$-tilings in which $T$ is not similar to $ABC$. 
\end{corollary}
\smallskip

\noindent{\em Proof}.  
  Consider the case $\gamma = \pi/2$
and $\alpha = \pi/4$, i.e. $T$ is an isosceles right triangle.  If two or more angles $\alpha$ occur
at a vertex of $ABC$, then that vertex angle is at least $\pi/2$.  Hence at most one vertex angle of $ABC$ can be split.
Hence two angles of $ABC$ are equal to $\pi/2$ or $\pi/4$.  We cannot have two angles equal to $\pi/2$;  if two angles are 
equal to $\pi/4$ then $ABC$ is similar to $T$, contradiction.  Hence one angle of $ABC$ is $\pi/2$ and one is $\pi/4$; but 
then triangle $ABC$ is similar to $T$, contradiction.  Hence the case $\gamma = \pi/2$ and $\alpha = \pi/4$ is ruled out.

Now consider the case $\gamma = \pi/2$ and 
$\alpha < \beta$; that is, the tile is a non-isosceles right triangle.  
By Theorem 5 of \cite{beeson1},  triangle $ABC$ must be isosceles and 
$T$ is similar to half of $ABC$.  By Theorem 5 of \cite{beeson1}, $N$ is twice a square, six times a square, or twice a sum of two squares, so conclusion (ii) of the theorem holds.

Therefore we may assume, for the rest of the proof, that $\gamma \neq \pi/2$.
In Lemma \ref{lemma:twopioverfive} we have ruled out the case $\gamma = 2\pi/5$.
Hence we can assume that $\gamma$ is not equal to
$2\pi/3$ , $\pi/2$, or $2\pi/5$.  Moreover, $\gamma$ is not equal to $\pi/3$, since in that case $T$ is equilateral, and 
in that case no vertex splitting can occur, so $T$ must be similar to $ABC$.  

We now continue with the proof of the theorem.
The inequality $R < P+Q$ means that,  at the vertices of $ABC$ (taken together),  there are more $\alpha$ and $\beta$
angles than $\gamma$ angles.   Then there exists a vertex $V$ at which there are more $\gamma$ angles that $\alpha$ and $\beta$
angles.   Let $n$, $m$, and $\ell$ be the number of $\alpha$ angles, $\beta$ angles, and $\gamma$ 
angles at $V$, respectively.  Then $\ell > n + m$.    
  Let the angle sum at $V$ be $k\pi$ (so $k=1$ or $k=2$).   It is not the case that both $n$ and $m$ are zero, 
since then $k\pi = \ell \gamma$; since $\gamma \ge \pi/3$ the only possibilities are $\gamma = 2\pi/3$, $\gamma = \pi/2$,
$\gamma = 2\pi/5$, and $\gamma = \pi/3$, all of which have been ruled out already.

If $\ell = 1$ then (since $\ell > n + m$) we have $n=m=0$, contradiction.  If $\ell = 2$ then 
$m+n=0$ or $m+n=1$; if $m+n=0$ then $m=n=0$, which we have ruled out above. If $\ell = 2$ and
and $m+n=1$ then $2\gamma + \alpha = k\pi$ or $2\gamma+\beta = k\pi$. If $k=1$
we can subtract $\alpha + \beta + \gamma = \pi$ to obtain $\gamma = \beta$ or $\gamma = \alpha$;
but since $ABC$ is not isosceles this is a contradiction. Hence $k=2$ and we have 
$2\gamma + \alpha = 2\pi$ or $2\gamma + \beta = 2\pi$.   Writing $2\pi = 2\alpha + 2\beta + 2\gamma$ and 
subtracting the left-hand side we have
$ 0 = \alpha + 2\beta$ or $0 = \beta + 2\alpha$, so $\alpha = \beta = 0$, which is a contradiction. 
That rules out the case $\ell = 2$.

Since $\gamma > \pi/3$,  we have $\ell \le 5$.   That leaves the possibilities 3, 4, and 
5 for $\ell$.  We will rule out each of these in turn;  each one requires a detailed argument.  

We first assume $\ell = 3$. Since $\gamma > \pi/3$,  we have $k=2$, i.e. the angle sum at $V$ is $2\pi$;
and we cannot have $n=m=0$, since then $\gamma$ would be $\pi/3$. 
Since $m+n < \ell$, we have $m+n \le 2$.
If $n=1$ and $m=1$ then 
\begin{eqnarray*}
2\pi &=& 3\gamma + \alpha + \beta \\
&=& 2\gamma + (\gamma + \alpha + \beta) \\
&=& 2\gamma + \pi 
\end{eqnarray*}
Hence $\gamma = \pi/2$, contradiction.  Hence $n=1$ and $m=1$ is impossible.

If $n=2$ and $m=0$, then
\begin{eqnarray*}
2\pi &=& 3 \gamma + 2 \alpha \\
&=& \gamma + 2\gamma + 2\alpha \\
2(\gamma + \alpha + \beta) &=& \gamma + 2\gamma + 2 \alpha \\
2\beta &=& \gamma  \\
\end{eqnarray*}
Assume $R=1$.  Then  by Lemma \ref{lemma:R1},
$\beta = (P-1)\alpha$.   Then since $\gamma = 2\beta$, we have 
\begin{eqnarray*}
\pi &=& \alpha + \beta + \gamma  \\
&=& \alpha + 3 \beta \\
&=& \alpha + 3(P-1)\alpha \\
&=& \alpha (3P-2)
\end{eqnarray*}
We can then solve for all three angles:
\begin{eqnarray*}
\alpha &=& \frac \pi {3P-2} \\
\beta &=&  \frac {P-1}{3P-2} \pi \\
\gamma &=& \frac {2(P-1)}{3P-2} \pi
\end{eqnarray*}
In that case the inequality 
$\gamma \le \pi/2$ is equivalent to $2(P-1)/(3P-2) \le 1/2$, which is equivalent to $P \le 2$.   Since 
$P \ge 4$ we have $\gamma > \pi/2$ in this case.  

If $R=0$ then we have 
\begin{eqnarray*} 
\pi &=& P\alpha + Q \beta  \\
2\pi &=& 3\gamma + 2\alpha \mbox{\qquad since $\ell = 3$, $n=2$, and $m=0$} \\
\gamma &=& 2\beta \mbox{\qquad by subtracting $2\gamma  + 2\alpha$} \\
\pi &=& \alpha + 3\beta \mbox{\qquad since $3\gamma + 2\alpha = 2\pi$ and $\gamma = 2\beta$}
\end{eqnarray*}
Since $P+Q+R \ge 5$, and $R=0$, we do not have $P=1$ and $Q=3$.  If $3P-Q=0$ and $(P,Q) \neq (1,3)$, then
the equations are contradictory.  Hence $3P-Q \neq 0$, and the equations are uniquely solvable:
\begin{eqnarray*}
\alpha &=& \frac {3-Q}{3P-Q} \pi \\
\beta &=& \frac {P-1}{3P-Q} \pi \\
\gamma &=& \frac {2P-2} {3P-Q} \pi
\end{eqnarray*}
Since $\alpha < \beta$, if $3P-Q < 0$ we have $ Q-3< 1-P$; hence $P+Q < 4$, contradiction.
Hence $3P-Q > 0$. 
We have $P-1 > 0$, hence $3P > Q$; hence $3-Q > 0$, hence $Q < 3$.  

If $n=1$ and $m=0$ then 
\begin{eqnarray*}
2\pi &=& \alpha +  3 \gamma  \mbox{\qquad since $\ell = 3$, $n=1$, and $m=0$} \\
2(\alpha + \beta + \gamma) &=&\alpha +  3\gamma \\
\gamma &=& \alpha + 2\beta \\
\alpha + \beta + (\alpha + 2\beta) &=& \pi  \mbox{\qquad since $\alpha + \beta + \gamma = \pi$} \\
2\alpha + 3 \beta &=& \pi
\end{eqnarray*}
If $R=1$, then  by Lemma \ref{lemma:R1}, $\beta = (P-1)\alpha$, so $(2 + 3(P-1))\alpha = \pi$, and we have 
\begin{eqnarray*}
\alpha &=& \frac{\pi}{3P-1} \\
\beta &=& \frac {P-1}{3P-1} \pi \\
\gamma &=& \frac {2P-1}{3P-1} \pi
\end{eqnarray*}
If $R=0$ then we have 
\begin{eqnarray} \label{eq:twoalphaplusthreebeta}
 2 \alpha + 3 \beta &=& \pi  \mbox{\qquad as shown above} \nonumber \\
P \alpha + Q \beta &=& \pi  \mbox{\qquad since $R=0$}
\end{eqnarray}
If $P=2$ and $Q=3$ then case (iii) of the theorem holds.  Therefore we 
may assume $(P,Q) \neq (2,3)$.
  It follows that 
$2Q-3P \neq 0$, since if $2Q-3P=0$ and $(P,Q) \neq (2,3)$, subtracting the two equations (\ref{eq:twoalphaplusthreebeta}) 
yields a contradiction.  Hence the determinant of the system (\ref{eq:twoalphaplusthreebeta})  is nonzero, and we have
\begin{eqnarray*}
\alpha &=& \frac {3-Q}{3P-2Q} \pi \\
\beta &=& \frac {P-2}{3P-2Q} \pi \\
\gamma &=& \frac {2P -Q - 1} {3P-2Q} \pi
\end{eqnarray*}
Since $\alpha < \beta$, if $3P-2Q < 0$ then we have $Q-3 < 2-P$, or $P+Q<5$, contradicting $P+Q+R \ge 5$.  Hence
$3P-2Q > 0$.   Then since $\alpha > 0$ we have $Q \le 2$, and since $\beta > 0$ we  have $P \ge 3$, 
and since $\alpha < \beta$ we have $3-Q < P-2$, or $P+Q > 5$.

If $n=0$ and $m=1$ then 
\begin{eqnarray*}
2\pi &=& \beta +  3 \gamma \\
2(\alpha + \beta + \gamma) &=&\beta +  3\gamma \\
\gamma &=& 2\alpha + \beta \\
\alpha + \beta + (2\alpha + \beta) &=& \pi  \mbox{\qquad since $\alpha + \beta + \gamma = \pi$} \\
3\alpha + 2 \beta &=& \pi
\end{eqnarray*}
If $R=1$, then  by Lemma \ref{lemma:R1}, $\beta = (P-1)\alpha$, so $(3 + 2(P-1))\alpha = \pi$, and we have 
\begin{eqnarray*}
\alpha &=& \frac{\pi}{2P+1} \\
\beta &=& \frac {P-1}{2P+1} \pi \\
\gamma &=& \frac {P+1}{2P+1} \pi
\end{eqnarray*}
Since $\alpha < \beta$ we have $2P+1 < P-1$, which is impossible.  Hence this case does not occur.
If $R=0$ then we have 
\begin{eqnarray*}
 3 \alpha + 2 \beta &=& \pi \\
P \alpha + Q \beta &=& \pi
\end{eqnarray*}
If $P=3$ and $Q=2$ then  case (iii) of the theorem holds.  Therefore we 
may assume $(P,Q) \neq (3,2)$.    Therefore $3Q-2P \neq 0$ and 
\begin{eqnarray*}
\alpha &=& \frac {2-Q}{2P-3Q} \pi \\
\beta &=& \frac {P-3}{2P-3Q} \pi \\
\gamma &=& \frac {P-2Q+1} {2P-3Q} \pi
\end{eqnarray*}
Since $\alpha < \beta$, if $2P-3Q < 0$, we have $Q-2 < 3-P$, or $P+Q < 5$, contradiction.
Hence $2P-3Q \ge 0$.  Then since $\alpha > 0$ we have $Q=0$ or $Q=1$.

If $m=2$ and $n=0$ we have
\begin{eqnarray*}
2\pi &=& 3\gamma + 2\beta \\
&=& 3(\pi-\alpha-\beta) + 2\beta \\
3\alpha +\beta &=& \pi 
\end{eqnarray*}
If $R=1$ then $\beta = (P-1)\alpha$ so $3\alpha + (P-1)\alpha = \pi$, and we have
\begin{eqnarray*}
\alpha &=& \frac \pi {P+2} \\
\beta &=& \frac {P-1}{P+2} \pi \\
\gamma &=& \frac 2 { P+2} \pi
\end{eqnarray*}
Since $\beta < \gamma$ we have $P-1 < 2$, or $P < 3$; since $Q=0$ this makes $P+Q+R < 4$, contradiction; so 
this case cannot occur.

If $R=0$ then we have
\begin{eqnarray*}
 3 \alpha + \beta &=& \pi \\
P \alpha + Q \beta &=& \pi
\end{eqnarray*}
We do not have $P=3$ and $Q=1$ since $P+Q \ge 5$.  Hence
\begin{eqnarray*}
\alpha &=& \frac {1-Q}{P-3Q} \pi \\
\beta &=& \frac {P-3}{P-3Q} \pi 
\end{eqnarray*}
Since $\alpha < \beta$, if $P-3Q < 0$ we have $Q-1 < 3-P$, or $P+Q < 4$, contradiction.
Hence $P-3Q > 0$.  Since $\alpha > 0$ we have $Q = 0$.  Hence 
\begin{eqnarray*}
\alpha &=& \frac \pi P  \\
\beta &=& \frac {P-3}{P} \pi \\
\gamma &=& \frac {2} {P} \pi
\end{eqnarray*}
Since $\alpha < \beta$, we have $P > 4$.
Since $\beta < \pi/2$, we have $P < 6$.  Hence $P=5$.   That makes $\alpha = \pi/5$ and $\beta = \gamma = 2\pi/5$, 
a case which has already been ruled out. 
The following table summarizes the results obtained above for $\ell = 3$:
\bigskip

\hskip 0.8 cm
\begin{tabular}{|c|c|c|c|c|c|c|c|c|} \hline
$n$ & $m$ & $R$  & $P$        & $Q$       & $\alpha/\pi$           & $\beta/\pi$        &$\gamma/\pi$ &info              \\ \hline
2   & 0   & 1    &            &  0        & $\frac 1 { 3P-2}$      & $\frac{P-1}{3P-2}$ & $\frac{2P-2}{3P-2}$&$\gamma = 2\beta $  \\ \hline
2   & 0   & 0    &            & $Q \le 2$ & $\frac {3-Q} { 3P-Q}$  & $\frac{P-1}{3P-Q}$ &$\frac{2P-2}{3P-Q}$ &$\gamma = 2\beta $  \\ \hline
1   & 0   & 1    &            &  0        & $\frac 1 { 3P-1}$      & $\frac{P-1}{3P-1}$ &$\frac{2P-1}{3P-1}$ &$\gamma = 2\beta + \alpha$  \\ \hline
1   & 0   & 0    & $P \ge 3$  & $Q \le 2$ & $\frac {3-Q} { 3P-2Q}$ & $\frac{P-2}{3P-2Q}$& $\frac{2P-Q-1}{3P-2Q}$ & $\gamma= 2\beta+ \alpha $   \\ \hline
0   & 1   & 0    & $P \ge 6-Q$&  0 or 1   & $\frac {2-Q} { 2P-3Q}$ & $\frac{P-3}{2P-3Q}$ &$\frac{P-2Q+1}{2P-3Q}$ &  $\gamma = 2\alpha + \beta$ \\ \hline
\end{tabular}
\bigskip

We have no further need of vertex $V$ and the associated numbers $n$ and $m$, so we reprint the table without the first two columns.
The point is, that the vertex-splitting numbers $P$, $Q$, and $R$  determine $\alpha$, $\beta$, and $\gamma$ uniquely.    
\bigskip

\hskip 1.5 cm
\begin{tabular}{|c|c|c|c|c|c|c|} \hline
 $R$  & $P$        & $Q$       & $\alpha/\pi$           & $\beta/\pi$        &$\gamma/\pi$ &info              \\ \hline
 1    &            &  0        & $\frac 1 { 3P-2}$      & $\frac{P-1}{3P-2}$ & $\frac{2P-2}{3P-2}$&$\gamma = 2\beta $  \\ \hline
 0    &            & $Q \le 2$ & $\frac {3-Q} { 3P-Q}$  & $\frac{P-1}{3P-Q}$ &$\frac{2P-2}{3P-Q}$ &$\gamma = 2\beta $  \\ \hline
 1    &            &  0        & $\frac 1 { 3P-1}$      & $\frac{P-1}{3P-1}$ &$\frac{2P-1}{3P-1}$ &$\gamma = 2\beta + \alpha$  \\ \hline
 0    & $P \ge 3$  & $Q \le 2$ & $\frac {3-Q} { 3P-2Q}$ & $\frac{P-2}{3P-2Q}$& $\frac{2P-Q-1}{3P-2Q}$ & $\gamma= 2\beta+ \alpha $   \\ \hline
 0    & $P \ge 6-Q$&  0 or 1   & $\frac {2-Q} { 2P-3Q}$ & $\frac{P-3}{2P-3Q}$ &$\frac{P-2Q+1}{2P-3Q}$ &  $\gamma = 2\alpha + \beta$ \\ \hline
\end{tabular}
\bigskip

Now that we see how $P$, $Q$, and $R$ determine the angles, we will show that each of these determinations leads to a contradiction.
Let the vertices of the tiling, other than $A$, $B$, and $C$,
 be $V_1, V_2, \ldots$;  let the angle sum at $V_i$ be $k_i \pi$, so that $k_i = 1$ for 
a non-strict or boundary vertex, and 2 for a strict interior vertex.
At each vertex $V_i$, let $n_i$, $m_i$, and $\ell_i$ be the number of $\alpha$, $\beta$, and $\gamma$ angles at that vertex.
Since there are $N$ tiles altogether, each having one $\alpha$, one $\beta$,
and one $\gamma$ angle, we have $P + \sum n_i = Q + \sum m_i = R + \sum n_i = N\pi$. 
Consider the quantity $q_i = 2\ell_i + m_i$.   We will show below that, for the first three rows of the table,  at each vertex $V_i$
we have $3n_i \ge q_i$.   Once that is proved, we finish the proof as follows:  Adding over all vertices, we have
$$ \sum 3n_i \ge \sum (2\ell_i + m_i) $$
The sum on the left is three times the total number of $\alpha$ angles at the vertices $V_i$.   This must be equal to $3(N-P)$.
The sum on the right is the total number of $\beta$ angles plus twice the number of $\gamma$ angles.   This 
is equal to $N-Q + 2(N-R) = 3N-Q-2R$.   Hence 
\begin{eqnarray*}
3(N-P) &\ge& 3N-Q-2R \\
3N-Q - 2R &\le& 3(N-P) \\
 -Q-2R &\le& -3P \\
3P &\le& Q + 2R
\end{eqnarray*}
But the table above shows that in every line of the table, we have $Q + 2R \le 2$.  Hence $3P \le 2$.  Since 
$P$ is a positive integer, this is a contradiction.

It remains to supply the proof that $3n_i \ge q_i$, for the cases given in the first three rows of the table.

We first claim that in case $k_i = 1$ (i.e. the angle sum at vertex $V_i$ is $\pi$), we can assume $\ell_i \le 1$,
i.e. if $\ell_i > 1$, the conclusion $3n_i \ge q_i$ holds.
We certainly have
$\ell_i \le 2$, since $3\gamma > \pi$.  Now assume that $\ell_i = 2$.  
Then $m_i =0 $, since $2\gamma + \beta  > \pi$.  Hence $q_i = 4$.  
Then $3n_i \ge q_i$ unless $n_i = 0$ or $n_i = 1$.  If $n_i = 0$, then 
the equation $n_i\alpha + m_i \beta + \ell_i \gamma = \pi$ becomes 
$2 \gamma = \pi$,  so $\gamma = \pi/2$; but we have proved that (we can assume that)
 $\gamma$ is not a right angle, so $n_i \neq 0$.   Hence we may assume $n_i = 1$.  Then, we have
$2\gamma + \alpha = \pi$.  But this, together with $\alpha + \beta + \gamma = \pi$,
implies $\beta = \gamma$, which contradicts the assumption that $T$ is not isosceles.
This contradiction shows that $n_i \neq 1$, and hence disposes of the case $k_i =1$ and $\ell_i = 2$.
That completes that proof of the claim that we can assume $\ell_i \le 1$ when $k_i = 1$.

For the first two rows of the table, we have $\gamma = 2\beta$.  In that case
$\pi = \alpha + \beta + \gamma = \alpha + 3\beta$, so $4\beta > \pi$.
Then $2\gamma + \beta = 5 \beta > \pi$ as well.

If $k_i = 1$ and $\ell_i = 0$ then $q_i = m_i \le 3$, since $4\beta  > \pi$.
If $k_i = 1$ and $\ell_i = 1$ then $q_i = 2 + m_i$, and since $\gamma = 2\beta$ and $4 \beta > \pi$,
we have $\gamma + 2\beta > \pi$, so $m_i \le 1$.  Hence $q_i \le 3$.  Since we showed above 
that we can assume $\ell_i \le 1$ when $k_i = 1$, we have proved that (for the first two rows of the table)
$q_i \le 3$ when $k_i = 1$  (or else the conclusion of the lemma holds).

Now we assume, for proof by contradiction, that $k_i = 2$ and the first row of the table applies.  Then
\begin{eqnarray}
q_i \beta + n_i \alpha &=& (2\ell_i + m_i)\beta + n_i \alpha \nonumber \\
&=& \ell_i \gamma + m_i\beta + n_i \alpha  \qquad\mbox{since $2\beta = \gamma$} \nonumber\\
q_i \beta + n_i \alpha&=& 2\pi \label{eq:44}
\end{eqnarray}
We wish to  bound $q_i$.   First we take up the first case (first row of the table).  Then in addition to $\gamma = 2\beta$
we have $\beta = (P-1)\alpha$.   We have 
\begin{eqnarray*}
6\beta + 2\alpha &=& 2\pi \qquad \mbox {since $6\beta = 2\beta + 2\gamma$} \\
q_i\beta + n_i \alpha &=& 2\pi \qquad \mbox{by (\ref{eq:44})} 
\end{eqnarray*}
Subtracting these two equations we have $(6-q_i)\beta = (n_i-2) \alpha$.
If $q_i > 6$ then $n_i < 2$.  Then $n_i=0$ or $n_i=1$, and $\beta = ((2-n_i)/(q_i-6)) \alpha $. But $\beta = (P-1) \alpha$,
so $(2-n_i)/(q_i-6) = P-1$, an integer at least 4.  This is impossible since the numerator $2-n_i$ is at most 2.  This 
shows that in the first row of the table, we cannot have $k_i = 2$.

Now we give a bound for $q_i$ in case $k_i=2$ and the second row of the table applies, so $\gamma = 2\beta $ 
and $\beta = ((P-1)/(3-Q)) \alpha$.  Assume $q_i > 6$. 
Then as in the previous paragraph, we have
$n_i = 0$ or $n_i=1$ and $\beta = ((2-n_i)/(q_i-6))\alpha$.  But now also $\beta = ((P-1)/(3-Q)) \alpha$.
Hence 
\begin{eqnarray*}
 \frac {2-n_i}{q_i-6} &=& \frac {P-1}{3-Q} \\
\frac {(2-n_i)(3-Q)} {q_i-6} &=& P-1
\end{eqnarray*}
On the left side, the  factor $2-n_i$  is either 1 or 2, and the factor $3-Q$ is 
either 1, 2, or 3, since $Q \le 2$.  Since $R=0$  in the second row of the table, we have  $P + Q \ge 5$.
If $Q=2$ then we have $(2-n_i)/(q_i-6) = P-1$ and $ P-1 \ge 2$; but then $(2-n_i)/(q_i-6)$ is an integer at least 2,
with numerator at most 2.  This is possible only if $q_i=7$ and $n_i=0$ and $(P,Q)=(3,2)$.  But then
  case (iii) of the theorem holds.  Therefore the case 
  $q_i > 6$ has led to a contradiction or to conclusion (iii) of the theorem.  Therefore we may assume
   $q_i \le 6$ in case $k_i=2$ and the 
second row of the table applies.

Now we will prove $q_i \le 3n_i$ for the first two rows of the table.
\smallskip

Case (i) $k_i=1$. We have 
\begin{eqnarray*}
q_i &=& 2\ell_i + m_i  \qquad\mbox{definition of $q_i$} 
\end{eqnarray*}
Multiplying by $\beta$ and adding $n_i \alpha$ we have
\begin{eqnarray*}
q_i\beta + n_i\alpha &=&  2\ell_i \beta + m_i \beta + n_i \alpha\\
&=& \ell_i \gamma + m_i\beta + n_i\alpha \qquad\mbox{since $\gamma = 2\beta$} \\
&=&\pi 
\end{eqnarray*}
Now if $n_i = 0$ then $q_i \beta = \pi$.  As shown above, we can assume $q_i \le 3$.
Then, since $q_i \beta = \pi$,  we have $\beta \ge \pi/3$.  Hence $\gamma = 2\beta \ge 2\pi/3$
and $\beta + \gamma \ge \pi$, contradiction.  Hence $n_i \neq 0$.  Therefore $n_i \ge 1$.  Since $q_i \le 3$,
we have $q_i \le 3 n_i$ as desired. 

Case (ii) $k_i=2$.  We showed above that this is impossible for the first row of the table, so the 
second row applies. 
We have $3\beta + \alpha = \pi$ and $q_i\beta + n_i \alpha = 2\beta + n_i\alpha = \pi$.
Therefore
\begin{eqnarray*}
2\beta + n_i\alpha &=&  3\beta + \alpha
\end{eqnarray*}
From the second row of the table, we have $\beta = ((P-1)/(3-Q)) \alpha$.  Hence we obtain 
\begin{eqnarray*}
 2\beta + n_i \alpha &=& 2\beta + \bigg(1 + \frac {P-1}{3-Q}\bigg) \alpha \\
n_i &=& \bigg(1 + \frac {P-1}{3-Q}\bigg)
\end{eqnarray*}
Since $Q \le 2$ the denominator is positive and at most 2; since $P+Q \ge 5$, we have $P \ge 3$ so the numerator is at least 2.
It follows that $n_i \ge 2$.  Since we proved above that $q_i \le 6$ in case $k_i=2$ and the second row of the table applies,
we have $3n_i \ge q_i$ as desired.   That completes the proof for the first two rows of the table.
\smallskip

Next we take up the third row of the table.  In this case we have $\gamma = 2\beta + \alpha$ instead of $\gamma = 2\beta$.
Again we consider $q_i = 2\ell_i + m_i$ and try to prove $3n \ge q$.  We start by observing that 
$\pi = \alpha + \beta + \gamma = 2\alpha + 3\beta$.  Hence $2\gamma = 4\beta + 2\alpha > \pi$, so $\gamma > \pi/2$.
Therefore, if  $k_i=1$, we have $\ell_i \le 1$, and if $k_i = 2$, we have $\ell_i \le 3$.

We argue by (more) cases.  

Case (iii).  $k_i = 1$ and $\ell_i =1$.
Then we have $m_i \le 1$, 
since $\gamma + 2\beta = 2\beta + \alpha + 2\beta$ $ = 4\beta + \alpha > 3\beta + 2\alpha = \pi$.  Then 
\begin{eqnarray*}
\pi &=&  \gamma + m_i\beta + n_i \alpha  \\ 
3\beta + 2\alpha &=& \gamma + m_i\beta + n_i \alpha  \\ 
3\beta + 2\alpha &=& (2+m_i)\beta + (n_i+1)\alpha  \mbox{\qquad since $\gamma = 2\beta + \alpha$}\\
(1-m_i) \beta &=& (n_i-1) \alpha \\
(1-m_i)(P-1)\alpha &=& (n_i-1)\alpha \\
n_i-1 &=&(1-m_i)(P-1)
\end{eqnarray*}
If $m_i = 1$, so the right side is zero, then $n_i = 1$ also, and we have $n_i = m_i = \ell_i = 1$.  In that case
$3n_i \ge q_i$ and we are finished.  Hence (since $m_i \le 1$) we may assume $m_i = 0$.  Then we have 
$n_i = P$.
Then $q = 2\ell + m_i = 2$, and $3n_i = 3P > q_i$ since $P \ge 2$, since $Q \le 2$ and $R\le1$ and $P+Q+R \ge 5$.
That disposes of Case (iii).

Case (iv) $k_i=1$ and $\ell_i = 0$ and $m_i > 4$.  This case can be ruled out, because 
$5\beta > 3\beta + 2\alpha = \pi$.  

Case (v) $k_i=1$ and $\ell_i = 0$ and $m_i = 4$.
Then $\pi = 4\beta + n_i \alpha > 3\beta + (n_i+1)\alpha$.  Since $3\beta + 2\alpha = \pi$, this implies $n_i=0$ and
$\beta = \pi/4$.   Then $\gamma = 2\beta + \alpha = \pi/2 + \alpha$.  Hence 
\begin{eqnarray*}
\pi &=& \alpha + \beta + \gamma \\
  &=& \alpha +  \pi/4 + \pi/2 + \alpha \\
\alpha   &=&  \pi/8 \\
\gamma  &=& 5\pi/8 
\end{eqnarray*}
For this special case of $T$, we cannot establish the desired bound.  But we can rule out this case by
cosidering the possibilities
for $P$ and $Q$, as follows.    We have assumed that
row 3  of the table applies.  Then $R=1$, so $P+Q \ge 4$, and  $\alpha= \pi/(3P-1)$, and since $\alpha = \pi/8$, 
we have $P=3$, and since $P+Q \ge 4$, we have $Q \ge 1$. But in row 3 of the table, we have $Q=0$, contradiction.
That disposes of this case. 

Case (vi) $k_i=1$ and $\ell_i = 0$ and $m_i \le 3$.  Then 
$\pi = m_i\beta + n_i \alpha$; and since $\pi = 3\beta + 2\alpha$, we have $n_i \ge 2$, and so $3n_i \ge 6$,
and  $q_i = 2\ell_i + m_i  = m_i \le 3$, so $3n_i > q_i$ as required.

This disposes of all cases where $k_i = 1$ and the third row of the table applies.  
Now suppose $k_i=2$ and the third row of the table applies, so we have $\gamma = 2\beta + \alpha$
and $\pi = 2\alpha + 3 \beta$.
Then we have
\begin{eqnarray*}
2\pi &=& \ell_i\gamma + m_i\beta + n_i \alpha \\
&=& \ell_i(2\beta + \alpha) +  m_i\beta + n_i \alpha \\
&=&  q_i \beta + (\ell_i + n_i) \alpha
\end{eqnarray*}
Subtracting from this $2\pi = 4\alpha + 6\beta$ (which is twice $\pi = 2\alpha + 3\beta$), 
we obtain
\begin{eqnarray}
0 &=& (q_i-6) \beta + (\ell_i+n_i-4) \alpha \label{eq:24}
\end{eqnarray}

Case (vii)  The third row of the table applies, and $k_i = 2$, and $q_i > 6$.
Then we have $\beta = (P-1)\alpha$, so we find 
\begin{eqnarray*}
0 &=& (q_i-6)(P-1) \alpha + (\ell_i+n_i-4) \alpha  \qquad\mbox{by (\ref{eq:24})} \\
4-\ell_i - n_i &=& (P-1)(q_i-6)
\end{eqnarray*}
In the third row of the table, we have $Q=0$, and since $P+Q + R \ge 5$ and $R=1$, we have $P \ge 4$,
so $P-1 \ge 3$.  

Since $q_i > 6$, the factor $q_i -6$ is positive. The left side, on the 
other hand, is an integer at most 4.  Hence, $q_i - 6$ must be 1, so $q_i = 7$,
and both sides are equal either to 3 or to 4.  Hence, looking at the left side
$4 - \ell_i -n_i$ is 3 or 4, so $\ell_i + n_i$ is 0 or 1.  In case it is zero,
we have $\ell_i = n_i = 0$ and $P=5$.  Since $q_i = 2\ell_i + m_i = 7$, we have $m_i = 7$,
so 
\begin{eqnarray*}
\beta &=& 2\pi/7 \\
\alpha &=& \frac \beta {P-1} \\
     &=& \beta/4\\ 
     &=& \pi/14 \\
\gamma &=& 2\beta + \alpha \\
     &=& 9 \pi/14
\end{eqnarray*}
This case is impossible, by Lemma \ref{lemma:special-14}.

Now suppose $\ell_i = 1$ and $n_i = 0$.
Then $q_i = 7$ implies $m_i = 1$, so the equation $\ell_i \gamma + m_i \beta + n_i \alpha = 2\pi$
becomes $\gamma + \beta  = 2\pi$; but since $\gamma + \beta < \pi$, that is a contradiction.
Suppose $\ell_i = 0$ and $n_i = 1$.  Then $q_i = 7$ implies $m_i = 7$, so the equation 
$\ell_i \gamma + m_i\beta + n_i \alpha = 2\pi$ becomes $7 \beta + \alpha = 2\pi$.
In this case we have $P=4$, so $\beta = (P-1)\alpha = 3\alpha$.  Hence 
\begin{eqnarray*}
2\pi &=& 7 \beta + \alpha \\
&=& 21 \alpha + \alpha \\
&=& 22 \alpha \\
\alpha &=& \frac \pi {11} \\
\beta &=&  \frac {3\pi} {11}\\
\gamma &=& \frac {7\pi} {11}
\end{eqnarray*}
Since $R=1$, the largest angle of $ABC$ must be at least $\gamma$,  but according to 
 Lemma \ref{lemma:special-11}, that is not possible.  

Case (viii) The third row of the table applies, and $k_i = 2$, and $q_i = 6$.
Then by (\ref{eq:24}) we have
\begin{eqnarray*}
0 &=& (\ell_i + n_i - 4) \alpha \\
\ell_i + n_i &=& 4 
\end{eqnarray*}
Since $q_i = 2\ell_i + m_i$, we have $\ell_i = (q_i - m_i)/2$.  Hence
\begin{eqnarray*}
q_i - m_i + 2 n_i &=& 8 \\
q_i &\le& 2n_i + m_i \\
  &\le& 3 n_i
\end{eqnarray*}
which disposes of this case.

Case (ix)  The third row of the table applies, and $k_i = 2$, and  $q_i < 6$. 
By (\ref{eq:24}) we have
\begin{eqnarray*}
4-\ell_i - n_i &=& (P-1)(q_i-6)
\end{eqnarray*}
The right side is negative, since $P \ge 5$ and $q_i < 6$. Hence
the left side $4-\ell_i - n_i$ must also be negative.
Since $\ell_i \le 3$, we must have $n_i \ge 2$.  But then $3n_i \ge 6 > q_i$,
and we have disposed of this case.

That completes all cases involving the third row of the table.

Now we take up the fourth row of the table.  Here we also have $\gamma = 2\beta + \alpha$,
but the formulas for $\beta$ and $\alpha$ in terms of $P$ and $Q$ are different, and we
have $R=0$.  It is simplest just to enumerate the possibilities for $\ell_i$, $m_i$, and 
$n_i$.  We have $3\gamma + \alpha = 2\pi$  (since $\ell = 3$, $m=0$, and $n=1$,
from the first table);  that gives the first row in the following table.  The 
others are obtained by using the equation $\gamma = 2\beta + \alpha$ and the 
equation $\beta = (P-2)/(Q-3) \alpha$ to ``trade in'' some larger angles for smaller ones.
The entries containing fractions apply only if the fraction shown is an integer.
This table assumes $k_i = 2$, so the angle sum is $2\pi$.  
\bigskip

\hskip 2 cm
\begin{tabular}{|c|c|c|c|c|} \hline
$\ell_i$ & $m_i$ &  $n_i$  & $2\ell_i + m_i$& $2m_i + \ell_i$ \\ \hline
3 & 0 & 1 & 6 & 3\\ \hline
2 & 2 & 2 & 6 & 6\\ \hline
2 & 0 & $1 + \frac{2(P-2)}{3-Q}$ & 4 & 2 \\ \hline
2 & 1 & $1 + \frac{P-2}{3-Q}$ & 1 & 4 \\ \hline
1 & 4 & 3 & 4 & 9\\ \hline
1 & 4-r & $3 + r \frac{P-2}{3-Q}$ & $6-r$& $9-2r$ \\ \hline
0 & 6 & 4 & 6 & 12\\ \hline
0 & 6-r & $3 + r \frac{P-2}{3-Q}$ & $6-r$& $12-2r$ \\ \hline
\end{tabular}
\bigskip

We note that the estimate $q_i = 2\ell_i + m_i \le 3 n_i$,  which we used for the first three rows
of the table, fails here, e.g. in the first row.   But we can replace it by $r_i = 2m_i + \ell_i \le 3 n_i$, shown in the 
last column of the table.   

Remembering that
the rows with fractions occur only if the fraction is an integer, and that $P-2$ and $Q-3$ are positive integers
since $Q \le 2$ and $P + Q \ge 5$,  inspection of the table shows that $r_i \le 3 n_i$ in each row.
This table applies only when $k_i = 2$, so we also need to consider the case $k_i = 1$.  Since $3\gamma + \alpha = 2\pi$,
we have $\gamma > \pi/2$, so when $k_i = 1$, we have $\ell_i = 0$ or 1.  If $\ell_i = 1$,  then since $\alpha \neq \beta$,
we have $m_i = n_i = 1$ (in which case $r_i = 3 \le 3 n_i$),  or $m_i = 0$ and $\beta$ is an integer multiple of $\alpha$.
In that case, $n_i$ is at least 3, since $\beta > \alpha$, and $r_i = 1$, so $r_i \le 3 n_i$.   Now suppose 
$k_i = 1$ and $\ell_i = 0$.  Since $\pi = (2\beta + \alpha) + \beta + \alpha = 3\beta + 2 \alpha$,
we have $m_i = 3$ and $n_i = 2$ as one possibility, in which case $r_i = 6 \le 3 n_i$;  and other possibilities
may arise if $3\beta$ is a multiple of $\alpha$, but if that happens, then $r_i < 6$ and $n_i > 2$, so we still
have $r_i \le 3 n_i$.

Hence, for the fourth row of Table 2, we always have $r_i \le 3 n_i$.
Now we sum over all vertices of the tiling that are not vertices of $ABC$.  We have
\begin{eqnarray*}
\sum r_i &=& \sum (2m_i + \ell_i) \\
&=& 2 \sum m_i + \sum \ell_i \\
&=& 2(N-Q) + N  \qquad \mbox{since $R=0$}\\
&=& 3N-2Q
\end{eqnarray*}
On the other hand we have 
\begin{eqnarray*}
3\sum n_i &=& 3(N-P) \\
\end{eqnarray*}
Since $\sum r_i \le 3 \sum n_i$, we have $3N-2Q \le 3N-3P$.  Hence $3P \le 2Q$.  But $Q \le 2$ and $P \ge 3$, so 
$3P \ge 9 > 4 \ge 2Q$, contradiction.   That completes the proof that the fourth row of Table 2 is impossible.

Now we turn to the fifth (and last) row of Table 2.   Here we have $\gamma = 2\alpha + \beta$.   The following 
table shows the possible vertices.   The rows containing fractions represent possible vertices only 
in case the fraction shown is actually an integer.  Since $Q \le 1$ in the fifth row of table 2, the 
denominator of the fractions shown is either 1 or 2.
\bigskip

\hskip 2 cm
\begin{tabular}{|c|c|c|c|} \hline
$k_i$ &$\ell_i$ & $m_i$ &  $n_i$   \\ \hline
2 & 3 & 1 & 0 \\ \hline
2 & 3 & 0 & $ \frac{P-3}{2-Q}$  \\ \hline
2 & 2 & 2 & 2  \\ \hline
2 & 2 & 1 & $2 + \frac{P-3}{2-Q}$   \\ \hline
2 & 2 & 0 & $2+2\frac{P-3}{2-Q}$   \\ \hline
2 & 1 & 3 & 6   \\ \hline
2 & 1 & 3-r & $6 + r \frac{P-3}{2-Q}$   \\ \hline
2 & 0 & 4 & 8 \\ \hline
2 & 0 & 4-r & $8 + r \frac{P-3}{2-Q}$   \\ \hline
1 & 1 & 1 & 1 \\ \hline
1 & 0 & 0 & $\frac{P-3}{2-Q}$  \\ \hline
1 & 0 & 1 & 2 \\ \hline
1 & 0 & 1 & $\frac{P-3}{2-Q}$  \\ \hline
1 & 0 & 0 & $2\frac{P-3}{2-Q}$  \\ \hline
\end{tabular}
\bigskip

For this row of Table 2, the estimates used in the previous rows do not work, because of
 the possibility shown in the first row of this table, namely 
$\ell_i = 3$ while $n_i = 0$.    Therefore we give a different argument.  Note that 
the first four rows of the table correspond to vertices
with at least two $\gamma$ angles.  Since there are altogether $N$ $\gamma$ angles,  there are at most $N/2$
such vertices.  Hence (and since $R=0$)  there are at least $N/2$ vertices other than those of $ABC$, each
contributing at most one $\gamma$ angle.  Now consider the difference $m_i - n_i$ at the two kinds of 
vertex.   At the vertices with two or three $\gamma$ angles, $m_i - n_i \le 1$; indeed only the first row
has $m_i - n_i > 0$, and in that row it is 1.   In all the other rows of the table, we have $n_i - m_i \ge 1$.
Hence, adding $n_i - m_i$ over all the vertices, we obtain at least one from each of at least $N/2$ vertices,
and at least $-1$ from the remaining (smaller number of) vertices,  so the sum of $n_i - m_i$ over all the 
vertices (not counting the vertices of $ABC$)  is non-negative.   But when we add in the contributions of the 
vertices of $ABC$, we add $P-Q$, which is a positive number (since $Q \le 1$ and $P+Q \ge 5)$.
   That is a contradiction, since altogether 
there are $N$ each of $\alpha$ angles and $\beta$ angles,  so the sum over all vertices, including those of $ABC$,
of $n_i - m_i$ must be zero.   That completes the proof that the fifth row of Table 2 is impossible, and 
it also completes the proof that $\ell = 3$ is impossible.  

We still have the cases $\ell = 4$ and $\ell = 5$ to deal with.   We now assume $\ell = 4$.  That is, 
there is at least one vertex $V$ where four $\gamma$ angles meet, and there are at most 3 other angles,
i.e. $n + m \le 3$.   
  Since $\gamma > \pi/3$,  we have $k=2$, i.e. 
the angle sum at $V$ is $2\pi$, not $\pi$.   If $n+m=0$, then $\gamma = \pi/2$, contradiction.  
If $n=0$ then we have 
\begin{eqnarray*} 
2\pi &=& 4\gamma + m \beta \\
2(\alpha + \beta + \gamma) &=& 4\gamma + m \beta \\
2\alpha + (2-m) \beta &=& 2\gamma
\end{eqnarray*}
Since $\alpha < \gamma$,  we must have $2-m > 0$, or $m < 2$.  Since $n=m=0$ is impossible, that leaves 
$m=1$ as the only possibility when $n=0$, and in that case $2\gamma = 2\alpha + \beta$.

We will prove that in this case ($\ell = 4$ and $n=0$) we have $R=0$.  By Lemma \ref{lemma:R2},
we have $R \le 1$.  Assume, for proof by contradiction, that $R=1$. Then
by Lemma \ref{lemma:R1} we have $Q=0$ and $\beta = (P-1)\alpha$, so 
$2\gamma = 2\alpha + \beta = (P+1)\alpha$.  We then have
\begin{eqnarray*}
\beta &<& \gamma \\
(P-1) \alpha &<& \frac {P+1} 2 \alpha \\
P-1 &<& \frac {P+1}2 \\
P &<& 3
\end{eqnarray*}
contradicting $P+Q\ge 5$, since $Q=0$.   That completes the proof that $R=0$.

Having proved that $R=0$, we then have 
\begin{eqnarray*}
\pi &=& \alpha + \beta + \gamma \\
&=& \alpha + \beta + \frac 1 2 (2\alpha + \beta)  \qquad\mbox{ since $2\gamma = 2\alpha + \beta$} \\
&=& 2\alpha + \frac 3 2 \beta \\
&=& P\alpha + Q\beta \\
\end{eqnarray*}
This gives us two equations in $\alpha$ and $\beta$.   The determinant of the 
system is (a multiple of) $4Q-3P$.   Assume, for proof by contradiction, that 
this determinant is zero.  Then 
\begin{eqnarray*}
\pi &=& P\alpha + Q \beta \\
4 \pi &=& 4 P \alpha + 4Q \beta \\
  &=& 4P \alpha + 3P \beta \qquad\mbox{since $4Q = 3P$}\\
&=& (4 \alpha + 3 \beta) P \\
&=& 2\pi P \qquad\mbox{since $4\alpha + 3\beta = 2\pi$}
\end{eqnarray*}
Thus $4\pi = 2\pi P$.  Hence $P = 2$.  Now we have 
\begin{eqnarray*}
2\pi &=& 4\alpha + 3 \beta \\
\pi &=& P \alpha + Q \beta \\
&=& 2 \alpha + Q \beta \qquad \mbox{since $P=2$}
\end{eqnarray*}
Multiplying the last equation by 2 and subtracting 
the first equation, we find $2Q = 3$, which is impossible since $Q$ is an integer.
This contradiction shows that the determinant is not zero.  Hence the equations are
solvable:
\begin{eqnarray*}
\alpha &=& \frac {2Q-3}{4Q-3P} \pi \\
\beta &=& \frac{4-2P}{4Q-3P} \pi  \\
\gamma &=& \alpha + \frac 1 2 \beta \\
       &=& \frac {2Q-P-1}{4Q-3P} \pi
\end{eqnarray*}
If $4Q > 3P$ then, since $\alpha < \beta$, we have 
\begin{eqnarray*}
2Q-3 &<& 4-2P \\
2(P+Q) &<& 7 
\end{eqnarray*}
contradicting $P+Q \ge 5$.  Hence $3P < 4Q$.   Writing the equations for 
$\alpha$, $\beta$, and $\gamma$ with positive denominator, we have 
\begin{eqnarray*}
\alpha &=& \frac {3-2Q}{3P-4Q} \pi \\
\beta &=& \frac{2P-4}{3P-4Q} \pi \\
\gamma &=& \frac {P+1-2Q}{3P-4Q} \pi
\end{eqnarray*}
Then the numerators must also be positive, so $2Q < 3$, which implies $Q=0$ or $Q=1$.
Assume, for proof by contradiction, that $Q=0$. Then $\alpha = \pi/P$ and $\beta = 2\pi/3 - 4\pi/(3P)$.  Since $\beta < \gamma$ we have
\begin{eqnarray*}
\beta &<& \gamma \\
\frac 2 3 \pi - \frac 4 {3P} \pi &<& \frac {P+1}{3P} \pi\\
\frac 2 3  &<& \frac {P+5} {3P} \\
6P &<& 3(P+5) \\
3P &<&  15\\
 P &<& 5 
\end{eqnarray*}
contradicting $P+Q \ge 5$.  Hence $Q \neq 0$.  Hence $Q=1$.   Then 
\begin{eqnarray*}
\alpha &=& \frac {1}{3P-4} \pi \\
\beta &=& \frac{2P-4}{3P-4} \pi \\
\gamma &=& \frac {P-1}{3P-4} \pi
\end{eqnarray*}
Since $\beta < \gamma$, we have $2P-4 < P-1$, which implies $P < 5$, contradicting $P + Q \ge 5$.  
Hence every subcase of $\ell = 4$, $n=0$ has led to a contradiction.  Hence $\ell = 4$ and $n=0$ cannot occur.

Next we take up the case $\ell =4$ and $n=1$. Then we have 
\begin{eqnarray*} 
2\pi &=& 4\gamma + m \beta + \alpha \\
2(\alpha + \beta + \gamma) &=& 4\gamma + m \beta  + \alpha\\
\alpha + (2-m) \beta &=& 2\gamma
\end{eqnarray*}
and since $\alpha < 2\gamma$, we must have $2-m > 0$, which implies $m < 2$.  The case $m=1$ can also
be ruled out, for then we would have $\alpha + \beta = 2 \gamma$, and 
adding $\gamma$ to both sides, we would have $\pi = 3 \gamma$,  contradicting $\gamma > \pi/3$. 
Hence $m=0$ when $n=1$, and we have $2 \gamma=\alpha + 2 \beta$.   Adding $\gamma$ to both sides we have
\begin{eqnarray*}
3 \gamma &=& \alpha + 2\beta + \gamma \\
&=& \pi + \beta \\
\gamma &=& \frac \pi 3 + \frac \beta 3 \\
\gamma &<& \frac \pi 2 \qquad \mbox{since $4 \gamma + \alpha = 2\pi$} \\
\end{eqnarray*}
Therefore
\begin{eqnarray*}
\beta + \frac \pi 3  &<& \frac \pi 2 
\end{eqnarray*}
Subtracting $\pi/3$ from both sides we have
\begin{eqnarray*}
\beta &<& \frac \pi 6 \\
\alpha &=& \pi - \gamma -\beta \\
   &>& \pi - \frac \pi 2 - \frac \pi 6 = \frac \pi 3  \qquad \mbox{since $\gamma < \pi/2$ and $\beta < \pi/6$}\\
\end{eqnarray*}
Therefore $\alpha > \pi/3$,  contradiction.  Hence the case $\ell = 4$ and $n=1$ is not possible.

If $\ell = 4$ and  $n=2$ then since $n+m \le 3$, we have $m=0$ or $m=1$ and 
\begin{eqnarray*} 
2\pi &=& 4\gamma + m \beta + 2\alpha \\
2(\alpha + \beta + \gamma) &=& 4\gamma + m \beta  + 2\alpha\\
(2-m) \beta &=& 2\gamma  \\
\beta &=& \frac 2 { 2-m} \gamma 
\end{eqnarray*}
Since $m$ is 0 or 1, the fraction on the right is either 2 or 1, contradicting $\beta < \gamma$.
Hence the case $n=2$ cannot occur.

If $\ell = 4$ and $n=3$ then since $n+m \le 3$ we have $m=0$, and 
\begin{eqnarray*} 
2\pi &=& 4\gamma + 3\alpha \\
2(\alpha + \beta + \gamma) &=& 4\gamma + 3\alpha\\
2\beta &=& 2\gamma + 3 \alpha  
\end{eqnarray*}
contradicting $\beta < \gamma$.  So the case $n=3$ cannot occur either.  Since $m+n < \ell$,
this is the last possibility with $\ell = 4$; hence $\ell = 4$ has been ruled out.

The last possibility for $\ell$ is $\ell = 5$.  Assume $\ell = 5$. 
 That is, 
there is at least one vertex $V$ where five $\gamma$ angles meet, and there are at most 4 other angles,
i.e. $n + m \le 4$.   
  Since $\gamma > \pi/3$,  we have $k=2$, i.e. 
the angle sum at $V$ is $2\pi$, not $\pi$.   If $n+m=0$, then $\gamma = 2\pi/5$, which has 
already been ruled out by Lemma \ref{lemma:twopioverfive}.
Assume, for proof by contradiction, that $n=0$.  Then we have 
\begin{eqnarray*} 
2\pi &=& 5\gamma + m \beta \\
2(\alpha + \beta + \gamma) &=& 5\gamma + m \beta \\
2\alpha + (2-m) \beta &=& 3\gamma \\
(2-m) \beta &=& 3\gamma - 2\alpha \\
&=& 2(\gamma -\alpha) + \gamma \\
&>& \gamma \qquad\mbox{since $\gamma > \alpha$}
\end{eqnarray*}
Since $\beta < \gamma$, this implies $2-m > 1$, which implies $m = 0$.  Since
we have assumed $n=0$, and we have proved $n+m \neq 0$,  that is a contradiction.
This contradiction completes the proof that $n \neq 0$.

Now assume $\ell = 5$ and $n=1$.  Then we have 
\begin{eqnarray*} 
2\pi &=& 5\gamma + m \beta + \alpha \\
2(\alpha + \beta + \gamma) &=& 5\gamma + m \beta  + \alpha\\
\alpha + (2-m) \beta &=& 3\gamma \\
(2-m) \beta &=& 3\gamma -\alpha \\
(2-m) \beta&>& 2\gamma \qquad\mbox{since $\alpha < \gamma$}
\end{eqnarray*}
Since $\beta < \gamma$, this implies $2-m > 2$, which is 
not possible, since $m \ge 0$.  Hence the case $\ell = 5$ and $n=1$ is impossible.

Now assume $\ell = 5$ and $n\ge 2$.  Then we have 
\begin{eqnarray*} 
2\pi &=& 5\gamma + m \beta + n\alpha \\
2(\alpha + \beta + \gamma) &=& 5\gamma + m \beta  + n\alpha\\
(2-m) \beta &=& 3\gamma + (n-2) \alpha \\
(2-m) \beta &\ge & 3\gamma \qquad\mbox{since $n\ge 2$}\\
\end{eqnarray*}
Since $\beta < \gamma$, this implies $2-m > 3$, which is 
not possible, since $m \ge 0$.  Hence the case $\ell = 5$ is impossible.
That was the last possibility for $\ell$, so the proof of the lemma is 
complete.

\section{Conclusions }
The aim of this series of papers is to classify the triples $(ABC, T, N)$ such that triangle 
$ABC$ can be $N$-tiled by $T$.  We have completed this classification, except for the case 
when the tile $T$  has a $120^\circ$ angle and is isosceles.    For example,  aside 
from that shape of tile, given an integer $N$, 
if there is any $N$-tiling, then $N$ is either a square, or a sum of two squares, or is 2,3, or 6 times
a square, or twice a sum of two squares.   For example, there are no $N$-tilings for $N= 7$, 11, or 19.   The following theorem gives more information about 
the possibilities for the shapes of the tile and the tiled triangle.  

\begin{theorem}[Main Theorem]
Suppose triangle $ABC$ is tiled by triangle $T$.  Suppose $ABC$ is not similar to $T$ and $T$ is not a 
right triangle.  Then one of the following conclusions holds: 
\smallskip

\noindent
(i)  $ABC$ is equilateral,  $T$ is isosceles with base angles $\pi/6$, and $N$ is three times a square, or
\smallskip

\noindent
(ii) $3\alpha + 2\beta = \pi$, where $\alpha$ and $\beta$ are the two smallest angles of the tile, in either order,
and $\alpha$ is not a rational multiple of $\pi$, and $\sin(\alpha/2)$ is rational (which implies that the sides of the 
tile have rational ratios), and two of the angles of $ABC$ are $2\alpha$ and $\beta$.
\smallskip 

\noindent
(iii) $\gamma = 2\pi/3$ (120 degrees) and $\alpha \neq \beta$ (the tile is not isosceles).  
\end{theorem}

Theorems covering the case when $ABC$ is similar to $T$ and the case when $T$ is a right triangle 
are in \cite{beeson1}.   Theorems covering the case in conclusion (ii)  are in \cite{beeson-triquadratics}, where 
a necessary and sufficient condition on $N$ is given for an $N$-tiling to exist in that case.
We conjecture there are no tilings falling under case (iii);  some partial results in that 
direction are in \cite{beeson120}.

\smallskip

%
%

\bibliography{TriangleTiling}
\bibliographystyle{plain-annote}
 
\end{document}